\documentclass[preprint,11pt,review]{elsarticle}
\usepackage{amssymb,amsmath,array}
\DeclareMathOperator*{\sumchannels}{sum-channels}
\DeclareMathOperator*{\conca2}{concat2}
\DeclareMathOperator*{\concat}{concat}
\journal{\ldots}

\usepackage[inline]{enumitem}
\usepackage{bm}
\usepackage{tikz}
\usetikzlibrary{decorations.markings}
\usetikzlibrary{positioning, automata}
\newcommand*{\Scale}[2][4]{\scalebox{#1}{\ensuremath{#2}}}%

\usepackage[margin=2cm]{geometry}
\usepackage{algorithm}
\usepackage{algorithmicx}
\usepackage{algpseudocode}
\algrenewcommand\algorithmiccomment[1]{!!~{\itshape #1}}
\algdef{SE}[SUBALG]{Indent}{EndIndent}{}{\algorithmicend\ }%
\algtext*{Indent}
\algtext*{EndIndent}
\usepackage{booktabs}
\usepackage{multirow}
\newcommand{\ra}[1]{\renewcommand{\arraystretch}{#1}}

\usepackage{setspace}

\usepackage{xcolor}
\usepackage{textcomp} 
\usepackage{siunitx}

\begin{document}

\begin{frontmatter}


%\ead{email address}

%\title{Title 1: The Space Filling Curve Convolutional Neural Network for use with data on Unstructured Meshes}
%\title{Title 2: Applying Convolutional Neural Networks to Data on Unstructured Meshes with Space-Filling Curves}
%\title{Title 3: Using Space-Filling Curves to Apply  Convolutional Neural Networks Directly to Data on Unstructured Meshes}
\title{Applying Convolutional Neural Networks to Data on Unstructured Meshes with Space-Filling Curves}

\author[IC]{C. E. Heaney\corref{cor1}} 
\ead{c.heaney@imperial.ac.uk}
\author[IC]{Y. Li} %liyuling2319@gmail.com
\author[Chem]{O. K. Matar} 
\author[IC,DA]{C. C. Pain} 

\cortext[cor1]{Corresponding author}
%\fntext[label2]{}

\address[IC]{Applied Modelling and Computation Group, Department of Earth Science and Engineering, Imperial College London, SW7 2AZ, UK}
\address[Chem]{Department of Chemical Engineering, Imperial College London, SW7 2AZ, UK}
\address[DA]{Data Assimilation Laboratory, Data Science Institute, Imperial College London, SW7 2AZ, UK}

\begin{abstract}
This paper presents the first classical Convolutional Neural Network (CNN) that can be applied directly to data from unstructured finite element meshes or control volume grids. CNNs have been hugely influential in the areas of image classification and image compression, both of which typically deal with data on structured grids. Unstructured meshes are frequently used to solve partial differential equations and are particularly suitable for problems that require the mesh to conform to complex geometries or for problems that require variable mesh resolution. Central to our approach are space-filling curves, which traverse the nodes or cells of a mesh tracing out a path that is as short as possible (in terms of numbers of edges) and that visits each node or cell exactly once. The space-filling curves (SFCs) are used to find an ordering of the nodes or cells that can transform multi-dimensional solutions on unstructured meshes into a one-dimensional (1D) representation, to which 1D convolutional layers can then be applied. Although developed in two dimensions, the approach is applicable to higher dimensional problems. 

To demonstrate the approach, the network we choose is a convolutional autoencoder (CAE), although other types of CNN could be used. The approach is tested by applying CAEs to data sets that have been reordered with a space-filling curve. Sparse layers are used at the input and output of the autoencoder, and the use of multiple SFCs is explored. We compare the accuracy of the SFC-based CAE with that of a classical CAE applied to two idealised problems on structured meshes, and then apply the approach to solutions of flow past a cylinder obtained using the finite-element method and an unstructured mesh.\\
\end{abstract}

%\includegraphics[width=\textwidth]{figures/Autoencoder.jpg}

%\item Use of space-filling curves to apply convolutional neural networks to data on unstructured meshes. 
%\item Use of multiple space-filling curves in order to improve the accuracy of the new convolutional neural networks. 
%\item Introduction of sparse smoothing layers at the input and output of the neural networks to reduce noise. 
%\item Incorporation of edge weights into a recurrent neural network based on mean-field theory to generate multiple space-filling curves.
%\item Application to fluid flow data on an unstructured mesh. 


\end{frontmatter}

%\linenumbers

\section{Introduction} 
% facial recognition, automated vehicles, driver-less vehicles, drones, augmented reality software, retail, medicine, predictive maintenance...
Also known as convnets, Convolutional Neural Networks (CNNs) have provided revolutionary advances in the fields of image compression~\cite{Minnen_2018,Mentzer_2020} and image classification~\cite{Krizhevsky_2012,Yim2015}. Although most often encountered in computer vision, technology from CNNs has recently shown promise in helping to solve governing equations for fluid motion~\cite{Gonzalez_2018}, as well as analysing these motions~\cite{lee2020model}. The accuracy of CNNs for image compression seems to carry over to these applications, however, CNN technology is generally based on structured %$n\times m$ grids or for %some applications $n\times %m\times k$ 3D grids
grids and is, therefore, not readily applicable to data held on unstructured meshes. Many computational fluid dynamics (CFD) problems rely on unstructured meshes to resolve complex geometries~\cite{Walton_2017} and  interfaces~\cite{Xie_2016}. These meshes typically have more resolution where needed, that is, where finer-scale features occur, such as in the vicinity of a viscous boundary layer, at the interface between fluids, at a density front or shock wave~\cite{Kampitsis_2020,Alauzet_2016}. 

Without a means of applying CNN techniques directly to data on unstructured meshes, a number of other methods have been developed. One may apply a 1D CNN directly to the vectorised nodal values of the solution variables. If the ordering of these variables has been optimised, to increase the efficiency of iterative or direct solvers for example, this might be more successful, however, we have found that this approach either does not work or has limited success. To improve this method, one could interpolate the data from an unstructured mesh to a structured mesh and then apply a CNN, although this will introduce errors due to interpolation and generally requires a much larger number of grid points to capture the features in the original data~\cite{Mack2020}. Of course, one could also restrict oneself to using structured meshes for the solution data and then apply CNNs directly to this~\cite{Gonzalez_2018,Carlberg2020,Xu2020}. However, for modelling certain geometries or physics, this is not feasible.

Graph-based methods can also be used to apply CNNs to arbitrary data that has been arranged in the form of a graph. Graph Neural Networks (GNNs) are networks that can be applied to such data~\cite{Wu2020}, for example, the data associated with an unstructured finite element mesh that has been transformed to a graph. % %\textcolor{red}{There is often structural data held within graphs which can be exploited by applying Graph Convolutional  Networks~\cite{kipf2017}.}%
As in classical CNNs, convolutional GNNs~\cite{Wu2020} or Graph Convolutional Networks~\cite{kipf2017} take a weighted average of the values at the vertices (or pixels or cells) adjacent to a given vertex and use this to form a new layer or feature map similar to classical CNNs. Recently developed within PyTorch~\cite{PyTorch_2019}, MeshCNN~\cite{Hanocka2019} is one example of this type of method. MeshCNN forms a graph from the discretisation and then applies convolutional filters directly to the data held on the unstructured mesh by using mesh coarsening operations, akin to those used in multi-grid methods~\cite{Dargaville_2020}. The resulting coarse meshes are then used to form pooling operations similar to those used in classical CNNs. MeshCNN~\cite{Hanocka2019} has a higher accuracy than other methods when tested on a variety of classification and segmentation problems. Another type of graph-based methods are Point Cloud methods~\cite{LU202042}, which form a graph from the nearest neighbouring points to a given point in the data cloud. Once the graph is constructed, convolutional filters can be applied as described above. Such point cloud methods may be very useful for working with data represented on a series of different meshes, for example, arising from the use of mesh adaptivity~\cite{Pain_2001}, as there may be no need to re-train the neural network. Graph networks (a type of GNN), which learn interactions between variables in a compressed space, have successfully been applied to problems in fluid and solid dynamics~\cite{Sanchez-Gonzalez2020,Pfaff2020}. Both approaches use an encoder and decoder to compress and reconstruct the data, and a time-stepping method for predicting the future solution variables. The former uses particle-based methods to solve the governing equations, and exploits a nearest-neighbour technique to form an underlying graph of interactions. The latter applies the finite-element method to solve the governing equations and simply uses the underlying finite-element discretisation as the graph, combined with a neighbourhood graph for modelling solid-cloth dynamics. In both cases, the graph is used to learn the dynamics of the systems.

By using Space-Filling Curves (SFCs), we demonstrate in this paper that one can apply the traditional CNN approaches (in 1D) \emph{directly} to data held on unstructured multi-dimensional meshes. The underlying ideas of space-filling curves are applied to produce a continuous curve through a finite element mesh, in which the curve traverses every node of the mesh. This space-filling curve is used to transform multi-dimensional data into 1D data, to which convolutional layers are then applied. Sparse layers are introduced to the network, which add some smoothing to the results. The use of multiple SFCs is also explored in this paper.

The first eponymous space-filling curve was discovered by Peano in 1890~\cite{Peano_1890}, when seeking a continuous mapping from the unit interval onto the unit square. One year later, Hilbert introduced a variant of this curve~\cite{Hilbert_1891} known as a Hilbert curve or Hilbert space-filling curve. Since then, many other SFCs have been developed including those by Moore, Lebesgue, Sierpi{\'n}ski~\cite{Sagan,Bader}. Ideal for transforming and reordering multi-dimensional data to 1D data, space-filling curves:
%label=\textbf{S.\arabic*}
\begin{enumerate}[label=\textbf{(P\arabic*)},ref=\textbf{P\arabic*}]
    \item generate a continuous curve in which neighbours of any node within that curve are nodes that are close to one another\label{prop1};
    \item have in-built hierarchies of features, for example, if one takes every $N$~points along the curve, the result is a coarser mesh across the whole domain\label{prop2};
    \item are Lebesgue-measure preserving (at least the classical SFCs have this property), so that sub-curves of equal length fill multi-dimensional regions of equal area or volume\label{prop3}.
\end{enumerate}
The properties of space-filling curves have already been exploited in a number of techniques used in the numerical solution of partial differential equations. For example, property~\ref{prop1} has led to the use of SFC orderings to manage cache memory when calculating matrix-matrix or matrix-vector products~\cite{Bohm2018}. Properties~\ref{prop1} and~\ref{prop2} have led to the use of SFCs to optimise the node numbering for the benefit of iterative solvers~\cite{Bungartz_2006}. Properties~\ref{prop2} and~\ref{prop3} have resulted in the extensive use of SFCs with unstructured meshes to find a partition suitable for parallel computing~\cite{Behrens_2000}. That is, if one slices the SFC into $I$ equal intervals, then the result is a domain decomposition that has $I$ partitions with minimal connectivity (few common edges) between the partitions.

Space-filling curves are a natural partner to CNNs because of properties~\ref{prop1}, \ref{prop2} and~\ref{prop3}, and have been applied to structured grids to reduce 3D structured grid data to 1D in image compression~\cite{Sprecher_2002,Papamarkos_1999}. They have also been applied to DNA sequencing and classification~\cite{Anjum_2019,Yin_2018}. The latter compares different mappings of 2D structured mesh data to 1D, to which a CNN is applied, and found that the mappings based on Hilbert SFCs performed the best. Analysis of SFCs for their data preserving abilities is investigated in \cite{Skubalska_1997}. Applying 1D CNNs to 2D data that has been transformed by an SFC has also displayed comparable accuracy and speed to classical 2D and 3D CNNs~\cite{Corcoran_2018}, again, for data on structured grids.

We present here a method of applying CNNs directly to data on unstructured meshes, such as solutions arising from finite element or control volume discretisations. The method is timely as researchers are currently looking to improve compression or dimensionality reduction methods within reduced-order model frameworks~\cite{Gonzalez_2018,Mack2020,Carlberg2020,Xu2020,Phillips_2020} by exploiting some of the useful properties of CNNs, such as the multiscale resolution, rotational and position invariance~\cite{Hinton_2011,Krizhevsky_2012,Gonzalez_2018}. Reduced-order models often rely on Singular Value Decomposition (SVD) for compression (or dimension reduction) through Proper Orthogonal Decomposition (POD)~\cite{Hesthaven_2018,Xiao_BE_2019,Xiao_CAF_2019,Ahmed2019}. SVD-based methods can be limited in their abilities to interpolate solutions, for example, suffering from Gibbs oscillations when there are abruptly changing fields~\cite{Carlberg2020} and poor accuracy for convection-dominated problems~\cite{Ahmed2019}. Going some way to tackle this, a combination of SVD and a fully-connected autoencoder is presented as a tool for dimension reduction in~\cite{Phillips_2020}, however, some of the issues of the SVD are still present. Thus, being able to apply CNNs to data on unstructured meshes could significantly improve the quality of such reduced-order models.  %Furthermore
In addition to the timely nature of this work, the elegant SFC-based approach outlined here, is expected to be computationally faster than the graph-based approaches (MeshCNN and GNNs), since the latter will be dominated by indirect addressing. In the method outlined here, there is only indirect addressing at the inputs and outputs of the space-filling curve CNN and then, only when multiple space-filling curves are used. 
%\textcolor{red}{For both the method outlined here and the approach of MeshCNN, the filters will be dependent on the connectivity and resolution, which is variable for unstructured meshes. More filters may be required to take account of this.}
%\textcolor{red}{There are also issues with the filters themselves used in the MeshCNN/GNN approach, as they will be dependent on the connectivity and resolution, which will be variable. This is the reason why the current MeshCNN approach uses the max-pooling operation to form coarsened meshes. This is also true of the unstructured mesh SFC-CNN approach. That is, the SFC will pass through areas of differing resolution which means that the filters may have to work harder (there may have to be more of them) with differing filters to take this into account. }

The remainder of this paper is as follows. The following section describes how space-filling curves are formed for an unstructured mesh. The architectures of the convolutional autoencoders used in this paper are described in Section~~\ref{sec:architectures}. The results are presented in Section~\ref{sec:results} which includes the application of SFC-based autoencoders to data on a structured mesh (for advection of a square wave and advection of a Gaussian function) and data on an unstructured mesh (for flow past a cylinder). %Section~\ref{sec:sfc_cnn_structured_applications} presents the application of the methods to structured mesh applications and Section~\ref{sec:sfc_cnn_unstructured_applications} presents results obtained by applying the methods to an unstructured mesh application; flow past a cylinder. 
Future work is then discussed and conclusions are drawn.  

\section{Determining Space-Filling Curves on Unstructured Meshes} \label{sec:determining_sfc_unstructured_meshes}
This section describes how to generate space-filling curves for unstructured meshes. The method involves ordering or numbering a nested sequence of partitions in such a way that, when following the numbering, the path through these partitions is continuous, and neighbouring vertices on the path will be close in the original graph. This is desirable, as, when applying CNNs to data, the filters search for structurally coherent features. To detect local features, the node ordering should be such that neighbouring nodes on the SFC are close to one another in physical space, so that any structural coherence in the data can be detected. The generation of space-filling curves for structured and unstructured meshes is described in Section~\ref{sec:SFC_numbering}. The method which creates the partitions used to generate the space-filling curve ordering is discussed in Section~\ref{sec:SFC_partitioning}. % %\textcolor{red}{Here we choose a recurrent neural network based on mean-field theory~\cite{Pain_1999}, which is able to generate partitions that have equal numbers of vertices (as far as possible) and that can use variable edge weights to influence the path, useful when generating multiple space-filling curves. Other methods could be used such as those based on agglomeration~\cite{Dargaville_2020}.} %
Finally, in Section~\ref{sec:SFC_multiple}, comments are made on the efficacy of multiple space-filling curves and how these are generated. 

\subsection{Forming the numbering of the space-filling curve} \label{sec:SFC_numbering}

%The classical Hilbert curve in 2D on a structured mesh is shown in figure~\ref{2d-3d-space_filling} for a $32\times 32$ grid. Two Hilbert curves are shown - one a 90 degree rotation of the other - and they are both used to form the SFC-CNN for structured meshes here. The approach for forming these curves is to form the fundamental shape - a U shape in the case of the Hilbert curve - and put this on the coarsest mesh then refine and place similar shapes where the cells/nodes are and rotate these shapes until they link with one another. Then repeat with a finer representation and so on until one visits each cell/node/vertex  of the structured mesh once. 

\begin{figure}[htpb]
\centering
\begin{tikzpicture}[scale=0.72]%[scale=0.9]
\usetikzlibrary{positioning}

% level 1 -------------------------------------------------------------------------------
\node[] at (2,-0.5) {\small (a)};
\draw[step=2cm,gray,thin] (0,0) grid (4,4);
%\draw[step=2cm,darkgray,thin] (0,0) grid (4,4);

\filldraw[red] (1,3) circle (2pt) node[anchor=east] {1};
\filldraw[red] (1,1) circle (2pt) node[anchor=east] {2};
\filldraw[red] (3,1) circle (2pt) node[anchor=west] {3};
\filldraw[red] (3,3) circle (2pt) node[anchor=west] {4};
\draw[red, ultra thick] (1,3) -- (1,1) -- (3,1) -- (3,3);

% level 1 and 2--------------------------------------------------------------------------
\node[] at (7,-0.5) {\small (b)};
\draw[step=1cm,gray,thin] (5,0) grid (9,4);
\draw[step=2cm,darkgray,thin] (5,0) grid (9,4);

\filldraw[red] (6,3) circle (2pt) node[anchor=east] {};
\filldraw[red] (6,1) circle (2pt) node[anchor=east] {};
\filldraw[red] (8,1) circle (2pt) node[anchor=west] {};
\filldraw[red] (8,3) circle (2pt) node[anchor=west] {};
\draw[red, ultra thick] (6,3) -- (6,1) -- (8,1) -- (8,3);

% level 2 - 1
\filldraw[blue] (5.5,3.5) circle (2pt) node[anchor=east] {};
\filldraw[blue] (5.5,2.5) circle (2pt) node[anchor=east] {};
\filldraw[blue] (6.5,2.5) circle (2pt) node[anchor=west] {};
\filldraw[blue] (6.5,3.5) circle (2pt) node[anchor=west] {};
\draw[blue, ultra  thick] (5.5,3.5) -- (5.5,2.5) -- (6.5,2.5) -- (6.5,3.5);

% level 2 - 2
\filldraw[blue] (5.5,1.5) circle (2pt) node[anchor=east] {};
\filldraw[blue] (5.5,0.5) circle (2pt) node[anchor=east] {};
\filldraw[blue] (6.5,0.5) circle (2pt) node[anchor=west] {};
\filldraw[blue] (6.5,1.5) circle (2pt) node[anchor=west] {};
\draw[blue, ultra  thick] (5.5,1.5) -- (5.5,0.5) -- (6.5,0.5) -- (6.5,1.5);

% level 2 - 3
\filldraw[blue] (7.5,1.5) circle (2pt) node[anchor=east] {};
\filldraw[blue] (7.5,0.5) circle (2pt) node[anchor=east] {};
\filldraw[blue] (8.5,0.5) circle (2pt) node[anchor=west] {};
\filldraw[blue] (8.5,1.5) circle (2pt) node[anchor=west] {};
\draw[blue, ultra  thick] (7.5,1.5) -- (7.5,0.5) -- (8.5,0.5) -- (8.5,1.5);

% level 2 - 4
\filldraw[blue] (7.5,3.5) circle (2pt) node[anchor=east] {};
\filldraw[blue] (7.5,2.5) circle (2pt) node[anchor=east] {};
\filldraw[blue] (8.5,2.5) circle (2pt) node[anchor=west] {};
\filldraw[blue] (8.5,3.5) circle (2pt) node[anchor=west] {};
\draw[blue, ultra  thick] (7.5,3.5) -- (7.5,2.5) -- (8.5,2.5) -- (8.5,3.5);

% level 1 and 2 rotated --------------------------------------------------------------------
\node[] at (12,-0.5) {\small (c)};
\draw[step=1cm,gray,thin] (10,0) grid (14,4);
%\draw[step=2cm,darkgray,thin] (10,0) grid (14,4);

\filldraw[red] (11,3) circle (2pt) node[anchor=east] {};
\filldraw[red] (11,1) circle (2pt) node[anchor=east] {};
\filldraw[red] (13,1) circle (2pt) node[anchor=west] {};
\filldraw[red] (13,3) circle (2pt) node[anchor=west] {};
\draw[red, ultra thick] (11,3) -- (11,1) -- (13,1) -- (13,3);

% level 2 - 1
\filldraw[blue] (10.5,3.5) circle (2pt) node[anchor=east] {};
\filldraw[blue] (10.5,2.5) circle (2pt) node[anchor=east] {};
\filldraw[blue] (11.5,2.5) circle (2pt) node[anchor=west] {};
\filldraw[blue] (11.5,3.5) circle (2pt) node[anchor=west] {};
\draw[blue, ultra  thick] (10.5,3.5) -- (11.5,3.5) -- (11.5,2.5) -- (10.5,2.5);

% level 2 - 2
\filldraw[blue] (10.5,1.5) circle (2pt) node[anchor=east] {};
\filldraw[blue] (10.5,0.5) circle (2pt) node[anchor=east] {};
\filldraw[blue] (11.5,0.5) circle (2pt) node[anchor=west] {};
\filldraw[blue] (11.5,1.5) circle (2pt) node[anchor=west] {};
\draw[blue, ultra  thick] (10.5,1.5) -- (10.5,0.5) -- (11.5,0.5) -- (11.5,1.5);

% level 2 - 3
\filldraw[blue] (12.5,1.5) circle (2pt) node[anchor=east] {};
\filldraw[blue] (12.5,0.5) circle (2pt) node[anchor=east] {};
\filldraw[blue] (13.5,0.5) circle (2pt) node[anchor=west] {};
\filldraw[blue] (13.5,1.5) circle (2pt) node[anchor=west] {};
\draw[blue, ultra  thick] (12.5,1.5) -- (12.5,0.5) -- (13.5,0.5) -- (13.5,1.5);

% level 2 - 4
\filldraw[blue] (12.5,3.5) circle (2pt) node[anchor=east] {};
\filldraw[blue] (12.5,2.5) circle (2pt) node[anchor=east] {};
\filldraw[blue] (13.5,2.5) circle (2pt) node[anchor=west] {};
\filldraw[blue] (13.5,3.5) circle (2pt) node[anchor=west] {};
\draw[blue, ultra  thick] (13.5,3.5) -- (12.5,3.5) -- (12.5,2.5) -- (13.5,2.5);

% level 1 and 2 joined  --------------------------------------------------------------------
\node[] at (17,-0.5) {\small (d)};
\draw[step=1cm,gray,thin] (15,0) grid (19,4);
%\draw[step=2cm,darkgray,thin] (15,0) grid (19,4);

\filldraw[red] (16,3) circle (2pt) node[anchor=east] {};
\filldraw[red] (16,1) circle (2pt) node[anchor=east] {};
\filldraw[red] (18,1) circle (2pt) node[anchor=west] {};
\filldraw[red] (18,3) circle (2pt) node[anchor=west] {};
\draw[red, ultra thick] (16,3) -- (16,1) -- (18,1) -- (18,3);

% level 2 - 1
% small footnotesize scriptsize
\filldraw[blue] (15.5,3.5) circle (2pt) node[anchor=east] {\scriptsize 1};
\filldraw[blue] (15.5,2.5) circle (2pt) node[anchor=east] {\scriptsize 4};
\filldraw[blue] (16.5,2.5) circle (2pt) node[anchor=north] {\scriptsize 3};
\filldraw[blue] (16.5,3.5) circle (2pt) node[anchor=south] {\scriptsize 2};
\draw[blue, ultra  thick] (15.5,3.5) -- (16.5,3.5) -- (16.5,2.5) -- (15.5,2.5);

% level 2 - 2
\filldraw[blue] (15.5,1.5) circle (2pt) node[anchor=east] {\scriptsize 5};
\filldraw[blue] (15.5,0.5) circle (2pt) node[anchor=east] {\scriptsize 6};
\filldraw[blue] (16.5,0.5) circle (2pt) node[anchor=north] {\scriptsize 7};
\filldraw[blue] (16.5,1.5) circle (2pt) node[anchor=south] {\scriptsize 8};
\draw[blue, ultra  thick] (15.5,1.5) -- (15.5,0.5) -- (16.5,0.5) -- (16.5,1.5);

% level 2 - 3
\filldraw[blue] (17.5,1.5) circle (2pt) node[anchor=south] {\scriptsize 9};
\filldraw[blue] (17.5,0.5) circle (2pt) node[anchor=north] {\scriptsize 10};
\filldraw[blue] (18.5,0.5) circle (2pt) node[anchor=west] {\scriptsize 11};
\filldraw[blue] (18.5,1.5) circle (2pt) node[anchor=west] {\scriptsize 12};
\draw[blue, ultra  thick] (17.5,1.5) -- (17.5,0.5) -- (18.5,0.5) -- (18.5,1.5);

% level 2 - 4
\filldraw[blue] (17.5,3.5) circle (2pt) node[anchor=south] {\scriptsize 15};
\filldraw[blue] (17.5,2.5) circle (2pt) node[anchor=north] {\scriptsize 14};
\filldraw[blue] (18.5,2.5) circle (2pt) node[anchor=west] {\scriptsize 13};
\filldraw[blue] (18.5,3.5) circle (2pt) node[anchor=west] {\scriptsize 16};
\draw[blue, ultra  thick] (18.5,3.5) -- (17.5,3.5) -- (17.5,2.5) -- (18.5,2.5);

% joins

\draw[blue, ultra  thick, dashed]  (15.5,2.5) -- (15.5,1.5) ;
\draw[blue, ultra  thick, dashed]  (16.5,1.5) -- (17.5,1.5) ;
\draw[blue, ultra  thick, dashed]  (18.5,1.5) -- (18.5,2.5) ;

\node[] at (22,-0.5) {\small (e)};
\draw[step=1cm,gray,thin] (20,0) grid (24,4);

\foreach \x in {0,...,3}
  \foreach \y [count=\yi] in {0,...,3}  
    \filldraw[color=green!60!black] (20.5+\x,0.5+\y) circle (3pt);
\draw[green!60!black, very thick] (20.5,0.5) -- (23.5,0.5);
\draw[green!60!black, very thick] (20.5,1.5) -- (23.5,1.5);
\draw[green!60!black, very thick] (20.5,2.5) -- (23.5,2.5);
\draw[green!60!black, very thick] (20.5,3.5) -- (23.5,3.5);
\draw[green!60!black, very thick] (20.5,0.5) -- (20.5,3.5);
\draw[green!60!black, very thick] (21.5,0.5) -- (21.5,3.5);
\draw[green!60!black, very thick] (22.5,0.5) -- (22.5,3.5);
\draw[green!60!black, very thick] (23.5,0.5) -- (23.5,3.5);

\end{tikzpicture}
\caption{Hilbert space-filling curves at level 1 and 2: (a)~shows the fundamental shape of the Hilbert curve, which also corresponds to the level~1 Hilbert curve. For level~2, (b)~locates the fundamental shape at the vertices of level~1; these shapes are rotated~(c) and linked~(d) to form a continuous curve. The graph of the level~2 discretisation is shown in~(e) with vertices (green) located at the centre of the control volumes, edges (green) and the grid corresponding to level~2 shown in grey. }
\label{fig:Hilbert_curve_level_1_and_2}
\end{figure}
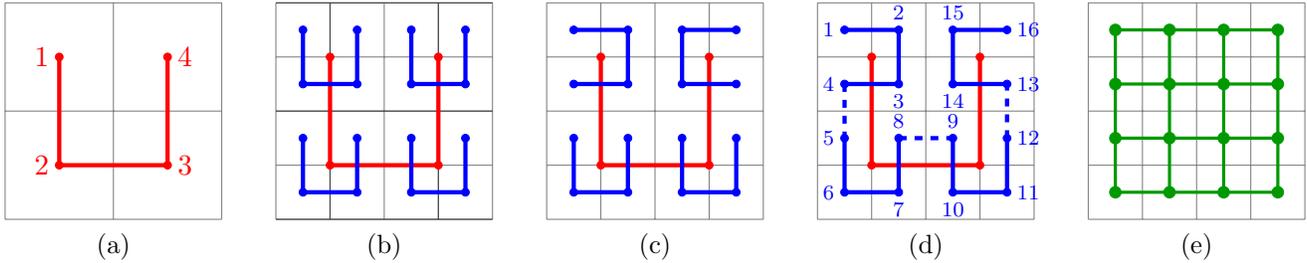

\begin{figure}[htbp]
\centering
\includegraphics[width=11cm]{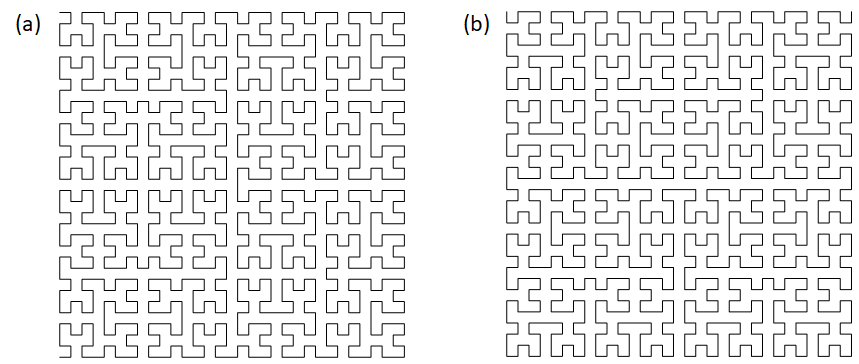}
\caption{Two 2D Hilbert space-filling curves for a ${32 \times 32}$ grid, where (b)~is a 90 degree clockwise rotation of~(a). Both are used to reorder data in the two examples for data held on structured grids in the results section.}
\label{2d-3d-space_filling}
\end{figure}

Given its fundamental shape of U, a Hilbert curve can be formed by placing the fundamental shape on the coarsest level (level~1) of a grid.  For a Hilbert curve, the level~one representation fits on a 2~by~2 grid, see Figure~\ref{fig:Hilbert_curve_level_1_and_2}(a). Once the starting vertex is chosen (from either vertex with a valency of~1), the vertices are numbered as ones traverses the path.  The numbering is not unique as the path could also be followed in the opposite direction. To form the level~2 Hilbert curve (on a 4~by~4 grid), the fundamental shape is centred on each vertex of the level~1 curve (Figure~\ref{fig:Hilbert_curve_level_1_and_2}(b)) and rotated (Figure~\ref{fig:Hilbert_curve_level_1_and_2}(c)) until the shapes can be linked with one another by connections (Figure~\ref{fig:Hilbert_curve_level_1_and_2}(d)) that (i)~are consistent with the discretisation stencil and (ii)~result in a continuous curve.  Here the stencil implied by the Hilbert curve’s fundamental shape has connections between horizontal and vertical nearest neighbours - known as a five point stencil and often used in solving differential equations. This process of locating the fundamental shape at vertices and rotating them until a continuous curve is found, is repeated until the desired level is reached. Neighbouring vertices on the space-filling curve in~(d) will be close in the graph for level~2 shown in~(e) but the reverse is not necessarily true. For example, vertices~2 and~15 in Figure~\ref{fig:Hilbert_curve_level_1_and_2}(d) are not close on the space-filling curve but are close in the original graph in~(e). To overcome potential inaccuracies caused by the disconnect between neighbouring vertices on the original graph, the generation of multiple space-filling curves is explored, the aim of which is that such vertices would be closer on a subsequent space-filling curve. More details are given about multiple space-filling curves throughout this section. As an example, two 2D Hilbert curves for a 32~by~32 grid can be seen in Figure~\ref{2d-3d-space_filling}; both are used to generate the results for the structured grid test cases that are presented later.

The approach we take to form space-filling curves for unstructured meshes is similar to the above. For the fundamental shape, we choose a straight line (two vertices and one edge). The two vertices represent two partitions and the edge between them represents an edge on the original graph which connects them. Generating a space-filling curve consists of three steps. First, the graph is constructed, based on the discretisation stencil associated with the nodes or cells of the mesh. In the graph, vertices represent the nodes or cells, and edges represent the connections between nodes or cells as determined by the stencil. We refer to this as the `original graph' to distinguish it from graphs of the partitions at intermediate levels of decomposition. Second, a hierarchy of partitions (and graphs) are created from the original graph using a domain decomposition approach based on nested bisection. The graph partitioner will be explained in the next section, but, important here is that it generates a series of partitions at different levels of refinement. Third, the partitions are ordered, beginning at the coarsest level and ending at the final level. Once complete, the space-filling curve ordering is obtained, either directly from the vertex ordering of the final level, or after an additional step if there is more than one vertex in any of the partitions.

%The domain is partitioned into two (as the fundamental shape has two vertices) and the fundamental shape is placed with one vertex in each partition which are numbered arbitrarily as partition~1 and~2. The graph of the mesh is consulted to ensure that there is an edge connecting the two partitions (this should always be true for the first decomposition). Each partition is then further decomposed into two partitions: thus, partition~1 (from level~1) now contains two partitions numbered~1 and~2, and partition~2 contains two partitions numbered~3 and~4. We know that there will be an edge connecting partitions~1 and~2 as these partitions originate from the same partition (one level higher); the same applies to partitions~3 and~4. However, we need to make sure that there is a connection between partition~2 and~3. This could be done in various ways. Here, we chose to define a functional, based on the number of edges traversed when travelling from partition~$p$ to partition~$q$. The numbering or ordering of partitions which gives the lowest value of the functional is accepted as the optimal path from partition~$p$ to partition~$q$. If two numberings give the same value for the functional, the most recent one to be calculated is accepted. The partition numbers are swapped in groups of four, with the functional calculating the the cost of traversing each path defined by these numberings only considering the\ldots , 

\begin{figure}[htbp]
\centering
%\begin{tikzpicture}[scale=0.6] % 0.6
\input{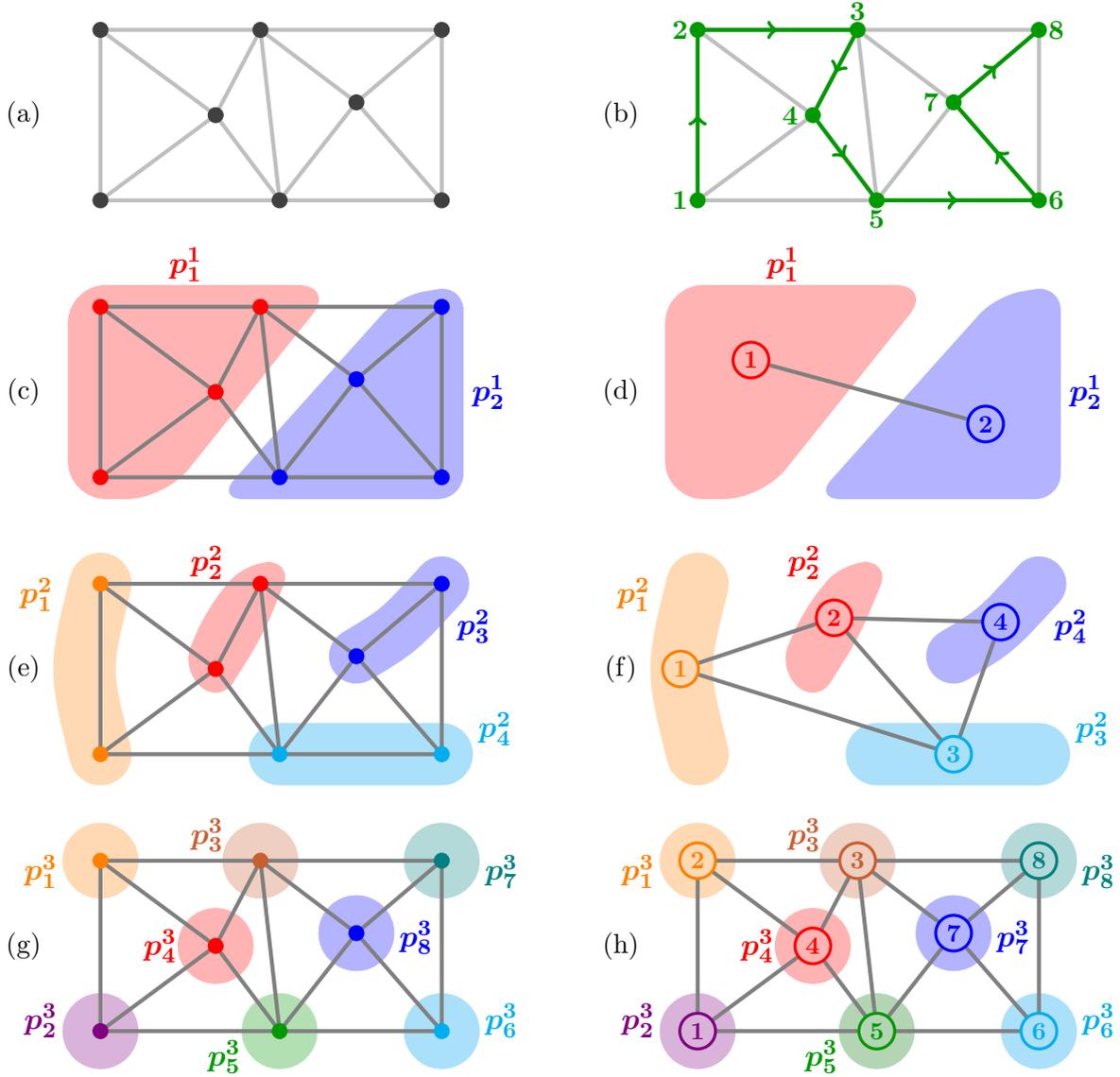}
\caption{The graph of the discretisation (2D linear finite elements) is shown in~(a). The partitions at levels~1, 2 and~3 are shown in (c), (e) and (g) respectively, labelled as $p_i^\ell$ where $i$ is the partition number and $\ell$~is the level. The partitions are indicated by a shaded region and coloured vertices. The graph of the discretisation is also shown in these plots. Shown in plots~(d), (f) and~(h) are the partitions (indicated as shaded regions) and graphs of their connectivities (shown in grey) for levels~1, 2 and~3 respectively. In the graphs, the vertices (indicated as numbers in circles) have been ordered according to the procedure outlined in the text.  The final space-filling curve is shown in green in plot~(b). }
\label{fig:sfc_unstructured_meshes}
\end{figure}

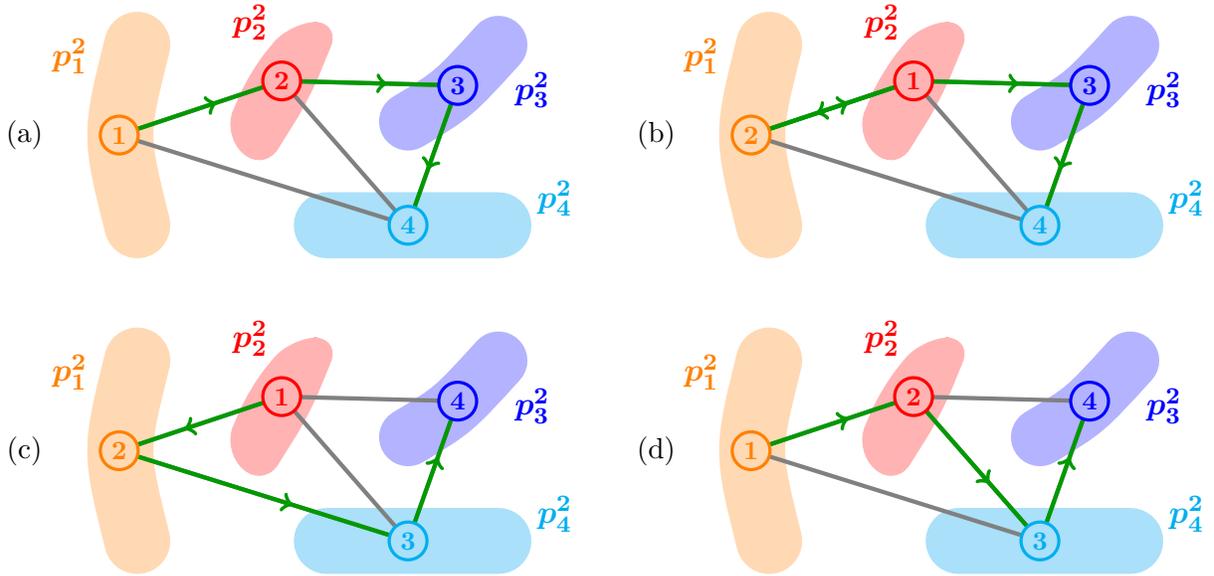
\begin{figure}[htbp]
    \centering
    \tikzset{every state/.style={minimum size=4pt, inner sep=2pt}}

\begin{tikzpicture}[scale=0.6] % 0.6

% level 2 graph v1 -----------------------------------------------------------
\begin{scope}[shift={(0,0)}]

    \coordinate (A) at (0, 0);
    \coordinate (B) at (0, 4);
    \coordinate (C) at (2+0.3+0.4, 2-0.2+0.2);
    \coordinate (D) at (4+0.2, 0);
    \coordinate (E) at (4-0.25, 4);
    \coordinate (F) at (6, 2+0.3);
    \coordinate (G) at (8, 0);
    \coordinate (H) at (8, 4);

%    \node[] at (-1.8,2) {(a)};
    \node[] at (-2.5,2) {(a)};

    \draw[line width=25pt, color=orange!30!white, cap=round] (A) .. controls (-0.5,2) .. (B) ; % .. controls
    \fill[color=red!30!white, draw=red!30!white, rounded corners = 10]  (1.9, 2) -- (2.8, 1.2) -- (3.8,2.8) -- (4.5, 4.5) -- (3.4, 4.5) -- (2.5, 3.4) -- cycle;
    \draw[line width=25pt, color=cyan!30!white, cap=round] (D) -- (G) ; 
    \draw[line width=22pt, color=blue!30!white, cap=round] (F) ..  controls  (7.0,2.9) .. (H) ; 

    \draw[gray, ultra thick] (-0.4,2) -- (3.2,3.2) -- (7.1,3.1) -- (6,0); % 2-1-3-4
    \draw[gray, ultra thick] (3.2,3.2) -- (6,0) -- (-0.4,2);              % 1-4-2

    \begin{scope}[very thick,decoration={
        markings,
        mark=at position 0.6 with {\arrow{>}}}
        ] 
        \draw[postaction={decorate}, green!60!black, ultra thick] (-0.4,2) -- (3.2,3.2);
        \draw[postaction={decorate}, green!60!black, ultra thick] (3.2,3.2) -- (7.1,3.1);
        \draw[postaction={decorate}, green!60!black, ultra thick] (7.1,3.1) -- (6.0,0.0);

    \end{scope}

    \node[state, color=orange, very thick, fill=orange!30!white] at (-0.4,2)  {\textbf{\small{1}}};
    \node[state, color=red, very thick, fill=red!30!white]       at (3.2,3.2) {\textbf{\small{2}}};
    \node[state, color=blue, very thick, fill=blue!30!white]     at (7.1,3.1) {\textbf{\small{3}}};
    \node[state, color=cyan, very thick, fill=cyan!30!white]     at (6,0)     {\textbf{\small{4}}};

    \node[] at (-1.5,3.75) { \textcolor{orange}{$\Scale[1.1]{\bm{p_1^2}}$}}; % (1,2)
    \node[] at (2.5,4.5) { \textcolor{red}{ $\Scale[1.1]{\bm{p_2^2}}$} }; 
    \node[] at (8.75,3) { \textcolor{blue}{  $\Scale[1.1]{\bm{p_3^2}}$}}; % (1,2)
    \node[] at (9.25,0.6) { \textcolor{cyan}{$\Scale[1.1]{\bm{p_4^2}}$} }; 

\end{scope}

% level 2 graph v2 -----------------------------------------------------------
\begin{scope}[shift={(14,0)}]

    \coordinate (A) at (0, 0);
    \coordinate (B) at (0, 4);
    \coordinate (C) at (2+0.3+0.4, 2-0.2+0.2);
    \coordinate (D) at (4+0.2, 0);
    \coordinate (E) at (4-0.25, 4);
    \coordinate (F) at (6, 2+0.3);
    \coordinate (G) at (8, 0);
    \coordinate (H) at (8, 4);

    \node[] at (-2.5,2) {(b)};

    \draw[line width=25pt, color=orange!30!white, cap=round] (A) .. controls (-0.5,2) .. (B) ; % .. controls
    \fill[color=red!30!white, draw=red!30!white, rounded corners = 10]  (1.9, 2) -- (2.8, 1.2) -- (3.8,2.8) -- (4.5, 4.5) -- (3.4, 4.5) -- (2.5, 3.4) -- cycle;
    \draw[line width=25pt, color=cyan!30!white, cap=round] (D) -- (G) ; 
    \draw[line width=22pt, color=blue!30!white, cap=round] (F) ..  controls  (7.0,2.9) .. (H) ; 

    \draw[gray, ultra thick] (-0.4,2) -- (3.2,3.2) -- (7.1,3.1) -- (6,0); % 2-1-3-4
    \draw[gray, ultra thick] (3.2,3.2) -- (6,0) -- (-0.4,2);              % 1-4-2

    \begin{scope}[very thick,decoration={
        markings,
        mark=at position 0.6 with {\arrow{>}}}
        ] 
        \draw[postaction={decorate}, green!60!black, ultra thick] (3.2,3.2) -- (-0.4,2);
        \draw[postaction={decorate}, green!60!black, ultra thick] (-0.4,2) -- (3.2,3.2);
        \draw[postaction={decorate}, green!60!black, ultra thick] (3.2,3.2) -- (7.1,3.1);
        \draw[postaction={decorate}, green!60!black, ultra thick] (7.1,3.1) -- (6.0,0.0);
    \end{scope}

    \node[state, color=orange, very thick, fill=orange!30!white] at (-0.4,2)  {\textbf{\small{2}}};
    \node[state, color=red, very thick, fill=red!30!white]       at (3.2,3.2) {\textbf{\small{1}}};
    \node[state, color=blue, very thick, fill=blue!30!white]     at (7.1,3.1) {\textbf{\small{3}}};
    \node[state, color=cyan, very thick, fill=cyan!30!white]     at (6,0)     {\textbf{\small{4}}};

    \node[] at (-1.5,3.75) { \textcolor{orange}{$\Scale[1.1]{\bm{p_1^2}}$}}; % (1,2)
    \node[] at (2.5,4.5) { \textcolor{red}{ $\Scale[1.1]{\bm{p_2^2}}$} }; 
    \node[] at (8.75,3) { \textcolor{blue}{  $\Scale[1.1]{\bm{p_3^2}}$}}; % (1,2)
    \node[] at (9.25,0.6) { \textcolor{cyan}{$\Scale[1.1]{\bm{p_4^2}}$} }; 

\end{scope}

% level 2 graph v3 -----------------------------------------------------------
\begin{scope}[shift={(0,-7)}]

    \coordinate (A) at (0, 0);
    \coordinate (B) at (0, 4);
    \coordinate (C) at (2+0.3+0.4, 2-0.2+0.2);
    \coordinate (D) at (4+0.2, 0);
    \coordinate (E) at (4-0.25, 4);
    \coordinate (F) at (6, 2+0.3);
    \coordinate (G) at (8, 0);
    \coordinate (H) at (8, 4);

    \node[] at (-2.5,2) {(c)};

    \draw[line width=25pt, color=orange!30!white, cap=round] (A) .. controls (-0.5,2) .. (B) ; % .. controls
    \fill[color=red!30!white, draw=red!30!white, rounded corners = 10]  (1.9, 2) -- (2.8, 1.2) -- (3.8,2.8) -- (4.5, 4.5) -- (3.4, 4.5) -- (2.5, 3.4) -- cycle;
    \draw[line width=25pt, color=cyan!30!white, cap=round] (D) -- (G) ; 
    \draw[line width=22pt, color=blue!30!white, cap=round] (F) ..  controls  (7.0,2.9) .. (H) ; 

    \draw[gray, ultra thick] (-0.4,2) -- (3.2,3.2) -- (7.1,3.1) -- (6,0); % 2-1-3-4
    \draw[gray, ultra thick] (3.2,3.2) -- (6,0) -- (-0.4,2);              % 1-4-2

    \begin{scope}[very thick,decoration={
        markings,
        mark=at position 0.6 with {\arrow{>}}}
        ] 
        \draw[postaction={decorate}, green!60!black, ultra thick] (3.2,3.2) -- (-0.4,2);
        \draw[postaction={decorate}, green!60!black, ultra thick] (-0.4,2) -- (6.0,0.0);
        \draw[postaction={decorate}, green!60!black, ultra thick] (6.0,0.0) -- (7.1,3.1);
    \end{scope}

    \node[state, color=orange, very thick, fill=orange!30!white] at (-0.4,2)  {\textbf{\small{2}}};
    \node[state, color=red, very thick, fill=red!30!white]       at (3.2,3.2) {\textbf{\small{1}}};
    \node[state, color=blue, very thick, fill=blue!30!white]     at (7.1,3.1) {\textbf{\small{4}}};
    \node[state, color=cyan, very thick, fill=cyan!30!white]     at (6,0)     {\textbf{\small{3}}};

    \node[] at (-1.5,3.75) { \textcolor{orange}{$\Scale[1.1]{\bm{p_1^2}}$}}; % (1,2)
    \node[] at (2.5,4.5) { \textcolor{red}{ $\Scale[1.1]{\bm{p_2^2}}$} }; 
    \node[] at (8.75,3) { \textcolor{blue}{  $\Scale[1.1]{\bm{p_3^2}}$}}; % (1,2)
    \node[] at (9.25,0.6) { \textcolor{cyan}{$\Scale[1.1]{\bm{p_4^2}}$} }; 

\end{scope}

% level 2 graph v4 -----------------------------------------------------------
\begin{scope}[shift={(14,-7)}]

    \coordinate (A) at (0, 0);
    \coordinate (B) at (0, 4);
    \coordinate (C) at (2+0.3+0.4, 2-0.2+0.2);
    \coordinate (D) at (4+0.2, 0);
    \coordinate (E) at (4-0.25, 4);
    \coordinate (F) at (6, 2+0.3);
    \coordinate (G) at (8, 0);
    \coordinate (H) at (8, 4);

    \node[] at (-2.5,2) {(d)};

    \draw[line width=25pt, color=orange!30!white, cap=round] (A) .. controls (-0.5,2) .. (B) ; % .. controls
    \fill[color=red!30!white, draw=red!30!white, rounded corners = 10]  (1.9, 2) -- (2.8, 1.2) -- (3.8,2.8) -- (4.5, 4.5) -- (3.4, 4.5) -- (2.5, 3.4) -- cycle;
    \draw[line width=25pt, color=cyan!30!white, cap=round] (D) -- (G) ; 
    \draw[line width=22pt, color=blue!30!white, cap=round] (F) ..  controls  (7.0,2.9) .. (H) ; 

    \draw[gray, ultra thick] (-0.4,2) -- (3.2,3.2) -- (7.1,3.1) -- (6,0); % 2-1-3-4
    \draw[gray, ultra thick] (3.2,3.2) -- (6,0) -- (-0.4,2);              % 1-4-2

    \begin{scope}[very thick,decoration={
        markings,
        mark=at position 0.6 with {\arrow{>}}}
        ] 
        \draw[postaction={decorate}, green!60!black, ultra thick] (-0.4,2) -- (3.2,3.2);
        \draw[postaction={decorate}, green!60!black, ultra thick] (3.2,3.2) -- (6.0,0.0);
        \draw[postaction={decorate}, green!60!black, ultra thick] (6.0,0.0) -- (7.1,3.1);
    \end{scope}

    \node[state, color=orange, very thick, fill=orange!30!white] at (-0.4,2)  {\textbf{\small{1}}};
    \node[state, color=red, very thick, fill=red!30!white]       at (3.2,3.2) {\textbf{\small{2}}};
    \node[state, color=blue, very thick, fill=blue!30!white]     at (7.1,3.1) {\textbf{\small{4}}};
    \node[state, color=cyan, very thick, fill=cyan!30!white]     at (6,0)     {\textbf{\small{3}}};

    \node[] at (-1.5,3.75) { \textcolor{orange}{$\Scale[1.1]{\bm{p_1^2}}$}}; % (1,2)
    \node[] at (2.5,4.5) { \textcolor{red}{ $\Scale[1.1]{\bm{p_2^2}}$} }; 
    \node[] at (8.75,3) { \textcolor{blue}{  $\Scale[1.1]{\bm{p_3^2}}$}}; % (1,2)
    \node[] at (9.25,0.6) { \textcolor{cyan}{$\Scale[1.1]{\bm{p_4^2}}$} }; 

\end{scope}

\end{tikzpicture}
    \caption{The four possible combinations of the vertex ordering at level~2, see Table~\ref{table:partition_numbering_level2}. The partitions are shown, labelled $p_i^\ell$ for partition~$i$ at level~$\ell=2$. The graph is shown in grey and the path corresponding to each ordering is highlighted in green. Once the ordering has been determined (in this case, that shown in~(d)) the partitions will be renumbered, see Figure~\ref{fig:sfc_unstructured_meshes}(f). }
    \label{fig:graphs_at_level_2}
\end{figure}

The starting point for generating the space-filling curve ordering or numbering is therefore a hierarchy of nested partitions. As a consequence of the fundamental shape, at level~1, the coarsest level, there will be two partitions; at level~2, there will be 4~partitions; at level~$\ell$ there will be $2^{\ell}$ partitions and so on. In addition to the original graph, at every level a graph of the partitions and their connectivities is constructed, from which the numbering is calculated. The method is now illustrated with reference to an example shown in Figure~\ref{fig:sfc_unstructured_meshes}. Here, a 2D linear Finite Element (FE) discretisation is assumed which uses three-noded elements. The original graph of the computational domain and discretisation can be seen in Figure~\ref{fig:sfc_unstructured_meshes}(a), which, for this stencil, coincides with the FE mesh. %\textcolor{red}{\textit{repetitious:} The vertices of the graph are the nodes of the discretisation and the stencil is such that the edges in the graph coincide exactly with the mesh.} 
Based on the original graph~(a), partitions at three levels, seen in~(c), (e) and~(g) respectively, are obtained from the graph partitioner. For each level, the ordering algorithm constructs a graph to represent the partitions and their connections seen in~(d), (f), and~(h). For level~1, the graph shown in~(d) consists of two vertices (one in each partition) joined by an edge. The vertices in partitions $p_1^{\ell=1}$ and $p_2^{\ell=1}$ can be numbered arbitrarily as these partitions have the same parent and will be connected by an edge on the original graph. Once the vertices are numbered, the partitions are renumbered to be consistent with the vertex numbering, i.e.~so that vertex~$k$ is located in partition $p_k^\ell$. For level~2, the first partition from level~1 ($p_{\ell=1}^1$) gives rise to partitions $p_1^{\ell=2}$ and $p_2^{\ell=2}$, and the second partition from level~1 gives rise to $p_3^{\ell=2}$ and $p_4^{\ell=2}$. Partitions $p_1^{\ell=2}$ and $p_2^{\ell=2}$ will be connected by an edge on the original graph, as will partitions $p_3^{\ell=2}$ and $p_4^{\ell=2}$, for the reason that they share the same parent partition. What must be found is the connection between either $p_1^{\ell=2}$ or $p_2^{\ell=2}$, and $p_3^{\ell=2}$ or $p_4^{\ell=2}$. Incorporating this edge into the space-filling curve gives it continuity, and  ideally, therefore, there would be just one edge needed to move from $p^{\ell=1}_1$ or $p^{\ell=1}_2$ to $p^{\ell=1}_3$ or $p^{\ell=1}_4$. This is so the path traced out by the space-filling curve traverses the vertices in as efficient a  manner as possible (given the stencil). At this level, there are four possible ways of numbering the vertices, resulting in four different paths through the graph. All combinations of numbering are tried and, for each, a functional is evaluated to find the best path or any one of the best paths. If one combination gives either the same value or an improved functional value, it is accepted. The functional is defined as the number of edges traversed in order to move from the first vertex to the fourth. In Figure~\ref{fig:graphs_at_level_2} we show the four possible combinations of numbering, also shown in Table~\ref{table:partition_numbering_level2} along with the corresponding functional values. In Figure~\ref{fig:graphs_at_level_2}(a) the vertices are ordered according to the partition number, which produces a path with 3~edges. Vertices~1 and~2 are swapped in Figure~\ref{fig:graphs_at_level_2}(b) and the number of edges in the path from vertex~1 to~4 is~4. This path is rejected as the functional value has not reduced. Vertices~3 and~4 are swapped in Figure~\ref{fig:graphs_at_level_2}(c) producing a path of 3~edges which is accepted. Finally vertices~1 and~2 are swapped back, shown in Figure~\ref{fig:graphs_at_level_2}(d), also  producing a path of 3~edges which is accepted. This is chosen as the final path and the partition numbers are updated to match the vertex numbers, shown in Figure~\ref{fig:sfc_unstructured_meshes}(f).

\begin{table}[htbp]
\centering \ra{1.3}
\begin{tabular}{l l llll l llc}
\toprule
%     & \multicolumn{4}{c}{level 2 partitions} &&&\\
      & \phantom{a} & $p_1^{\ell=2}$ & $p_2^{\ell=2}$ & $p_3^{\ell=2}$ & $p_4^{\ell=2}$ & \phantom{a} & $f^{1\rightarrow 4}$ & path & Figure\\
\midrule
\multirow{4}{*}{vertex ordering} & & 1 & 2 & 3 & 4 & & 3 & accepted     & \ref{fig:graphs_at_level_2}(a) \\
                       & & 2 & 1 & 3 & 4 & & 4 & not accepted & \ref{fig:graphs_at_level_2}(b) \\
                       & & 2 & 1 & 4 & 3 & & 3 & accepted     & \ref{fig:graphs_at_level_2}(c) \\
                       & & 1 & 2 & 4 & 3 & & 3 & accepted     & \ref{fig:graphs_at_level_2}(d) \\
\bottomrule
\end{tabular}
\caption{The four possible combinations of the vertex ordering at level~2. In the first row, the first of the four paths is shown, where the first vertex in the path is that in partition~$p_{1}^{2}$, the second vertex is~$p_{2}^{2}$ and so on. Each path is illustrated in Figure~\ref{fig:graphs_at_level_2} from which the functional $f^{v\rightarrow w}$ can be calculated. The functional is defined as the minimum number of edges that are traversed in going from vertex~$v$ to vertex~$w$.}
\label{table:partition_numbering_level2}
\end{table}

For level~3, there are 8~partitions, see Figure~\ref{fig:sfc_unstructured_meshes}(g). The graph for this level can be seen in Figure~\ref{fig:sfc_unstructured_meshes}(h). At this stage, there is only one vertex in each partition, so the graph resembles the FE mesh. As the number of vertices increases, the number of combinations grows exponentially, so at level~3, if the entire path is considered, that would result in~16 possible orderings. Instead of considering all the partitions simultaneously, groups of 4~partitions are considered in turn, as shown in Table~\ref{table:partition_numbering_level3}. First, the four combinations of ordering the vertices in partitions~1 to 4 are considered by evaluating the functional from vertex~1 to~4, $f^{1\rightarrow 4}$. Once the optimal ordering of these vertices has been found, vertices~1 and~2 are fixed and the vertices in partitions~3 to~6 are considered. The functional is now evaluated from the second to the sixth vertex, $f^{2\rightarrow 6}$ as the location of the third vertex relative to the second will affect the optimality of the path. When the ordering has been found, vertices~3 and~4 are fixed and now the vertices in partitions~5 to~8 are considered. The algorithm sweeps forwards and backwards through the partitions, considering four vertices at a time and fixing pairs of vertices based on the optimal path calculated between the appropriate number of vertices (4, 5 or 6) depending on whether the vertices either-side of the path have been fixed (see Table~\ref{table:partition_numbering_level3}). For simple cases, only a small number of iterations are sufficient. For complex cases, more iterations are likely to be needed. In any case, it is not known whether the entire path from vertex~1 to vertex~$N^{\text{vert}}$ ($N^{\text{vert}}$ is the total number of vertices or nodes in the graph) produced by this method will be optimal, however, using 10 forwards and backwards `sweeps' was found to be sufficient for the challenging example of an unstructured mesh with a discontinuous Galerkin stencil.  
Once the final level is reached, if there are multiple vertices in a partition, a series of vertex swapping operations are performed until an optimal path is found within that partition: that is a fully-connected path that links to the surrounding partitions and which has the shortest path. Hence, the space-filling curve ordering is obtained, of which, one possible solution is shown in Figure~\ref{fig:sfc_unstructured_meshes}(b) 

%For the first four partitions, when evaluating the functional the path between the four combinations of vertices~1 to~4 is calculated. Next the path between vertices~2 to~6 is calculated. Vertex~2 remains fixed so there are still 4 combinations or vertices~3 to~6, however, vertex~2 must be included as the path starts here. Next the path between vertices~4 and~8 is considered. Once the optimal path is found, the iteration sweeps backwards through the partitions (groups in fours, every two). Now the end points of the path are also included in the functional evaluations. 
%In practice, we sweep through all the partitions in a particular level swapping over the vertex numbers going from the first to the last partition in each level and then sweeping back from the last to the first, see table~\ref{table:partition_numbering_level3}. We then repeat this 10 times or until we do not change the ordering of the partitions. 
%If, after the final decomposition, there are multiple vertices in a partition, a series of vertex swapping operations are performed until an optimal path is found within that partition, that is a fully-connected path that links to the surrounding partitions and which has the shortest path.
%This process is summarised in Algorithm~\ref{alg:sfc_numbering}.

\begin{table}[htbp]
\centering \ra{1.1}
\begin{tabular}{l l llllllll l l}
\toprule
%     & \multicolumn{4}{c}{level 2 partitions} &&&\\
      & \phantom{a} & $p_1^{\ell=3}$ & $p_2^{\ell=3}$ & $p_3^{\ell=3}$ & $p_4^{\ell=3}$ & $p_5^{\ell=3}$ & $p_6^{\ell=3}$ & $p_7^{\ell=3}$ & $p_8^{\ell=3}$ & \phantom{a} & $f^{\text{path}}$ \\
\midrule
 & & 1 & 2 & 3 & 4 & \multirow{4}{*}{5} & \multirow{4}{*}{6} & \multirow{4}{*}{7} & \multirow{4}{*}{8} & & $f^{1\rightarrow 4}$ \\
                       & & 2 & 1 & 3 & 4 &  &  &  &  & & $f^{1\rightarrow 4}$ \\
                       
                       & & 1 & 2 & 4 & 3 &  &  &  &  & & $f^{1\rightarrow 4}$ \\
\multirow{4}{*}{\rotatebox[origin=c]{90}{\textit{\ldots~forward sweep~\ldots}}} & & 2 & 1 & 4 & 3 &  &  &  &  & & $f^{1\rightarrow 4}$ \\
                        \cmidrule{3-12} 
 & & \multirow{4}{*}{2} & \multirow{4}{*}{1} & 4 & 3 & 5 & 6 & \multirow{4}{*}{7} & \multirow{4}{*}{8} & & $f^{2\rightarrow 6}$ \\
                       & &  &  & 3 & 4 & 5 & 6 &  &  & & $f^{2\rightarrow 6}$ \\
                       & &  &  & 3 & 4 & 6 & 5 &  &  & & $f^{2\rightarrow 6}$ \\
                       & &  &  & 4 & 3 & 6 & 5 &  &  & & $f^{2\rightarrow 6}$ \\
                       \cmidrule{3-12}
  & & \multirow{4}{*}{2} & \multirow{4}{*}{1} & \multirow{4}{*}{4} & \multirow{4}{*}{3} & 5 & 6 & 7 & 8 & & $f^{4\rightarrow 8}$ \\
                    & &  &  &  &  & 6 & 5 & 7 & 8 & & $f^{4\rightarrow 8}$ \\
                    & &  &  &  &  & 6 & 5 & 8 & 7 & & $f^{4\rightarrow 8}$ \\
\multirow{4}{*}{\rotatebox[origin=c]{90}{\textit{\ldots~backward sweep~\ldots}}}    & &  &  &  &  & 5 & 6 & 8 & 7 & & $f^{4\rightarrow 8}$ \\
                        \cmidrule{3-12} 
 & & \multirow{4}{*}{2} & \multirow{4}{*}{1} & 4 & 3 & 5 & 6 & \multirow{4}{*}{8} & \multirow{4}{*}{7} & & $f^{2\rightarrow 7}$ \\
                       & &  &  & 3 & 4 & 5 & 6 &  &  & & $f^{2\rightarrow 7}$ \\
                       & &  &  & 3 & 4 & 6 & 5 &  &  & & $f^{2\rightarrow 7}$ \\
                       & &  &  & 4 & 3 & 6 & 5 &  &  & & $f^{2\rightarrow 7}$ \\
                       \cmidrule{3-12}
 & & 2 & 1 & 4 & 3 & \multirow{4}{*}{6} & \multirow{4}{*}{5} & \multirow{4}{*}{8} & \multirow{4}{*}{7} & & $f^{1\rightarrow 5}$ \\
                       & & 1 & 2 & 4 & 3 &  &  &  &  & & $f^{1\rightarrow 5}$ \\
                       & & 1 & 2 & 3 & 4 &  &  &  &  & & $f^{1\rightarrow 5}$ \\
                       & & 2 & 1 & 3 & 4 &  &  &  &  & & $f^{1\rightarrow 5}$ \\
\bottomrule
\end{tabular}
\caption{Determining the ordering or numbering of the vertices for the partitions after three levels of bisection. The partitions are labelled $p_i^\ell$ for partition index~$i$ and level~$\ell$. The first row indicates that the first path to be considered is partitions~1, 2, 3 then~4. The number of edges traversed in the graph when following this path is counted, that is, the functional from vertices~1 to~4 is evaluated. The second path to be considered is partitions~2, 1, 3 then~4. The functional for this path is evaluated. A certain number of iterations are carried out, which sweep forward and backwards through the partitions, at which point the space-filling curve path is obtained. }
\label{table:partition_numbering_level3}
\end{table}

%\textcolor{blue}{previous text:\\
%Effectively, the partition numbers are swapped until the best numbering is found according to the functional. This effectively tries every combination of vertex numbering in each partition and chooses the optimal node numbering - defined by the shortest path length in terms of number of edges in the graph traversed by the node numbering. }

%\begin{algorithm}[!htb]
%    \begin{algorithmic}
%    \While{\textcolor{red}{\textbf{some condition}}}
%        \For{each partition~$p$}
%        \State split the partition in two and number arbitrarily as partitions~$2p-1$ and~$2p$
%        \State if partitions~$2p-1$ and~$2p$ are not connected by an edge on the graph, swap the partition numbers $2p-1$ and~$2p$  until the best numbering is found (continue to do this until for lower numbered partitions)
%        \EndFor
%        \EndWhile
%        \label{alg:sfc_numbering}
%        \caption{Algorithm caption}
%    \end{algorithmic}
%\end{algorithm}

This method of generating space-filling curves is applied to an $8\times 8$ structured grid and the resulting space-filling curve is shown in Figure~\ref{unstructured_method_regular}. The graph for the $8\times 8$ grid was generated by assuming a 5~point finite-difference stencil: the same stencil that was used to generate the graph seen in Figure~\ref{fig:Hilbert_curve_level_1_and_2}(e). In Figure~\ref{unstructured_method_regular} it can be seen that a satisfactory  space-filling curve has been found for the grid, whilst using a method that can be applied to fully unstructured meshes. The curve found here has four instances where the path is longer than it needs to be (indicated by the diagonal lines). However, this path is close to optimal and deemed good enough for use here. For this grid, a Hilbert curve would represent an optimal solution. The number of edges traversed by the Hilbert curve would be~63. Our method returns a value of~67 as every diagonal line would count as two edges traversed on the stencil.

% % figure of the SFC approach applied to the 8 by 8 grid and DG stencil combined
\begin{figure}[htbp]
\centering
\begin{minipage}{0.6\textwidth}
\centering
%\documentclass{article}
%\usepackage{tikz}
%\usepackage[margin=2cm]{geometry}
%\begin{document}
\begin{tikzpicture}[scale=0.55] % 0.8

\node[] at (0,0) {\small{27}};
\node[] at (1,0) {\small{26}};
\node[] at (2,0) {\small{23}};
\node[] at (3,0) {\small{22}};
\node[] at (4,0) {\small{10}};
\node[] at (5,0) {\small{9}};
\node[] at (6,0) {\small{8}};
\node[] at (7,0) {\small{7}};

\node[] at (0,1) {\small{28}};
\node[] at (1,1) {\small{25}};
\node[] at (2,1) {\small{24}};
\node[] at (3,1) {\small{21}};
\node[] at (4,1) {\small{11}};
\node[] at (5,1) {\small{12}};
\node[] at (6,1) {\small{5}};
\node[] at (7,1) {\small{6}};

\node[] at (0,2) {\small{29}};
\node[] at (1,2) {\small{30}};
\node[] at (2,2) {\small{19}};
\node[] at (3,2) {\small{20}};
\node[] at (4,2) {\small{14}};
\node[] at (5,2) {\small{13}};
\node[] at (6,2) {\small{4}};
\node[] at (7,2) {\small{3}};

\node[] at (0,3) {\small{32}};
\node[] at (1,3) {\small{31}};
\node[] at (2,3) {\small{18}};
\node[] at (3,3) {\small{17}};
\node[] at (4,3) {\small{16}};
\node[] at (5,3) {\small{15}};
\node[] at (6,3) {\small{1}};
\node[] at (7,3) {\small{2}};

\node[] at (0,4) {\small{33}};
\node[] at (1,4) {\small{34}};
\node[] at (2,4) {\small{47}};
\node[] at (3,4) {\small{48}};
\node[] at (4,4) {\small{49}};
\node[] at (5,4) {\small{50}};
\node[] at (6,4) {\small{62}};
\node[] at (7,4) {\small{63}};

\node[] at (0,5) {\small{36}};
\node[] at (1,5) {\small{35}};
\node[] at (2,5) {\small{46}};
\node[] at (3,5) {\small{45}};
\node[] at (4,5) {\small{52}};
\node[] at (5,5) {\small{51}};
\node[] at (6,5) {\small{61}};
\node[] at (7,5) {\small{64}};

\node[] at (0,6) {\small{37}};
\node[] at (1,6) {\small{39}};
\node[] at (2,6) {\small{44}};
\node[] at (3,6) {\small{43}};
\node[] at (4,6) {\small{54}};
\node[] at (5,6) {\small{53}};
\node[] at (6,6) {\small{60}};
\node[] at (7,6) {\small{59}};

\node[] at (0,7) {\small{38}};
\node[] at (1,7) {\small{40}};
\node[] at (2,7) {\small{41}};
\node[] at (3,7) {\small{42}};
\node[] at (4,7) {\small{55}};
\node[] at (5,7) {\small{56}};
\node[] at (6,7) {\small{57}};
\node[] at (7,7) {\small{58}};

\draw[gray!70!white, ultra thick] (-0.5,-0.5) -- (7.5,-0.5);
\draw[gray!70!white, ultra thick] (-0.5, 0.5) -- (7.5, 0.5);
\draw[gray!70!white, ultra thick] (-0.5, 1.5) -- (7.5, 1.5);
\draw[gray!70!white, ultra thick] (-0.5, 2.5) -- (7.5, 2.5);
\draw[gray!70!white, ultra thick] (-0.5, 3.5) -- (7.5, 3.5);
\draw[gray!70!white, ultra thick] (-0.5, 4.5) -- (7.5, 4.5);
\draw[gray!70!white, ultra thick] (-0.5, 5.5) -- (7.5, 5.5);
\draw[gray!70!white, ultra thick] (-0.5, 6.5) -- (7.5, 6.5);
\draw[gray!70!white, ultra thick] (-0.5, 7.5) -- (7.5, 7.5);

\draw[gray!70!white, ultra thick] (-0.5, -0.5) -- (-0.5, 7.5);
\draw[gray!70!white, ultra thick] ( 0.5, -0.5) -- ( 0.5, 7.5);
\draw[gray!70!white, ultra thick] ( 1.5, -0.5) -- ( 1.5, 7.5);
\draw[gray!70!white, ultra thick] ( 2.5, -0.5) -- ( 2.5, 7.5);
\draw[gray!70!white, ultra thick] ( 3.5, -0.5) -- ( 3.5, 7.5);
\draw[gray!70!white, ultra thick] ( 4.5, -0.5) -- ( 4.5, 7.5);
\draw[gray!70!white, ultra thick] ( 5.5, -0.5) -- ( 5.5, 7.5);
\draw[gray!70!white, ultra thick] ( 6.5, -0.5) -- ( 6.5, 7.5);
\draw[gray!70!white, ultra thick] ( 7.5, -0.5) -- ( 7.5, 7.5);

\begin{scope}[shift={(10,0)}]

\draw[gray!70!white, ultra thick] (-0.5,-0.5) -- (7.5,-0.5);
\draw[gray!70!white, ultra thick] (-0.5, 0.5) -- (7.5, 0.5);
\draw[gray!70!white, ultra thick] (-0.5, 1.5) -- (7.5, 1.5);
\draw[gray!70!white, ultra thick] (-0.5, 2.5) -- (7.5, 2.5);
\draw[gray!70!white, ultra thick] (-0.5, 3.5) -- (7.5, 3.5);
\draw[gray!70!white, ultra thick] (-0.5, 4.5) -- (7.5, 4.5);
\draw[gray!70!white, ultra thick] (-0.5, 5.5) -- (7.5, 5.5);
\draw[gray!70!white, ultra thick] (-0.5, 6.5) -- (7.5, 6.5);
\draw[gray!70!white, ultra thick] (-0.5, 7.5) -- (7.5, 7.5);

\draw[gray!70!white, ultra thick] (-0.5, -0.5) -- (-0.5, 7.5);
\draw[gray!70!white, ultra thick] ( 0.5, -0.5) -- ( 0.5, 7.5);
\draw[gray!70!white, ultra thick] ( 1.5, -0.5) -- ( 1.5, 7.5);
\draw[gray!70!white, ultra thick] ( 2.5, -0.5) -- ( 2.5, 7.5);
\draw[gray!70!white, ultra thick] ( 3.5, -0.5) -- ( 3.5, 7.5);
\draw[gray!70!white, ultra thick] ( 4.5, -0.5) -- ( 4.5, 7.5);
\draw[gray!70!white, ultra thick] ( 5.5, -0.5) -- ( 5.5, 7.5);
\draw[gray!70!white, ultra thick] ( 6.5, -0.5) -- ( 6.5, 7.5);
\draw[gray!70!white, ultra thick] ( 7.5, -0.5) -- ( 7.5, 7.5);

\foreach \x in {0,...,7}
  \foreach \y [count=\yi] in {0,...,7}  
    \filldraw[color=green!60!black] (\x,\y) circle (3pt);

\draw[color=green!60!black, ultra thick] (6,3) -- (7,3) -- (7,2) -- (6,2) -- (6,1) -- (7,1) -- (7,0) -- (6,0);
\draw[color=green!60!black, ultra thick] (6,0) -- (5,0) -- (4,0) -- (4,1) -- (5,1) -- (5,2) -- (4,2) -- (5,3) -- (4,3);
\draw[color=green!60!black, ultra thick] (4,3) -- (3,3) -- (2,3) -- (2,2) -- (3,2) -- (3,1) -- (3,0) -- (2,0) -- (2,1);
\draw[color=green!60!black, ultra thick] (2,1) -- (1,1) -- (1,0) -- (0,0) -- (0,1) -- (0,2) -- (1,2) -- (1,3) -- (0,3);
\draw[color=green!60!black, ultra thick] (0,3) -- (0,4) -- (1,4) -- (1,5) -- (0,5) -- (0,6) -- (0,7) -- (1,6) -- (1,7);
\draw[color=green!60!black, ultra thick] (1,7) -- (2,7) -- (3,7) -- (3,6) -- (2,6) -- (3,5) -- (2,5) -- (2,4) -- (3,4);
\draw[color=green!60!black, ultra thick] (3,4) -- (4,4) -- (5,4) -- (5,5) -- (4,5) -- (5,6) -- (4,6) -- (4,7) -- (5,7);
\draw[color=green!60!black, ultra thick] (5,7) -- (6,7) -- (7,7) -- (7,6) -- (6,6) -- (6,5) -- (6,4) -- (7,4) -- (7,5);

\end{scope}

\end{tikzpicture}
%\end{document}
\caption{The space-filling curve approach applied to a structured 8~by~8 grid in 2D. Left: the resulting vertex numbers or space-filling curve ordering, right: the path of the space-filling curve shown in green. The underlying grid is shown in grey.}
\label{unstructured_method_regular} 
\end{minipage}\hspace{0.05\textwidth}
\begin{minipage}{0.3\textwidth}
\centering
\begin{tikzpicture}[scale=1.2, rotate=270] % 2

\coordinate (A1) at (0, 0+0.1);
\coordinate (A2) at (-0.5, 0.85+0.1);
\coordinate (A3) at ( 0.5, 0.85+0.1);

\draw[green!60!black, ultra thick] (0,0.35) circle (1.5pt);         %(:,:+0.25)
\filldraw[green!60!black] ( 0.5-0.2, 0.85+0.1-0.12) circle (1.5pt); %(:-0.2,:-0.12)
\filldraw[green!60!black] (-0.5+0.2, 0.85+0.1-0.12) circle (1.5pt); %(:+0.2,:-0.12)

\coordinate (B1) at (0-0.1, 0);
\coordinate (B2) at (-0.5-0.1, 0.85); 
\coordinate (B3) at (-1-0.1,0);

\filldraw[green!60!black] (0-0.1-0.2, 0+0.12) circle (1.5pt);         
\filldraw[green!60!black] (-0.5-0.1, 0.85-0.25) circle (1.5pt); % (-0.25,:)
\filldraw[green!60!black] (-1-0.1+0.2,0+0.12) circle (1.5pt);

\coordinate (C1) at (0+0.1, 0);
\coordinate (C2) at (0.5+0.1, 0.85);
\coordinate (C3) at (1+0.1,0);

\filldraw[green!60!black] (0+0.1+0.2, 0+0.12) circle (1.5pt);         
\filldraw[green!60!black] ( 0.5+0.1, 0.85-0.25) circle (1.5pt); % (-0.25,:)
\filldraw[green!60!black] (1+0.1-0.2,0+0.12) circle (1.5pt);

\draw[gray!60!, ultra thick] (A1) -- (A2) -- (A3) -- cycle;
\draw[gray!60!, ultra thick] (B1) -- (B2) -- (B3) -- cycle;
\draw[gray!60!, ultra thick] (C1) -- (C2) -- (C3) -- cycle;

\coordinate (D1) at (0, 0+0.85+0.85+0.1+0.1);
\coordinate (D2) at (-0.5, 0.85+0.1+0.1);
\coordinate (D3) at (0.5, 0.85+0.1+0.1);

\filldraw[green!60!black] (0, 0+0.85+0.85+0.1+0.11-0.25) circle (1.5pt);         
\filldraw[green!60!black] (-0.5+0.2, 0.85+0.1+0.11+0.12)     circle (1.5pt); % (-0.25,:)
\filldraw[green!60!black] ( 0.5-0.2, 0.85+0.1+0.11+0.12) circle (1.5pt);

\draw[gray!60!, ultra thick] (D1) -- (D2) -- (D3) -- cycle;

\coordinate (E1) at (0, 0-0.1-0.1);
\coordinate (E2) at (-0.5, -0.85-0.1-0.1);
\coordinate (E3) at (0.5, -0.85-0.1-0.1);

\filldraw[gray!60!] (0, 0-0.1-0.1-0.25)             circle (1.5pt); % (-0.25,:)        
\filldraw[gray!60!] (-0.5+0.2, -0.85-0.1-0.1+0.12)  circle (1.5pt);  %(:+0.2,:+0.12)
\filldraw[gray!60!] ( 0.5-0.2, -0.85-0.1-0.1+0.12)  circle (1.5pt);  %(:-0.2,:+0.12)

\coordinate (F1) at (0-0.1, 0-0.1);
\coordinate (F2) at (-0.5-0.1, -0.85-0.1);
\coordinate (F3) at (-1-0.1,0-0.1);

\filldraw[gray!60!] (0-0.1-0.2, 0-0.1-0.12)     circle (1.5pt); % (:+0.2,:+0.12)
\filldraw[gray!60!] (-0.5-0.1, -0.85-0.1+0.25)  circle (1.5pt);  % (+0.25,:) 
\filldraw[gray!60!] (-1-0.1+0.2,0-0.1-0.12)    circle (1.5pt);  %(:-0.2,:+0.12)

\coordinate (G1) at (0+0.1, 0-0.1);
\coordinate (G2) at (0.5+0.1, -0.85-0.1);
\coordinate (G3) at (1+0.1,0-0.1);

\filldraw[gray!60!] (0+0.1+0.2, 0-0.1-0.12)     circle (1.5pt); % (:+0.2,:+0.12)
\filldraw[gray!60!] (0.5+0.1, -0.85-0.1+0.25)  circle (1.5pt);  % (+0.25,:) 
\filldraw[gray!60!] (1+0.1-0.2,0-0.1-0.12)    circle (1.5pt);  %(:-0.2,:+0.12)

\draw[gray!60!, ultra thick] (E1) -- (E2) -- (E3) -- cycle;
\draw[gray!60!, ultra thick] (F1) -- (F2) -- (F3) -- cycle;
\draw[gray!60!, ultra thick] (G1) -- (G2) -- (G3) -- cycle;

\coordinate (I1) at (0-0.1, 0-0.1+1.7+0.2+0.2);
\coordinate (I2) at (-0.5-0.1, -0.85-0.1+1.7+0.2+0.2);
\coordinate (I3) at (-1-0.1,0-0.1+1.7+0.2+0.2);

\draw[gray!60!, ultra thick] (I1) -- (I2) -- (I3) -- cycle;

\filldraw[gray!60!] (0-0.1-0.2, 0-0.1-0.12+1.7+0.2+0.2)     circle (1.5pt); % (:+0.2,:+0.12)
\filldraw[gray!60!] (-0.5-0.1, -0.85-0.1+0.25+1.7+0.2+0.2)  circle (1.5pt);  % (+0.25,:) 
\filldraw[gray!60!] (-1-0.1+0.2,0-0.1-0.12+1.7+0.2+0.2)    circle (1.5pt);  %(:-0.2,:+0.12)

\coordinate (J1) at (0+0.1, 0-0.1+1.7+0.2+0.2);
\coordinate (J2) at (0.5+0.1, -0.85-0.1+1.7+0.2+0.2);
\coordinate (J3) at (1+0.1,0-0.1+1.7+0.2+0.2);

\draw[gray!60!, ultra thick] (J1) -- (J2) -- (J3) -- cycle;

\filldraw[gray!60!] (0+0.1+0.2, 0-0.1-0.12+1.7+0.2+0.2)     circle (1.5pt); % (:+0.2,:+0.12)
\filldraw[gray!60!] (0.5+0.1, -0.85-0.1+0.25+1.7+0.2+0.2)  circle (1.5pt);  % (+0.25,:) 
\filldraw[gray!60!] (1+0.1-0.2,0-0.1-0.12+1.7+0.2+0.2)    circle (1.5pt);  %(:-0.2,:+0.12)

\coordinate (K1) at (0-1.2, 0+0.1);
\coordinate (K2) at (-0.5-1.2, 0.85+0.1);
\coordinate (K3) at ( 0.5-1.2, 0.85+0.1);

\draw[gray!60!, ultra thick] (K1) -- (K2) -- (K3) -- cycle;

\filldraw[gray!60!] (0-1.2,0.35) circle (1.5pt);         %(:,:+0.25)
\filldraw[gray!60!] ( 0.5-0.2-1.2, 0.85+0.1-0.12) circle (1.5pt); %(:-0.2,:-0.12)
\filldraw[gray!60!] (-0.5+0.2-1.2, 0.85+0.1-0.12) circle (1.5pt); %(:+0.2,:-0.12)

\coordinate (L1) at (0+1.2, 0+0.1);
\coordinate (L2) at (-0.5+1.2, 0.85+0.1);
\coordinate (L3) at ( 0.5+1.2, 0.85+0.1);

\draw[gray!60!, ultra thick] (L1) -- (L2) -- (L3) -- cycle;

\filldraw[gray!60!] (0+1.2,0.35) circle (1.5pt);         %(:,:+0.25)
\filldraw[gray!60!] ( 0.5-0.2+1.2, 0.85+0.1-0.12) circle (1.5pt); %(:-0.2,:-0.12)
\filldraw[gray!60!] (-0.5+0.2+1.2, 0.85+0.1-0.12) circle (1.5pt); %(:+0.2,:-0.12)

\coordinate (M1) at (0-1.2, 0+0.1+0.85+0.85+0.1);
\coordinate (M2) at (-0.5-1.2, 0.85+0.1+0.1);
\coordinate (M3) at ( 0.5-1.2, 0.85+0.1+0.1);

\draw[gray!60!, ultra thick] (M1) -- (M2) -- (M3) -- cycle;

\filldraw[gray!60!] (0-1.2,0.35+0.85+0.85-0.4) circle (1.5pt);         %(:,:+0.25)
\filldraw[gray!60!] ( 0.5-0.2-1.2, 0.85+0.1-0.12+0.32) circle (1.5pt); %(:-0.2,:-0.12)
\filldraw[gray!60!] (-0.5+0.2-1.2, 0.85+0.1-0.12+0.32) circle (1.5pt); %(:+0.2,:-0.12)

\coordinate (N1) at (0+1.2, 0+0.1+0.85+0.85+0.1);
\coordinate (N2) at (-0.5+1.2, 0.85+0.1+0.1);
\coordinate (N3) at ( 0.5+1.2, 0.85+0.1+0.1);

\draw[gray!60!, ultra thick] (N1) -- (N2) -- (N3) -- cycle;

\filldraw[gray!60!] (0+1.2,0.35+0.85+0.85-0.4) circle (1.5pt);         %(:,:+0.25)
\filldraw[gray!60!] ( 0.5-0.2+1.2, 0.85+0.1-0.12+0.32) circle (1.5pt); %(:-0.2,:-0.12)
\filldraw[gray!60!] (-0.5+0.2+1.2, 0.85+0.1-0.12+0.32) circle (1.5pt); %(:+0.2,:-0.12)

\end{tikzpicture}
\caption{A group of linear finite elements with a discontinuous Galerkin stencil shown for one node (hollow green circle). The stencil of this node connects it to all the other green nodes.}
\label{fig:DGstencil}
\end{minipage}
\end{figure}

\begin{figure}[htbp]
     \centering
\hbox{ 
\hspace{-1.cm} 
      \includegraphics[width=9.5cm, trim=0mm 0mm 0mm 0mm, clip]{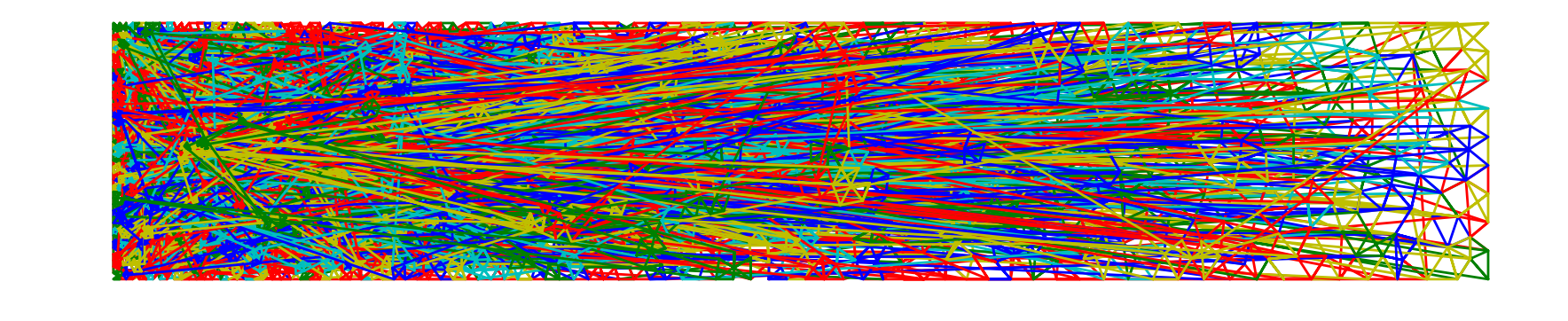}
\hspace{-0.5cm} 
      \includegraphics[width=9.5cm, trim=0mm 0mm 0mm 0mm, clip]{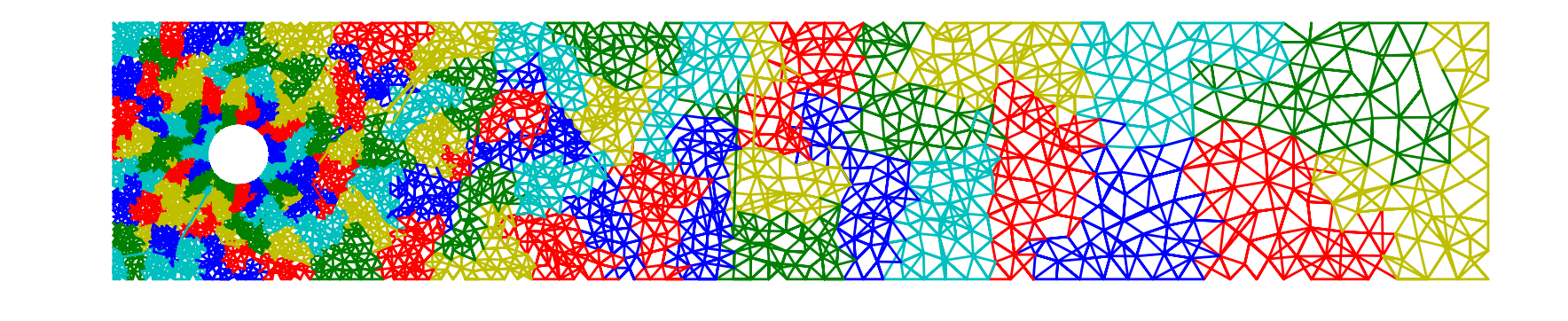}
}
\hbox{ 
\hspace{-0.6cm} 
      \includegraphics[width=8.7cm, trim=0mm 0mm 0mm 0mm, clip]{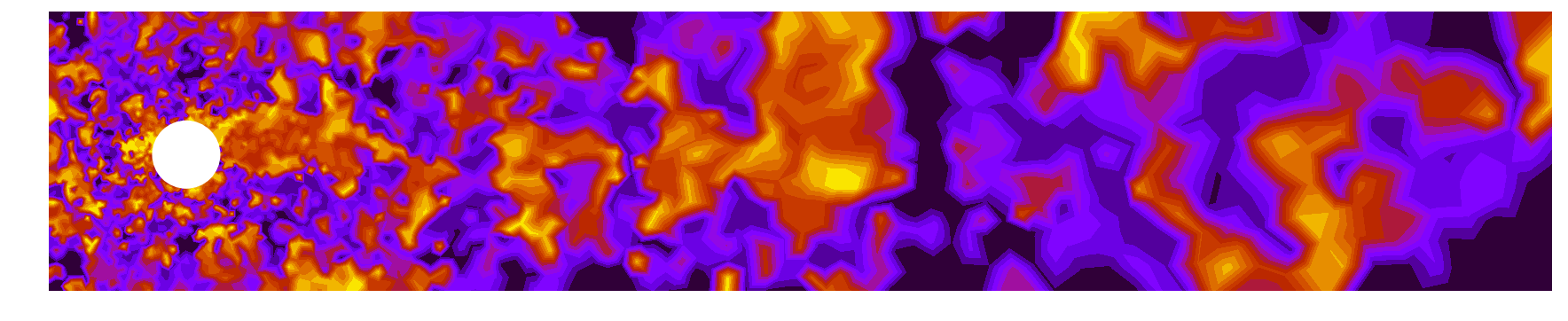}
\hspace{0.3cm} 
      \includegraphics[width=8.7cm, trim=0mm 0mm 0mm 0mm, clip]{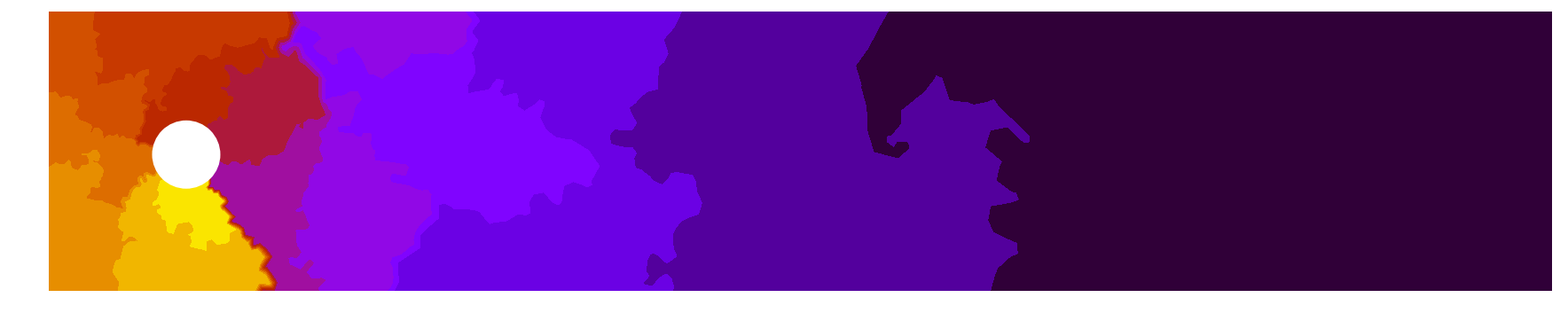}
}
\caption{The ordering based on a space-filling curve ordering applied to an unstructured mesh and DG discretisation for flow past a cylinder. Top two plots show the lines drawn between nodes with consecutive ordering: with the original FE ordering as determined by the meshing software (left) and with the SFC ordering (right). Notice that nodes that are close to one another are physically close also. The path between each group of 200 consecutive nodes are shown in same colour which forms a series of subdomains because of the properties of SFCs. Bottom two pictures: contour plot of the node numbers according the FEM ordering (left) and SFC ordering (right).   }
\label{2curves_2d_fem_cylinder_space_filling} 
\end{figure}

The discontinuous Galerkin (DG) discretisation has a particularly complex stencil (and graph) associated with it, see Figure~\ref{fig:DGstencil} for the stencil associated with one node. For a DG node in element~$k$, this node is connected (through the stencil) to all nodes in element~$k$ and all nodes in elements which share an edge (2D) or face (3D) with element~$k$. The edges for the node shown in Figure~\ref{fig:DGstencil} would be constructed by drawing a line from the this node to each other green node. To draw the graph for the entire mesh shown here, this process would need to be done for every node in the mesh. The result of applying this method to an unstructured mesh can be seen in Figure~\ref{2curves_2d_fem_cylinder_space_filling}, where a linear DG finite element stencil is used. There are 20,550~DG nodes in this mesh and 6,850~elements. In the top plots of this figure (top left and right), we show the path linking the first 200~DG nodes in one colour. The path connecting the next 200~DG nodes are plotted with another colour and so on. The plot on the left uses the node ordering as dictated by the meshing software, the plot on the right uses the node ordering provided by the space-filling curve. The ability of the space-filling curve to group together neighbouring nodes is clearly demonstrated. By grouping 200 DG nodes together in this way, a coarse-mesh representation of the problem is formed, illustrating one of the key reasons why combining space-filling curves and convolutional networks is so effective. The top right plot in Figure~\ref{2curves_2d_fem_cylinder_space_filling} could be interpreted as coarsening of the finite element mesh. Suppose that one level of coarsening was applied by taking every 4th node along the space-filling curve (reducing by a factor of~2 in each dimension), and this was done 4 times, then this is equivalent to taking every 256th node. %Therefore, one could say the coarsening shown in this plot is similar to applying 4 levels of coarsening (as 200$\sim$256), just as one could do if applying convolutional layers with a stride of 2 for each dimension.
Therefore, one could say the coarsening shown in this plot, by colouring groups of 200~DG nodes, is similar to applying 4 levels of coarsening which reduces the resolution by a factor~2 in each dimension (as 200$\sim$256), which is similar to applying 4 convolutional layers with a stride of 2 for each dimension. The lower plots in Figure~\ref{2curves_2d_fem_cylinder_space_filling} show contour plots of the node numbers, with low numbers in black through to the highest node numbers in yellow. The plot on the left corresponds to the default node numbering based on the numbering provided by the meshing software, and the plot on the right corresponds to the node numbering from the space-filling curve. 

%\clearpage
\subsection{Partitioning method} \label{sec:SFC_partitioning}
Many problems can be abstracted to the partitioning of a graph, such as the decomposition of computational domains, social network problems and travel network problems. Therefore much effort has been applied to the fundamental challenge of developing a method which simultaneously achieves an equal distribution of vertices/load between the partitions whilst minimising communication between partitions. A comprehensive survey of methods can be found here~\cite{Buluc2016}, which include spectral partitioning, geometric partitioning~\cite{Fitzgerald2019}, multi-grid methods~\cite{Dargaville_2020},  methods based on the mean-field theorem~\cite{Kawamoto2019}, and also nested bisection approaches~\cite{Karypis1998_JPDC, Karypis1998_SIAM, Moulitsas2001}, which have seen success for partitioning large graphs.

%Here we choose a recurrent neural network based on mean-field theory~\cite{Pain_1999}, which is able to generate partitions that have equal numbers of vertices (as far as possible) and that can use variable edge weights to influence the path, useful when generating multiple space-filling curves. Other methods could be used such as those based on agglomeration~\cite{Dargaville_2020}
A recurrent neural network based on mean-field theory (MFT-RNN)~\cite{Pain_1999} is chosen to be the graph partitioner for the current work, as it (i)~has a facility to set variable edge weights in the graph and (ii)~can balance  exactly (as far as possible) the number of vertices in each partition at each level. When generating multiple space-filling curves, there is a need to avoid cutting certain edges in the graph, so that different space-filling curves traverse distinct paths through the graph; which leads to the first requirement. The second requirement arises from the desire to obtain a final set of partitions, of which as many as possible  contain at most one vertex. Partitions with more than one vertex can be dealt with, as can empty partitions, but the performance of the space-filling curve numbering algorithm will suffer in the former case. %
%A recurrent neural network is combined with mean-field theory which allows perfect balancing of the numbers of vertices in each sub-graph (to within one vertex) and edge weights can be easily included in the formulation. Other graph partitioning methods could be used as long as they have properties~(i) and~(ii). 
This multi-grid partitioning approach forms a hierarchy or nested sequence of partitions akin to that seen in multi-grid finite element methods~\cite{Pain_1999,Dargaville_2020}. Other graph partitioning methods could be used as long as properties~(i) and~(ii) are attainable, such as methods that use vertex agglomeration by repeatedly merging adjacent vertices in a similar approach to multi-grid agglomeration methods~\cite{Dargaville_2020}.

More details of the MFT-RNN partitioner can be found in~\cite{Pain_1999}, however the main points are given here, along with details of the modifications made when generating multiple space-filling curves.  The number of partitions created at each level is represented by~$\mathcal{S}$, which, for the particular fundamental shape used in this paper, is necessarily ${\mathcal{S}=2}$.  This recurrent neural network associates $\mathcal{S}$ neurons with each vertex in the graph (i.e.~one neuron for each desired partition), whose values represent the probability that a vertex~$i$ is in partition~$\mu$:
\begin{equation}
    z^\mu_i \in [0,1]\ , \qquad \text{subject to the constraint } \sum_{\mu=1}^{\mathcal{S}} z^\mu_i = 1\quad \forall i\ .
\end{equation}
The constraint has led to this method being referred to as normalised mean-field theory. For each partition, a vector of probabilities can be defined as 
\begin{equation}
    \bm{z}^{\mu} = \left( z^\mu_1, z^\mu_2, \ldots ,z^\mu_i, \ldots, z^\mu_{N^{\text{vert}}} \right)^T \qquad \forall \mu\in\{1,2,\ldots,\mathcal{S}\}\ ,
\end{equation}
where $N^{\text{vert}}$ is the number of vertices in the graph. For the special case of $\mathcal{S}=2$, only one neuron need be included for every vertex, as, given the probability of the vertex being in one partition, the probability of being the other can be calculated directly from the constraint. However, the code used to generate the examples in this paper was written only for the general case. Representing the quality of the partitions, the functional comprises two terms: the first describes the communication cost between vertices in different (proposed) partitions and the second describes how well the vertices are distributed between the partitions. It is written as follows:
\begin{equation}\label{gp4}
F=\tfrac{1}{2}\,\bm{z}^T\bm{K}\bm{z} + \tfrac{1}{2} \, \alpha \bm{z}^T\bm{C}\bm{z} \,, 
\end{equation}
where $\bm{z}^T = \left( (\bm{z}^{1})^T, (\bm{z}^{2})^T,\ldots,(\bm{z}^{\mu})^T,\ldots,(\bm{z}^{\mathcal{S}})^T \right)$, and the scalar~$\alpha$ controls the balance between the importance of the distribution of vertices throughout the partitions and the importance of reducing the communication costs. The matrices~$\bm{K}$ and~$\bm{C}$ represent communication costs and balancing respectively, and are now explained in detail. The matrix $\bm{K}$ consists of sub-matrices $\bm{K}^{\mu\nu}$ for partition numbers~$\mu$ and~$\nu$, where  ${\mu,\nu\in\{1,2,\ldots,\mathcal{S}\}}$. For $\mu=\nu$, the sub-matrices consist entirely of zeros. For $\mu\neq\nu$, the entries of the sub-matrices are defined as
\begin{equation}\label{gp3----}
K_{ij}^{\mu\nu} = \left\{ \begin{array}{l} h_{ij} \quad \text{if $i\neq j$ and vertices~$i$ and~$j$ are adjacent in the graph,} \\ 0 \quad \text{otherwise.}\end{array} \right. 
\end{equation}
The value $h_{ij}$ represents the cost of vertex~$i$ communicating with vertex~$j$ and as well as determining the communication cost in the functional, this value is also taken as the $ij$th edge weight in the graph. When generating one space-filling curve, ${h_{ij}=1\ \forall i,j}$. When generating more than one space-filling curve, this parameter is set as
\begin{equation}
    h_{ij} = \left\{ \begin{array}{lll} \displaystyle{\max_{m}} \,\lvert s^m_i - s^m_j\rvert^\gamma & \qquad & \text{if vertices~$i$ and~$j$ are adjacent} \\ 0 & & \text{otherwise},
    \end{array}\right.
\end{equation}
where the heuristically obtained exponent~$\gamma$ is set to be~0.2 and $s^m_i$ is the new space-filling curve node numbering  %\textcolor{red}{node number \textit{(Chris check this)}} 
of the $i$th vertex in any one of the pre-existing space-filling curves labelled by~$m$. The idea behind introducing variable edge weights to the graph is to inject a means with which to discourage the graph to partition across edges that have large differences in space-filling curve numbers. This helps to ensure that a subsequent space-filling curve will avoid edges traversed by pre-existing curves, and encourages subsequent space-filling curves to find previously undiscovered paths through the graph. So, if two points are far apart on the first space-filling curve but close in Cartesian space, a second space-filling curve will attempt to find a path on which these points are closer.

Similar to~$\bm{K}$, the matrix $\bm{C}$ consists of sub-matrices $\bm{C}^{\mu\nu}$ for partition numbers~$\mu$ and~$\nu$, where ${\mu,\nu\in\{1,2,\ldots,\mathcal{S}\}}$. To balance the vertices equally between the partitions (as far as possible), a diffusion matrix is used, the entries of which are defined as 
\begin{equation}
  C^{\mu\nu}_{ij} = \left\{ \begin{array}{cll} \dfrac{\mathcal{S}-1}{\mathcal{S}} & \qquad & \text{if $i= j$}\\[3mm] \dfrac{-1\ }{\mathcal{S}}  & & \text{if $i\neq j$} \end{array}\right.    
\end{equation}
for all partition numbers $\mu,\nu\in\{1, 2, \ldots, \mathcal{S}\}$ and for all vertices~$i,j\in\{1,2,\ldots,N^{\text{vert}}\}$. In order to balance the number of vertices in each of the partitions at a particular level, the probabilities $z_i^\mu$ of vertex~$i$ belonging to partition~$\mu$ are exploited. If there is a need to transfer some vertices from one partition to another, the most suitable neurons are chosen from other partitions according to their probabilities. In this way, the uncertainty near the partition boundaries of the MFT-RNN is exploited, and thus, good quality partitions are maintained as well as near perfect balance of vertices in each partition.

The two terms in the functional tends to compete against one another during optimisation, so it is generally advantageous to set~$\alpha$ as small as possible, although just large enough so that the vertices are equally distributed between the partitions. If the value of~$\alpha$ is too large, the partitions will have equal numbers of vertices but the path traced out by the space-filling curve will not be as continuous as it could have been, resulting in a poor space-filling curve. In this paper, an optimised value of~$\alpha$ is used as derived in~\cite{Pain_1999}: 
\begin{equation}
\alpha=\frac{1}{\mathcal{S} N^{\text{vert}}}  \sum_{j=1}^{N^{\text{vert}}} \sum_{i=1}^{N^{\text{vert}}} K_{ij} \  . 
\label{gp14}
\end{equation}
Extensive numerical experiments have shown that the value of~$\alpha$ obtained from this equation    results in good load balancing as well as good quality partitioning that minimises the sum of the weights between the different partitions.   

%\textcolor{red}{\textit{Chris:} General procedure and updating - to be completed}\\
In summary, the graph partitioning algorithm takes a graph and decomposes this into partitions, using a nested bisection approach. This continues until almost all of the partitions contain one vertex. This is done by calculating the number of levels of decomposition, $L$, by finding the smallest integer~$L$ that satisfies
\begin{equation}
2^{L} \geqslant \text{the total number of vertices in the original graph.}
\end{equation}
The functional, Equation~\eqref{gp4}, is minimised using the neural network updating algorithm associated with the MFT neural network as described in~\cite{Pain_1999} to obtain the desired graph partitioning. No training of this network is required as the weights are determined by the functional. %
%  Empty partitions are allowed, as are partitions with more then one vertex. 
The space-filling curve can now be constructed from the numbered partitions.

\subsection{Multiple space-filling curves} \label{sec:SFC_multiple}

\begin{figure}[htbp]
    \centering
    \begin{tikzpicture}[scale=1.1]%[scale=0.9]
\usetikzlibrary{positioning}
% sfc 1 -------------------------------------------------------------------------------

\node[] at (2,-0.5) {\small (a)};
\draw[step=1cm,darkgray,thin] (0,0) grid (4,4);

\foreach \x in {0,...,3}
  \foreach \y [count=\yi] in {0,...,3}  
    \filldraw[color=cyan!75!black] (0.5+\x,0.5+\y) circle (3pt);
\draw[cyan!75!black,very thick] (0.5,2.5) -- (0.5,3.5) -- (1.5,3.5)  -- (1.5,2.5)   -- (1.5,1.5) -- (0.5,1.5) -- (0.5,0.5) -- (1.5,0.5) -- (2.5,0.5) -- (3.5,0.5) -- (3.5,1.5) -- (2.5,1.5) -- (2.5,2.5) -- (3.5,2.5) -- (3.5,3.5) -- (2.5,3.5) ;

% sfc 2  --------------------------------------------------------------------
\node[] at (7,-0.5) {\small (b)};

\draw[step=1cm,darkgray,thin] (5,0) grid (9,4);

\foreach \x in {0,...,3}
  \foreach \y [count=\yi] in {0,...,3}  
    \filldraw[color=orange] (5.5+\x,0.5+\y) circle (3pt);
%\draw[green!60!black, very thick] (3.5,0.5) -- (3.5,3.5);
\draw[orange,very thick] (8.5,3.5) .. controls (8.85,2) .. (8.5,0.5) -- (8.5,1.5) -- (8.5,2.5) -- (5.5,0.5) .. controls (5.15,2) .. (5.5,3.5) -- (6.5,0.5) -- (7.5,0.5) -- (7.5,1.5) -- (6.5,1.5) -- (5.5,1.5) -- (5.5,2.5) -- (6.5,2.5) -- (7.5,2.5) -- (7.5,3.5) --(6.5,3.5);

% sfc 1 and 2  -------------------------------------------------------------------
\node[] at (12,-0.5) {\small (c)};
\draw[step=1cm,gray,thin] (10,0) grid (14,4);

\draw[gray!90!black,very thick] (10.5,2.5) -- (10.5,3.5) -- (11.5,3.5)  -- (11.5,2.5)   -- (11.5,1.5) -- (10.5,1.5) -- (10.5,0.5) -- (11.5,0.5) -- (12.5,0.5) -- (13.5,0.5) -- (13.5,1.5) -- (12.5,1.5) -- (12.5,2.5) -- (13.5,2.5) -- (13.5,3.5) -- (12.5,3.5) ;

\draw[gray!90!black,very thick] (13.5,3.5) .. controls (13.85,2) .. (13.5,0.5) -- (13.5,1.5) -- (13.5,2.5) -- (10.5,0.5) .. controls (10.15,2) .. (10.5,3.5) -- (11.5,0.5) -- (12.5,0.5) -- (12.5,1.5) -- (11.5,1.5) -- (10.5,1.5) -- (10.5,2.5) -- (11.5,2.5) -- (12.5,2.5) -- (12.5,3.5) --(11.5,3.5);

\draw[red,ultra thick] (10.5,1.5) -- (11.5,1.5);
\draw[red,ultra thick] (11.5,0.5) -- (12.5,0.5);
\draw[red,ultra thick] (13.5,0.5) -- (13.5,1.5);

\draw[blue!80!cyan,ultra thick] (11.5,0.5) -- (11.5,1.5);

\foreach \x in {0,...,3}
  \foreach \y [count=\yi] in {0,...,3}  
    \filldraw[color=gray!90!black] (10.5+\x,0.5+\y) circle (3pt);

\end{tikzpicture}
\caption{Here we show the paths traversed by the first~(a) and second~(b) space-filling curves for a ${4\times 4}$ uniform grid with a 5~point stencil. In the plot on the right, we superpose both curves (grey). Red indicates edges that appear in both space-filling curves and blue indicates the only edge that is in neither of the space-filling curves. 
}
\label{2SFCs-structured} 
\end{figure}
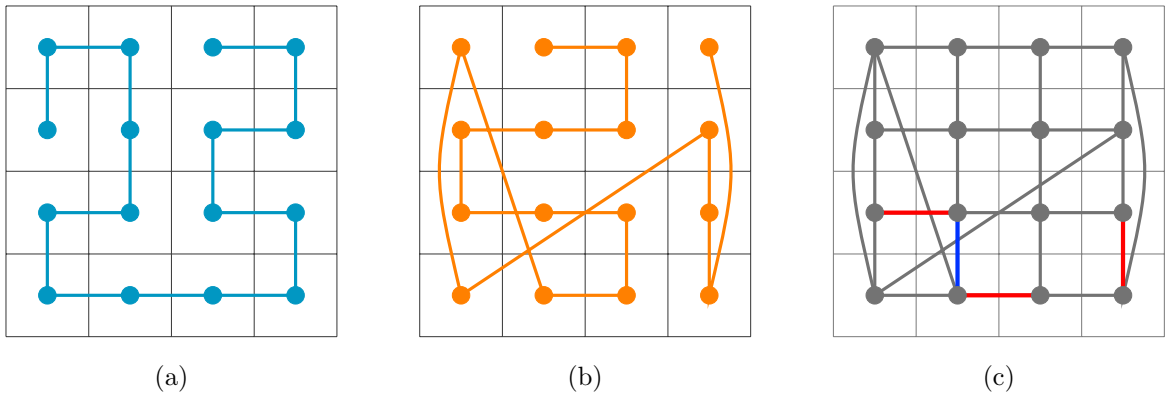

One problem with using the space-filling curve approach is that there can be little connectivity between certain regions of the domain, because although points close together on the space-filling curve are also close together in the physical space, the reverse is not true. A simple example of this can be seen in Figure~\ref{fig:Hilbert_curve_level_1_and_2}(d), where vertices~2 and~15 are not close on the space-filling curve but are close in the original graph in~(e). % A simple example of this is show in figure \ref{2SFCs-structured-start-finish-the-same} (left which shows a SFC for a $4\times 4$ regular grid) in which node number $3$ in the space filling curve ordering is next to node number $16$ -  near the maximum node number difference in the space filling curve  which is $15$. 
A possible remedy is to use more than one space-filling curve, where the second curve should be discouraged from having the same edges as the first. %
%To achieve this goal in the first approach these edges are deleted from the graph that the secondary SFC uses. If 3 space filling curves are needed (e.g. for a 3D unstructured mesh) then one excludes both the primary and secondary SFC edges in the graph that is used to produce the tertiary SFC and so on if we want more than 3 SFCs. We have applied this approach to a regular $4\times 4$ grid with a 5 point stencil defining the edges of the graph and show the primary and figure \ref{2SFCs-structured}. Notice that an optimal SFC is produced for the primary SFC. This is not always true but it will be near optimal in some sense, as seen in figure  \ref{unstructured_method_regular}. The secondary SFC shown in the bottom-middle of the figure  \ref{2SFCs-structured} has some long distance links between the nodes.  This is partly because the 4 corner nodes of the $4\times 4$ grid don't have any neighbours in the graph that is used to form the secondary SFC. The first and second SFC are shown on the same plot in the bottom right of figure \ref{2SFCs-structured}. Notice that only one edge has neither the primary nor secondary SFC going through it. That is almost all edges will be traversed by either the primary or secondary SFC, in this case. 
One way of achieving this is though the graph partitioning 
procedure, by weighting the graph edges as described in Equation~\eqref{gp3----}. This effectively discourages partitions from being made across the same edges in the previous SFCs. Nothing else is needed as the edge weights determine the communication between areas of the domain in the any additional SFCs. An example of two space-filling curves is shown in Figure~\ref{2SFCs-structured}. The first SFC (left) is optimal in that it uses the minimum number of edges to traverse all the cells in the grid. The second SFC (centre) is less optimal, as it jumps to nodes that are not neighbours in the stencil, shown as curved or diagonal lines in~(b). In order to do this, the SFC will pass through nodes and edges more than once. This is partly due to the small size of the example. However, from~(c), the combination of the two space-filling curves is shown, which shows that all but one of the edges on the original graph (shown in Figure~\ref{fig:Hilbert_curve_level_1_and_2}(e)) are connected and only three edges appear in both SFCs.

%\begin{figure}[htbp]
%     \centering
%%    \hbox{\centering
%%      \hspace{-1.5cm} 
%%      \includegraphics[width=7.5cm]{./diagrams/pic-conv1.pdf}
%%      \hspace{1.75cm} 
%%      \includegraphics[width=7.5cm]{./diagrams/pic-conv2.pdf}
%%    }
%    \hbox{ \hspace{1.5cm} %\centering
%      \includegraphics[width=5cm, trim = 60mm 0mm 50mm 0mm, clip]{./diagrams/cnn_paper_pics_p7.pdf}
%      \hspace{-0.5cm} 
%      \includegraphics[width=5cm, trim = 60mm 0mm 50mm 0mm, clip]{./diagrams/cnn_paper_pics_p8.pdf}
%      \hspace{-0.5cm} 
%      \includegraphics[width=5cm, trim = 60mm 0mm 50mm 0mm, clip]{./diagrams/cnn_paper_pics_p9.pdf}
%  }
%\caption{Here we show the paths traversed by the first and second space-filling curves for a ${4\times 4}$ uniform grid with a 5~point stencil or graph. The leftmost plot shows the first space-filling curve, the central plot shows the second SFC, and on the right we superpose both curves. In the right plot the black lines show the first SFC, the red lines the second SFC, the yellow edges indicate edges that are both on the first and second SFC, and the black dotted edge is the edge that is in neither the primary nor the secondary SFC. %\textcolor{red}{Claire: redraw this?}
%}
%\label{2SFCs-structured} 
%\end{figure}

%The primary SFC node/cell ordering are shown in figures  \ref{2curves_2d_fem_cylinder_space_filling}. 
%These show that there is significant overlap between the SFC grouping of 200 FEM nodes 
%which helps to evaluate the features right across the domain and produce good connectivity of the 
%information captured by the CNN. 

We have experimented with an alternative method, which deletes the edges of the first SFC from the original graph used by the graph partitioner. Subsequent SFCs will then be forced to find new paths through the mesh. With this method, the quality of the second and subsequent SFCs was not as good as with the previous method, possibly because it excessively reduced the number of graph edges available to subsequent SFCs. 
%Hard constraint.

We have also experimented with introducing a constraint on the second SFC so that this space-filling curve starts at the final node of the first SFC. In this way, the second SFC follows the first SFC resulting in a  continuous curve, which is advantageous to the performance of the 1D CNN. This constraint could be introduced in the SFC algorithm, simply, by always ensuring that the first partition on each level of the nested bisection contains this node (that is, the end node of the first SFC). For three space-filling curves one adds a similar constraint to the third SFC so that it start at the node next to the final node of the second SFC. %
%See figure \ref{2SFCs-structured-start-finish-the-same} for an example of the use of this. 
A possible disadvantage of using this start-end constraint is that it can compromise the quality of the second and subsequent SFCs, in terms of producing a continuous curve, as, once again, the options of the SFC algorithm are restricted. A second disadvantage is that the number of weights is greater when the SFCs are combined into one continuous curve. For a 2D data set with 2 SFCs, having two separate SFCs results in a factor of two saving, in terms of weights in the 1D convolutional layers for the same total number of channels. A similar saving is also obtained in terms of arithmetic operations. For 3D problems this approach results in a saving of 3~times the number of weights and approximately similar saving for the CPU requirements associated with the convolutional layers. Thus, we chose here to have separate SFCs in the SFC-based CNNs.

\section{Architectures of the convolutional networks}\label{sec:architectures}
The architectures of the CNNs that will be used in the results section are described in detail in this section. To demonstrate the application of convolutional networks to data on both structured and unstructured meshes, a convolutional autoencoder is chosen. This type of network attempts to learn a compressed representation of data due to a bottleneck as its central layer. Consisting of an encoder and a decoder, the encoder compresses the data to a predetermined number of latent variables, known as the dimension of the latent space. The decoder reconstructs the data from the latent variables. For a convolutional autoencoder, the encoder has a number of convolutional layers typically followed by some fully-connected layers, and the decoder is the reverse of this. %Consisting of an encoder for the compression and a decoder for the reconstruction, the autoencoder lends itself to comparison with singular value decomposition, another method that is often used for compressing data. %\textcolor{red}{Comment on not being as good for picture compression as jpg. However generative aspects of VAE, removal of noise, and dimensionality reduction are interesting properties and focus of research.}

In addition to the layers found in a classical convolutional autoencoder, the SFC-based convolutional autoencoders have: (1)~a transformation from multi-dimensional data to 1D data and vice versa, based on SFCs; (2)~sparse layers at either end of the network; and (3)~nearest-neighbour smoothing (included in some of the networks). The transformation from multi-dimensional data to 1D data is determined by one or more space-filling curves and occurs in the first and final layers of the SFC-based autoencoders. Immediately after the first layer and before the final layer are the sparse layers. The function of the sparse layers is to apply some smoothing to the results and to decide the priority of the feature maps when recombining the data. The amount of smoothing provided by these layers can be increased by including a node's nearest neighbours on the space-filling curve within the smoothing.

On transforming the multi-dimensional data to 1D data with the SFCs, there may be some associated CPU speed advantages as 1D arrays generally map straightforwardly onto the various memory hierarchies of modern CPUs or GPUs. This could make the SFC-based CNNs computationally faster than the classical multi-dimensional CNNs. However, we have found, through considerable trial and error, that it is necessary to use roughly uniform numbers of channels in the convolutional layers of the SFC-based networks. This is in contrast with the classical CNNs, which seem to work best when the number of channels increases in the encoder and decreases in the decoder.  The more channels, the greater the computational speed and memory requirements, especially on layers with more neurons. Thus the classical CNN has a computational advantage over the SFC-based CNN in this regard.

To motivate the sparse layers further, consider a node on an SFC which is close in physical space to another node (or, more precisely, close on the graph associated with the discretisation stencil), but far away in terms of node numbering on the SFC. These nodes are disconnected, in some sense, which could potentially lead to a large variation in values occurring. This can result in discontinuities in the outputs from the autoencoder where there should be none. Having multiple SFCs allows such nodes to have low associated weights associated in the mapping to the final output layer, thus  increasing the accuracy of the values of the output neurons. Therefore, for multiple space-filling curves, this layer can be thought of as a means of letting the network decide how to optimally combine results from two space-filling curves to form the final output result neural network. The transpose of this architecture is used for the encoder (as is often done with autoencoders). This effectively takes a weighted average of a node on a SFC with its two neighbouring nodes and forms another neuron from this as well as combining (in some optimal way) SFC results to form feature maps that are fed into the rest of the autoencoder.

The notation used to describe the networks is given in Section~\ref{sec:notation}. In the subsequent sections, the autoencoders used in this paper are described in detail. For the structured data on a ${128\times 128}$ grid, we construct (1)~a classical 2D autoencoder given in Section~\ref{sec:structured_CAE}; (2)~an autoencoder based on one space-filling curve described in Section~\ref{sec:structured_CAE_SFC}; (3)~an autoencoder based on two space-filling curves described in Section~\ref{sec:structured_CAE_2SFC}; (4)~an autoencoder based on two space-filling curves with smoothing from the nearest neighbours given in Section~\ref{sec:structured_CAE_2SFC_NN}. For the unstructured data, given on a mesh of 20,550 nodes, each with two velocity components, we construct (1)~an autoencoder based on one space-filling curve with nearest-neighbour smoothing given in Section~\ref{sec:unstructured_CAE_SFC_NN}; and (2)~an autoencoder based on two space-filling curves with nearest-neighbour smoothing described in Section~\ref{sec:unstructured_CAE_2SFC_NN}. 

\begin{figure}[htbp]
\centering %[width=15cm
\includegraphics[width=16cm]{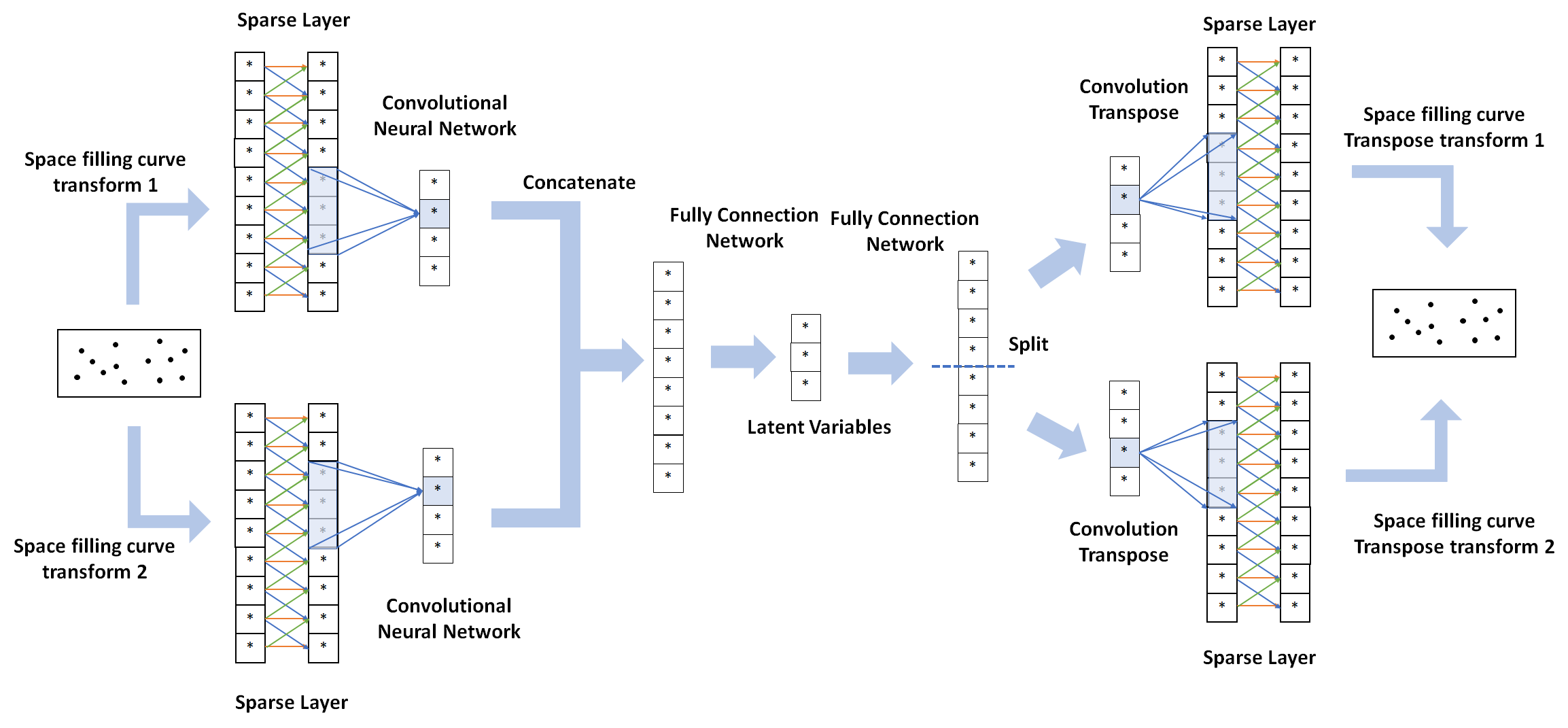}
\caption{The architecture of the SFC-based autoencoder used within this work. This diagram shows the autoencoder based on two space-filling curves with nearest-neighbour smoothing on the first and final layers. The architecture of the autoencoder based on one space-filling curve is similar, but has only one SFC branch so is half the size of that shown. Notice that there are convolutional layers associated with the encoder and decoder, and a fully-connected multi-layer percepton in the centre, which reduces the number of variables down to the required quantity.}
\label{fig:Autoencoder}
\end{figure}

\subsection{Notation used to describe autoencoder architectures}\label{sec:notation}
The notation that describes the architecture of the convolutional autoencoders used in the following equations and in Tables~\ref{ta:ArchituctureOfOne}, \ref{ta:ArchituctureOfTwo}, \ref{one-SFC-CNN-cylinder} and~\ref{two-SFC-CNN-cylinder} is explained.
%
%Layer-ordering is described as L-YYY$\cal C$X. L means L-th layer. YYY is the type of data storage: YYY=GRID means data is presented in structured $n {\times} m $ grids, YYY=SFC means data is ordered by space filling curve, while YYY=FEM means data is on a finite element mesh,  $\cal C$ indicates which SFC is being used ($\cal C$ = 1 or 2 in our cases), X is optional and is used when compressing velocity data (for flow past a cylinder here) and manipulating velocity components u and v say and thus X$\in\{$u,v$\}$. 
%
The layer ordering is described as L$\,$-[Type][SFC index][Component], where L refers to the layer number. `Type' refers to how the data is stored and can be `GRID/grid' (stored on an ${m\times n}$ grid), `SFC/sfc' (ordered by space-filling curve) or `FEM/fem' (stored on a finite-element mesh). `SFC index' indicates which SFC is being used (either~1 or~2 in our case) and `Component' is an optional description, which, when compressing velocities, indicates that component as~$u$ or~$v$. 
`Identity' in the CNN tables refers to the use of no activation function (that is, the identity function itself) or simply mapping the data from the GRID/grid or FEM/fem representation to the SFC/sfc ordering or vice versa. A similar notation is used to indicate that the neurons $\bm{x}$ are associated with either the GRID/grid, SFC/sfc or FEM/fem  within the equations describing the neuron updates within the layers. 

In the tables, `1 variable' or `3 variables' for the kernel size represents the identical operations (using kernel sizes of~1 or~3) to the standard filter approach except for the fact that the filter weights that applied to the neurons are not the same but vary for each neuron. The values are optimised as part of the training process. These filters are used to form the sparse layers, which apply smoothing of the SFC-based autoencoder  solutions and increase the accuracy of these solutions. Each node in the output is determined by one or three nodes in the input (`1~variable' or `3~variables'). Connecting one node in the output with three nodes in the input is referred to here as nearest-neighbour smoothing.  We call these sparse layers, because the number of parameters is much less than if the layers were fully connected. We have found their use to be essential in order to reduce the noise one would otherwise see in the SFC-based autoencoder results. 

To identify neighbouring nodes in the SFC ordering, for layer~$\ell$ and node number~$i$, we use ${{}^{sfc}\bm{x^{+}_{\ell}}}$ to represent the neighbour (with an increase in SFC node number) of ${{}^{sfc}\bm{x_\ell}}$ (that is ${}^{sfc}x^+_{\ell,i}={}^{sfc}x_{\ell,i+1}$), and use ${{}^{sfc}\bm{x^{-}_\ell}}$ to represent the neighbour (with a decrease in SFC node number) of ${{}^{sfc}\bm{x}}$ (that is ${}^{sfc}x^-_{\ell,i}={}^{sfc}x_{\ell,i-1}$).

\subsection{Convolutional autoencoders for data on structured meshes}
\subsubsection{Architecture of a classical 2D convolutional autoencoder}\label{sec:structured_CAE}
Table~\ref{classicCNN} shows details of the classical 2D convolutional autoencoder applied to data on a structured ${128\times 128}$ grid. Layers~1 to~4 and~11 to~14 are convolutional layers, and the fully-connected layers are at the centre of the network, in layers~5 to~10. As is often done, the number of channels increases through the enocder and decreases through the decoder.

\begin{table}[htbp]
\centering
%\resizebox{\textwidth}{!}{%
\begin{tabular}{|c|c|c|c|c|c|c|c|}
\hline
layer&input size&kernel size&channels&stride&padding&output size&activation\\
\hline
1& (1, 128, 128)& 5 $\times$ 5 & 2 & 2 $\times$ 2 & 2 & (2, 64, 64) & ReLU\\
\hline
2& (2, 64, 64) & 5 $\times$ 5 & 4 & 2 $\times$ 2 & 2 & (4, 32, 32)& ReLU\\
\hline
3& (4, 32, 32) & 5 $\times$ 5 & 8 & 2 $\times$ 2 & 2 & (8, 16, 16)& ReLU\\
\hline
4& (8, 16, 16) & 5 $\times$ 5 & 16 & 2 $\times$ 2 & 2 & (16, 8, 8)& ReLU\\
\hline
\multicolumn{8}{|c|}{Flatten the data to 1D, then use fully-connected layers to compress further}\\
\hline
5& 1024 (= 16$\times$8$\times$8) & \multicolumn{4}{c|}{ } & 256& ReLU\\
\hline
6& 256 & \multicolumn{4}{c|}{ } & 64& ReLU\\
\hline
7& 64 & \multicolumn{4}{c|}{ } & 16& ReLU\\
\hline
8& 16 & \multicolumn{4}{c|}{ } & 64& ReLU\\
\hline
9& 64 & \multicolumn{4}{c|}{ } & 256& ReLU\\
\hline
10& 256 & \multicolumn{4}{c|}{ } & 1024& ReLU\\
\hline
\multicolumn{8}{|c|}{Convert data from 1024 to (16, 8, 8) then use convolutional layers to reconstruct the data}\\
\hline
11& (16, 8, 8)& 5 $\times$ 5 & 8 & 2 $\times$ 2 & 2 & (8, 16, 16) & ReLU\\
\hline
12& (8, 16, 16) & 5 $\times$ 5 & 4 & 2 $\times$ 2 & 2 & (4, 32, 32)& ReLU\\
\hline
13& (4, 32, 32) & 5 $\times$ 5 & 2 & 2 $\times$ 2 & 2 & (2, 64, 64)& ReLU\\
\hline
14& (2, 64, 64) & 5 $\times$ 5 & 1 & 2 $\times$ 2 & 2 & (1, 128, 128)& ReLU\\
\hline
\end{tabular}
%}
\caption{The structure of the classical 2D convolutional autoencoder.}
\label{classicCNN}
\end{table}

\subsubsection{Architecture of the convolutional autoencoder based on one space-filling curve}\label{sec:structured_CAE_SFC}
% (see table~\ref{ta:ArchituctureOfOne})
As shown in Table~\ref{ta:ArchituctureOfOne}, the autoencoder based on one SFC first transforms the multi-dimensional data to 1D data with the space-filling curve ordering. After this, are the sparse layers (layers~1 and~2). Then follows a classical 1D convolutional autoencoder. At the output of this, there are more sparse layers (layers~17 and~18) followed by a transformation of the data from 1D to the original multi-dimensional form. Here we describe the sparse layers. The filters of the sparse layers are defined using weight vectors, as these filters change at each node. This results in an efficient implementation of the sparse layers. These weights effectively provide smoothing of the outputs of the SFC-based autoencoders, without which, would lead to noisy results.

\paragraph{Layers 1 and 2 (sparse layers)}\hfill\\
%\noindent Layers 1 and 2 (sparse layer - Space Filling Curve (SFC) 1):\\
The input,  ${{}^{grid}\bm{x}\in \mathbb{R}^{128\times 128}}$, is transformed 
using the Hilbert curve  mapping to a 1D vector ${}^{sfc}\bm{x_1}\in \mathbb{R}^{16384}$ ($128\times 128=16384$) in layer~1. Then the output of the sparse layer~2 is: 
\begin{equation}
%\label{eq:straight_line}
{}^{sfc}\bm{x_2}=f({}^{sfc}\bm{w}_{\bm{2}} \odot {}^{sfc}\bm{x_1}+{}^{sfc}\bm{b}_{\bm{2}})\,, 
\end{equation}
where ${{}^{sfc}\bm{w}_{\bm{2}}\in \mathbb{R}^{16384}}$ is the weight vector (subject to neural network training), ${{}^{sfc}\bm{b}_{\bm{2}}\in \mathbb{R}^{16384}}$ is the bias (again subject to training), $f$ is the ReLU activation function, and $\odot$ is the Hadamard product indicating entry-wise multiplication. Note that ${{}^{sfc}\bm{x_2}\in \mathbb{R}^{16384}}$ is the input for layer~3 of the neural network. 
See Table~\ref{ta:ArchituctureOfOne} for a description of the  convolutional layers, which take input from layer~2. 

\begin{table}[htbp]
\centering
\resizebox{\textwidth}{!}{%
\begin{tabular}{|c|c|c|c|c|c|c|c|}
\hline
layer&input size&kernel size&channels&stride&padding&output size&activation\\
\hline
1-GRID& (1, 16384, GRID)& 1 Identity & 1 & 1 & 0 & (1, 16384, SFC1) & Identity\\
\hline
2-SFC1& (1, 16384, SFC1)& 1 Variable (1$\times$ 16384)& 1 & 1 & 0 & (1, 16384, SFC1) & ReLU\\
\hline
3-SFC1& (1, 16384, SFC1)& 32 & 2 & 4 & 16 & (2, 4097, SFC1) & ReLU\\
\hline
4-SFC1& (2, 4097, SFC1)& 32 & 4 & 4 & 16 & (4, 1025, SFC1) & ReLU\\
\hline
5-SFC1& (4, 1025, SFC1)& 32 & 8 & 4 & 16 & (8, 257, SFC1) & ReLU\\
\hline
6-SFC1& (8, 257, SFC1)& 32 & 16 & 4 & 16 & (16, 65, SFC1) & ReLU\\
\hline
\multicolumn{8}{|c|}{Flatten the data to 1D, then use fully-connected layers to compress further}\\
\hline
7-SFC1& 1040 (= 16$\times$65) & \multicolumn{4}{c|}{ } & 256& ReLU\\
\hline
8-SFC1& 256 & \multicolumn{4}{c|}{ } & 64& ReLU\\
\hline
9-SFC1& 64 & \multicolumn{4}{c|}{ } & 16& ReLU\\
\hline
10-SFC1 & 16 & \multicolumn{4}{c|}{ } & 64& ReLU\\
\hline
11-SFC1& 64 & \multicolumn{4}{c|}{ } & 256& ReLU\\
\hline
12-SFC1& 256 & \multicolumn{4}{c|}{ } & 1024& ReLU\\
\hline
\multicolumn{8}{|c|}{Convert data from 1024 to (16, 64, SFC1) then use convolutional layers to reconstruct the data}\\
\hline
13-SFC1& (16, 64, SFC1)& 32 & 8 & 4 & 14 & (8, 256, SFC1) & ReLU\\
\hline
14-SFC1& (8, 256, SFC1)& 32 & 4 & 4 & 14 & (4, 1024, SFC1) & ReLU\\
\hline
15-SFC1& (4, 1024, SFC1)& 32 & 2 & 4 & 14 & (2, 4096, SFC1) & ReLU\\
\hline
16-SFC1& (2, 4096, SFC1)& 32 & 1 & 4 & 14 & (1, 16384, SFC1) & ReLU\\
\hline
17-SFC1& (1, 16384, SFC1)& 1 Variable (1$\times$ 16384)& 1 & 1 & 0 & (1, 16384, SFC1) & Identity\\
\hline
18-SFC1& (1, 16384, SFC1)& 1 Identity  & 1 & 1 & 0 & (1, 16384, GRID) & ReLU\\
\hline
\end{tabular}
}
\caption{The structure of the SFC-based CNN with one Hilbert curve for a $128\times 128$ structured grid. The filter size for the smoothing layers at the start and end of the auto-encoder is 1 but replace `1 Variable' with `3 Variables' for the nearest-neighbour smoothing approach. The number of neurons in the bottleneck layer can be adjusted from 16 to any suitable value less than or equal to 64. Whenever a non-identity activation term is used a bias is added to each neuron.}
\label{ta:ArchituctureOfOne}
\end{table}

\paragraph{Layers 17 and 18 (sparse layers)}\hfill\\
Given the input ${{}^{sfc}\bm{x_{16}}\in \mathbb{R}^{16384}}$, the output of the sparse layer 17 (which is 
the final output of the SFC-based CNN) can be written as:
\begin{equation}
%\label{eq:straight_line}
{}^{sfc}\bm{x_{17}}=f({}^{sfc}\bm{w_{17}}\odot {}^{sfc}\bm{x_{16}}+{}^{sfc}\bm{b_{17}})\, , 
\end{equation}
where ${{}^{sfc}\bm{w_{17}}\in \mathbb{R}^{16384}}$ is the weight vector, ${{}^{sfc}\bm{b_{17}}\in \mathbb{R}^{16384}}$ is the bias, and $f$ is the ReLU activation function. 
Using the inverse space-filling  curve mapping, the 1D data,  ${}^{sfc}\bm{x_{17}}$, is transformed back to 2D data on the structured grid, 
${}^{grid}{\bm{x_{18}}\in \mathbb{R}^{128 \times 128}}$.

\subsubsection{Architecture of a convolutional autoencoder based on two space-filling curves}\label{sec:structured_CAE_2SFC}
% (see table~\ref{ta:ArchituctureOfTwo}) 
The architecture of autoencoder based on two space-filling curves is given in Table~\ref{ta:ArchituctureOfTwo}. When using two space-filling curves, we keep the neuron values associated with these curves separate in the convolutional layers, but bring this information together within the fully-connected laters at the centre, as well as the input and output layers. The sparse layers are now described in detail.

\paragraph{Layers~1 and~2 (sparse layers)}\hfill\\
The data on the grid, ${}^{grid}{\bm{x}\in \mathbb{R}^{128\times 128}}$,  is transformed to two 1D vectors using the mappings from the two space-filling curves to produce ${{}^{sfc1}\bm{x}_{\bm{1}}, \, {}^{sfc2}\bm{x}_{\bm{1}} \in \mathbb{R}^{16384}}$ in layer~1. The output of the sparse layer~2 is:
\begin{equation}
{}^{sfc1}\bm{x_{2}}=f({}^{sfc1}\bm{w_{2}}\odot{}^{sfc1}\bm{x_{1}}+{}^{sfc1}\bm{b_{2}})\, ,
\end{equation}
\begin{equation}
{}^{sfc2}\bm{x_{2}}=f({}^{sfc2}\bm{w_{2}}\odot{}^{sfc2}\bm{x_{1}}+{}^{sfc2}\bm{b_{2}})\, ,
\end{equation}
where ${}^{sfc1}\bm{w}_{2}, {}^{sfc2}\bm{w}_{2} \in \mathbb{R}^{16384}$ are the weight vectors, ${{}^{sfc1}\bm{b_{2}}}, {{}^{sfc2}\bm{b_{2}}}\in \mathbb{R}^{16384}$ are the bias vectors, and $f$~is the ReLU activation function. The input for layer~3 takes the form ${{}^{sfc1}\bm{x_{2}}},{{}^{sfc2}\bm{x_{2}}} \in \mathbb{R}^{16384}$. See Table~\ref{ta:ArchituctureOfTwo} for a description of the following convolutional layers of the network. 

\paragraph{Layers~18 and~19 (sparse layers)}\hfill\\
Given the inputs ${{}^{sfc1}\bm{x_{17}}\in \mathbb{R}^{16384}}$ and ${{}^{sfc2}\bm{x_{17}}\in \mathbb{R}^{16384}}$, the output of the sparse layer~18 is:
\begin{equation}
%\label{eq:straight_line}
{}^{sfc1}\bm{x_{18}}={}^{sfc1}\bm{w_{18}}\odot{}^{sfc1}\bm{x_{17}}\, ,
\end{equation}
\begin{equation}
%\label{eq:straight_line}
{}^{sfc2}\bm{x_{18}}={}^{sfc2}\bm{w_{18}}\odot{}^{sfc2}\bm{x_{17}}\, , 
\end{equation}
where ${}^{sfc1}\bm{w_{18}}, {}^{sfc2}\bm{w_{18}} \; \in \mathbb{R}^{16384}$ 
are the weight vectors. 
Using the inverse mappings from the first and second space-filling curves, the vectors ${}^{sfc1}\bm{x_{18}}$ and ${}^{sfc2}\bm{x_{18}}$ are transformed to obtain ${}^{grid1}\bm{x_{18}}\in\mathbb{R}^{128 \times 128}$ and ${}^{grid2}\bm{x_{18}}\in\mathbb{R}^{128 \times 128}$. Added to this is the bias ${}^{grid}\bm{b_{19}}\in \mathbb{R}^{128 \times 128}$. The output of the network is 
 ${{}^{grid}\bm{x_{19}} \in \mathbb{R}^{128 \times 128}}$  and is obtained from: 
\begin{equation}
%\label{eq:straight_line}
{}^{grid}\bm{x_{19}}=f({}^{grid1}\bm{x_{18}} + {}^{grid2}\bm{x_{18}} + {}^{grid}\bm{b_{19}}), 
\end{equation}
in which $f$ is the ReLU activation function.

\begin{table}[htbp]
\centering
\resizebox{\textwidth}{!}{%
\begin{tabular}{|c|c|c|c|c|c|c|c|}
\hline
layer ordering&input size \&  ordering&kernel size&channels&stride&padding&output size \& ordering&activation\\
\hline
1-GRID& (1, 16384, GRID)& 1 Identity & 1 & 1 & 0 & (1, 16384, SFC1) & Identity\\
\hline
2-SFC1& (1, 16384, SFC1)& 1 Variable ($1\times 16384$)& 2 & 1 & 0 & (1, 16384, SFC1) & ReLU\\
\hline
3-SFC1& (1, 16384, SFC1)& 32 & 2 & 4 & 16 & (2, 4097, SFC1) & ReLU\\
\hline
4-SFC1& (2, 4097, SFC1)& 32 & 4 & 4 & 16 & (4, 1025, SFC1) & ReLU\\
\hline
5-SFC1& (4, 1025, SFC1)& 32 & 8 & 4 & 16 & (8, 257, SFC1) & ReLU\\
\hline
6-SFC1& (8, 257, SFC1)& 32 & 16 & 4 & 16 & (16, 65, SFC1) & ReLU\\
\hline
1-GRID& (1, 16384, GRID)& 1 Identity & 1 & 1 & 0 & (1, 16384, SFC2) & Identity\\
\hline
2-SFC2& (1, 16384)& 1 Variable ($1\times 16384$)& 1 & 1 & 0 & (1, 16384, SFC2) & ReLU\\
\hline
3-SFC2& (1, 16384)& 32 & 2 & 4 & 16 & (2, 4097, SFC2) & ReLU\\
\hline
4-SFC2& (2, 4097)& 32 & 4 & 4 & 16 & (4, 1025, SFC2) & ReLU\\
\hline
5-SFC2& (4, 1025)& 32 & 8 & 4 & 16 & (8, 257, SFC2) & ReLU\\
\hline
6-SFC2& (8, 257)& 32 & 16 & 4 & 16 & (16, 65, SFC2) & ReLU\\
\hline
\multicolumn{8}{|c|}{Flatten the output data of layer 6-SFC1 and 6-SFC2 to 1D - concatenate 2 sequences}\\
\hline
7& 2080 (= 16$\times$65 + 16$\times$65) & \multicolumn{4}{c|}{ } & 512& ReLU\\
\hline
8& 512 & \multicolumn{4}{c|}{ } & 128& ReLU\\
\hline
9& 128 & \multicolumn{4}{c|}{ } & 16& ReLU\\
\hline
10& 16 & \multicolumn{4}{c|}{ } & 128& ReLU\\
\hline
11& 128 & \multicolumn{4}{c|}{ } & 512& ReLU\\
\hline
12& 512 & \multicolumn{4}{c|}{ } & 2048& ReLU\\
\hline
\multicolumn{8}{|c|}{Split the data into 2 sequence as the input of layer 11-SFC1 and 11-SFC2, convert from 1024 to (16, 64)}\\
\hline
13-SFC1& (16, 64, SFC1)& 32 & 8 & 4 & 14 & (8, 256, SFC1) & ReLU\\
\hline
14-SFC1& (8, 256, SFC1)& 32 & 4 & 4 & 14 & (4, 1024, SFC1) & ReLU\\
\hline
15-SFC1& (4, 1024, SFC1)& 32 & 2 & 4 & 14 & (2, 4096, SFC1) & ReLU\\
\hline
16-SFC1& (2, 4096, SFC1)& 32 & 1 & 4 & 14 & (1, 16384, SFC1) & ReLU\\
\hline
17-SFC1& (2, 16384, SFC1)& 1 Variable ($1\times 16384$)& 1 & 1 & 0 & (1, 16384, SFC1) & Identity\\
\hline
18-SFC1& (1, 16384, SFC1)& 1 Identity & 1 & 1 & 0 & (1, 16384, GRID1) & Identity\\
\hline
13-SFC2& (16, 64, SFC2)& 32 & 8 & 4 & 14 & (8, 256, SFC2) & ReLU\\
\hline
14-SFC2& (8, 256, SFC2)& 32 & 4 & 4 & 14 & (4, 1024, SFC2) & ReLU\\
\hline
15-SFC2& (4, 1024, SFC2)& 32 & 2 & 4 & 14 & (2, 4096, SFC2) & ReLU\\
\hline
16-SFC2& (2, 4096, SFC2)& 32 & 1 & 4 & 14 & (1, 16384, SFC2) & ReLU\\
\hline
17-SFC2& (2, 16384, SFC2)& 1 Variable ($1\times 16384$)& 1 & 1 & 0 & (1, 16384, SFC2) & Identity\\
\hline
18-SFC2& (1, 16384, SFC2)& 1 Identity & 1 & 1 & 0 & (1, 16384, GRID2) & Identity\\
\hline
19-GRID& ($1$, 16384)& 1 Identity GRID1 + GRID2 & 1 & 1 & 0 & (1, 16384, GRID) & ReLU\\
\hline
\end{tabular}
%}
}
\caption{The structure of the SFC-based CNN with two space-filling  curves for a ${128\times 128}$ structured grid. The filter size for the smoothing layers at the start and end of the autoencoder is~1, however replacing `1 Variable' with `3 Variable' would give the nearest-neighbour smoothing.  The number of neurons in the bottleneck layer can be adjusted from~16 to any suitable value.  Whenever a non-identity activation term is used a bias is added to each neuron.}
\label{ta:ArchituctureOfTwo}
\end{table}

\subsubsection[Architecture of a CAE based on two space-filling curves with neighbours]{Architecture of a convolutional autoencoder based on two space-filling curves with nearest-neighbour smoothing}\label{sec:structured_CAE_2SFC_NN}
% (see table~\ref{ta:ArchituctureOfTwo})

The only difference between this SFC-based CNN and the previous one is the introduction of the nearest neighbours in the sparse layers. This effectively reduces the noise in the output of the SFC-based CNN. 

\paragraph{Layers 1 and 2 (sparse layers)}\hfill\\
The input to the network, ${{}^{grid}\bm{x}\in \mathbb{R}^{128\times 128}}$, is transformed 
to ${{}^{sfc1}\bm{x}_{\bm{1}}, {}^{sfc2}\bm{x}_{\bm{1}}\in \mathbb{R}^{16384}}$ by the two space-filling curve mappings in layer~1. In SFC ordering, ${{}^{sfc{\cal C}}\bm{x^{+}_{\bm{1}}}\in \mathbb{R}^{16384}}$ is the neighbour (with an increase in SFC node number) of ${{}^{sfc{\cal C}}\bm{x}_{\bm{1}}}$, and ${{}^{sfc{\cal C}}\bm{x^{-}_{\bm{1}}}\in \mathbb{R}^{16384}}$ is the neighbour (with an decrease in SFC node number) of ${{}^{sfc{\cal C}}\bm{x}_{\bm{1}}}$,  in which ${\cal C}\in\{1,2\}$ is the number of the SFC. The output of the sparse layer~2 is given as:
\begin{equation}
{}^{sfc1}\bm{x_{2}}=f({}^{sfc1}\bm{w_{2}}\odot{}^{sfc1}\bm{x_{1}}+{}^{sfc1}\bm{w_{2}^{+}}\odot{}^{sfc1}\bm{x^{+}_{1}}+{}^{sfc1}\bm{w_{2}^{-}}\odot{}^{sfc1}\bm{x^{-}_{1}}+{}^{sfc1}\bm{b_{2}})\, ,
\end{equation}
\begin{equation}
{}^{sfc2}\bm{x_{2}}=f({}^{sfc2}\bm{w_{2}}\odot{}^{sfc2}\bm{x_{1}}+{}^{sfc2}\bm{w_{2}^{+}}\odot{}^{sfc2}\bm{x^{+}_{1}}+{}^{sfc2}\bm{w_{2}^{-}}\odot{}^{sfc2}\bm{x^{-}_{1}}+{}^{sfc2}\bm{b_{2}})\, ,
\end{equation}
where ${}^{sfc1}\bm{w_{2}}, \, {}^{sfc2}\bm{w_{2}}, \, {}^{sfc1}\bm{w_{2}^{+}}, \, {}^{sfc2}\bm{w_{2}^{+}}, \, {}^{sfc1}\bm{w_{2}^{-}}, \, {}^{sfc2}\bm{w_{2}^{-}}\in \mathbb{R}^{16384}$ are the weight vectors, ${}^{sfc1}\bm{b_{2}}$, ${}^{sfc2}\bm{b_{2}}$ $\in \mathbb{R}^{16384}$ are the bias vectors, and $f$ is the ReLU activation function. The input to layer~3 is ${}^{sfc1}\bm{x_{2}}$ and ${}^{sfc2}\bm{x_{2}} $. See Table~\ref{ta:ArchituctureOfTwo} for a description of the convolutional layers.

\paragraph{Layers 18 and 19 (sparse layers)}\hfill\\
Given the inputs ${{}^{sfc{\cal C}}\bm{x_{17}}\in \mathbb{R}^{16384}}$ $\forall {\cal C}\in\{1,2\}$, in the SFC ordering ${{}^{sfc{\cal C}}\bm{x^{+}_{\bm{17}}}\in \mathbb{R}^{16384}}$ is the neighbour (with an increase in SFC node number) of ${{}^{sfc{\cal C}}\bm{x}_{\bm{17}}}$, and ${{}^{sfc{\cal C}}\bm{x^{-}_{\bm{17}}}\in \mathbb{R}^{16384}}$ is the neighbour (with an decrease in SFC node number) of ${{}^{sfc{\cal C}}\bm{x}_{\bm{17}}}$. The output of layer~18 is:
\begin{equation}
%\label{eq:straight_line}
{}^{sfc1}\bm{x_{18}}={}^{sfc1}\bm{w_{18}}\odot{}^{sfc1}\bm{x_{17}}+{}^{sfc1}\bm{w_{18}^{+}}\odot{}^{sfc1}\bm{x^{+}_{17}}+{}^{sfc1}\bm{w_{18}^{-}}\odot{}^{sfc1}\bm{x^{-}_{17}}\, ,
\end{equation}
\begin{equation}
%\label{eq:straight_line}
{}^{sfc2}\bm{x_{18}}={}^{sfc2}\bm{w_{18}}\odot{}^{sfc2}\bm{x_{17}}+{}^{sfc2}\bm{w_{18}^{+}}\odot{}^{sfc2}\bm{x^{+}_{17}}+{}^{sfc2}\bm{w_{18}^{-}}\odot{}^{sfc2}\bm{x^{-}_{17}}\, ,
\end{equation}
where ${}^{sfc1}\bm{w_{18}}, \, {}^{sfc2}\bm{w_{18}}, \, {}^{sfc1}\bm{w_{18}^{+}}, \, {}^{sfc2}\bm{w_{18}^{+}}, \, {}^{sfc1}\bm{w_{18}^{-}}, \, {}^{sfc2}\bm{w_{18}^{-}}\in \mathbb{R}^{16384}$ are the weight vectors. 
Using the inverse of both space-filling curve mappings, transform ${}^{sfc1}\bm{x_{18}}$ and ${}^{sfc2}\bm{x_{18}}$ 
to ${}^{grid1}\bm{x_{18}}\in \mathbb{R}^{128 \times 128}$ and ${}^{grid2}\bm{x_{18}}\in \mathbb{R}^{128 \times 128}$ respectively, then add the bias ${{}^{grid}\bm{b_{19}}\in \mathbb{R}^{128 \times 128}}$ to obtain: 
\begin{equation}
{}^{grid}\bm{x_{19}}=f({}^{grid1}\bm{x_{18}} + {}^{grid2}\bm{x_{18}} + {}^{grid}\bm{b_{19}}), 
\end{equation}
in which $f$ is the ReLU activation function and ${{}^{grid}\bm{x_{19}}\in \mathbb{R}^{128 \times 128}}$ is the output of this SFC-based CNN on the original grid.

\subsection{Architectures of convolutional autoencoders for data on unstructured meshes}
\subsubsection[Architecture of a CAE based on one space-filling curve with neighbours]{An autoencoder based on one space-filling curve with nearest-neighbour smoothing}\label{sec:unstructured_CAE_SFC_NN}

Since flow past a cylinder is a more demanding compression problem than the previous two idealised cases, two channels are used for the sparse layer filters. If only one channel were used to be used, as in the idealised cases,  the SFC-based CNN outputs tended to be noisy. In addition, within the sparse layers, the velocity components~$u$ and~$v$ are kept separate, in order to reduce the number of weights required in these layers and also to apply the smoothing to the velocity components separately. 
These modifications to the CNN are also applied to the SFC-based CNN with two SFCs,  also used for flow past a cylinder.  See Table~\ref{ta:ArchituctureOfTwo} for a description of this autoencoder. We now describe the sparse layers in detail.

\paragraph{Layers 1 and 2 (sparse layers)}\hfill\\
The input to the neural network is ${{}^{fem}\bm{x}\in \mathbb{R}^{20550\times 2}}$ and contains two channels; one for the first velocity component,~$u$, and the second for the other velocity component,~$v$. These are held in the FEM DG ordering of the nodal values of the velocities. The input, ${}^{fem}\bm{x}$, is mapped to the SFC ordering, ${}^{sfc1}\bm{x_1}$,  using the SFC mapping in layer~1. By splitting the channels of ${}^{sfc1}\bm{x_1}$, we obtain ${}^{sfc1u}\bm{x_1}$ and ${}^{sfc1v}\bm{x_1}$, associated with the two velocity components respectively. 
In the SFC ordering, ${{}^{sfc1X}\bm{x^{+}_1}\in \mathbb{R}^{20550\times 1}}$ $\forall X\in\{u,v\}$ are the neighbours (with an increase in SFC node number) of ${{}^{sfc1X}\bm{x_1}}$ (that is ${}^{sfc1X}x^+_{i,k}={}^{sfc1X}x_{i+1,k}$), ${{}^{sfc1X}\bm{x^{-}_1}\in \mathbb{R}^{20550\times 1}}$ is the neighbour (with an decrease in SFC node number) of ${{}^{sfc1X}\bm{x_1}}$ (that is ${}^{sfc1X}x^-_{i,k}={}^{sfc1X}x_{i-1,k}$). The output of the sparse layer 2 is then:
\begin{eqnarray}
{}^{sfc1u}\bm{x_{2}}&=&f\left({}^{sfc1u}\bm{w_{2}}\odot \conca2({}^{sfc1u}\bm{x_1})+{}^{sfc1u}\bm{w_{2}^{+}}\odot 
\conca2({}^{sfc1u}\bm{x^{+}_1})\right.\nonumber\\
&& \hspace{6.5cm}\left.+{}^{sfc1u}\bm{w_{2}^{-}}\odot \conca2({}^{sfc1u}\bm{x^{-}_1})+{}^{sfc1u}\bm{b_{2}}\right)\, , \\[3mm]
{}^{sfc1v}\bm{x_{2}}&=&f\left({}^{sfc1v}\bm{w_{2}}\odot \conca2({}^{sfc1v}\bm{x_1})+{}^{sfc1v}\bm{w_{2}^{+}}\odot \conca2({}^{sfc1v}\bm{x^{+}_1})\right. \nonumber\\
&& \hspace{6.5cm} \left.+{}^{sfc1v}\bm{w_{2}^{-}}\odot \conca2({}^{sfc1v}\bm{x^{-}_1})+{}^{sfc1v}\bm{b_{2}}\right)\, ,
\end{eqnarray}
where ${}^{sfc1u}\bm{w_{2}}, \, {}^{sfc1v}\bm{w_{2}}, \, , {}^{sfc1u}\bm{w_{2}^{+}}, \, {}^{sfc1v}\bm{w_{2}^{+}}, \, {}^{sfc1u}\bm{w_{2}^{-}}, \, {}^{sfc1v}\bm{w_{2}^{-}}\in \mathbb{R}^{20550\times 2}$ are the weight vectors, ${}^{sfc1u}\bm{b_{2}}$, ${}^{sfc1v}\bm{b_{2}} \in \mathbb{R}^{20550\times 2}$ are the bias vectors, and $f$ is the tanh activation function. Then: 
\begin{equation}
{}^{sfc1}\bm{x_{2}}=\concat( {}^{sfc1u}\bm{x_{2}}, {}^{sfc1v}\bm{x_{2}} ) ,
\end{equation}
where ${{}^{sfc1}\bm{x_{2}}\in \mathbb{R}^{20550\times 4}}$ is the input for layer~3. The notation $\concat(\bm{a}, \bm{b})=(\bm{a}^T , \bm{b}^T)^T$ is used to represent the concatenation of two vectors into one, and the notation $\conca2(\bm{a})=(\bm{a}^T , \bm{a}^T)^T$ represents the concatenation of one vector with itself. %See table~\ref{ta:ArchituctureOfTwo} for a description of the SFC-based CNN including its convolutional layers that ${}^{sfc1}\bm{x_{2}}$ inputs into. 

\paragraph{Layers 13, 14 and 15 (sparse layers)}\hfill\\
The vector ${{}^{sfc1}\bm{x_{12}}\in \mathbb{R}^{20550 \times 4}}$ (in SFC ordering), 
is disconcatenated or split into two; ${}^{sfc1u}\bm{x_{12}}, {}^{sfc1v}\bm{x_{12}} \in \mathbb{R}^{20550 \times 2}$. 
Similar to previously ${{}^{sfc1X}\bm{x_{12}^{+}}\in \mathbb{R}^{20550 \times 2}}$ $\forall X\in\{u,v\}$ are the neighbours (with an increase in SFC node number) of ${{}^{sfc1X}\bm{x_{12}}}$, and ${{}^{sfc1X}\bm{x_{12}^{-}}\in \mathbb{R}^{20550 \times 2}}$ are the neighbours (with a decrease in SFC node number) of ${{}^{sfc1X}\bm{x_{12}}}$. Thus: 
\begin{equation}%\label{eq:straight_line}
{}^{sfc1u}\bm{x_{13}}={\sumchannels}(
{}^{sfc1u}\bm{w_{13}}\odot{}^{sfc1u}\bm{x_{12}} + {}^{sfc1u}\bm{w_{13}^{+}}\odot{}^{sfc1u}\bm{x_{12}^{+}} + {}^{sfc1u}\bm{w_{13}^{-}}\odot{}^{sfc1u}\bm{x_{12}^{-}} 
),
\end{equation}
\begin{equation}%\label{eq:straight_line}
{}^{sfc1v}\bm{x_{13}}=\sumchannels(
{}^{sfc1v}\bm{w_{13}}\odot{}^{sfc1v}\bm{x_{12}} + {}^{sfc1v}\bm{w_{13}^{+}}\odot{}^{sfc1v}\bm{x_{12}^{+}} + {}^{sfc1v}\bm{w_{13}^{-}}\odot{}^{sfc1v}\bm{x_{12}^{-}} 
),
\end{equation}
where ${{}^{sfc1u}\bm{w_{13}}}, \, {{}^{sfc1v}\bm{w_{13}}}, \, {}^{sfc1u}\bm{w_{13}^{+}}, \, {}^{sfc1v}\bm{w_{13}^{+}}, \, {}^{sfc1u}\bm{w_{13}^{-}}, \, {}^{sfc1v}\bm{w_{13}^{-}}$ $\in\mathbb{R}^{20550 \times 2}$ are the weight vectors and $ {}^{sfc1u}\bm{x_{13}}$, ${}^{sfc1v}\bm{x_{13}}\in\mathbb{R}^{20550 \times 1}$. 
The operator $\sumchannels$ simply adds together the vectors in each channel to produce a one-channel output. Using the inverse SFC mappings, 
${}^{sfc1u}\bm{x_{13}}$ and ${}^{sfc1v}\bm{x_{13}}$ are then transformed to ${}^{fem1u}\bm{x_{14}}$ and ${}^{fem1v}\bm{x_{14}}$ in the original FEM mesh ordering. The output neurons in FEM ordering are: 
\begin{equation}
{}^{fem}\bm{x_{15}}=f\left(\concat({}^{fem1u}\bm{x_{14}}, {}^{fem1v}\bm{x_{14}}) + {}^{fem}\bm{b_{15}}\right), 
\end{equation}
in which ${}^{fem}\bm{b_{15}}\in \mathbb{R}^{20550 \times 2}$ 
is the bias, $f$ is the tanh activation function and the two channels of ${}^{fem}\bm{x_{15}} \in \mathbb{R}^{20550 \times 2}$ contain the~$u$ and~$v$ nodal velocity components, respectively, on the FEM DG mesh.

\begin{table}[htbp]
\centering
\resizebox{\textwidth}{!}{%
\begin{tabular}{|c|c|c|c|c|c|c|c|}
\hline
layer&input size&kernel size&channels&stride&padding&output size&activation\\
\hline
1-FEM& (2, 20550, FEM)& 1 Identity & 2 & 1 & 0 & (2, 20550, SFC1) & Identity\\
\hline
\multicolumn{8}{|c|}{Split (2, 20550, SFC1) into two:(1, 20550, SFC1$u$), (1, 20550, SFC1$v$) }\\
\hline
2-SFC1$u$& (1, 20550, SFC1$u$)& 3 Variable ($3\times20550$)& 2 & 1 & 0 & (2, 20550, SFC1$u$) & tanh\\
\hline
2-SFC1$v$& (1, 20550, SFC1$v$)& 3 Variable ($3\times20550$)& 2 & 1 & 0 & (2, 20550, SFC1$v$) & tanh\\
\hline
\multicolumn{8}{|c|}{concatenate:((2, 20550, SFC1$u$), (2, 20550, SFC1$v$)) forming (4, 20550, SFC1) }\\
\hline
3-SFC1& (4, 20550, SFC1)& 32 & 16 & 4 & 16 & (16, 5138, SFC1) & tanh\\
\hline
4-SFC1& (16, 5138, SFC1)& 32 & 16 & 4 & 16 & (16, 1285, SFC1) & tanh\\
\hline
5-SFC1& (16, 1285, SFC1)& 32 & 16 & 4 & 16 & (16, 322, SFC1) & tanh\\
\hline
6-SFC1& (16, 322, SFC1)& 32 & 16 & 4 & 16 & (16, 81, SFC1) & tanh\\
\hline
\multicolumn{8}{|c|}{Flatten the output data of layer 6 to 1D}\\
\hline
7& 1296 (= 16$\times$81) & \multicolumn{4}{c|}{ } & 128& tanh\\
\hline
8 & 128 & \multicolumn{4}{c|}{ } & 1296& tanh\\
\hline
\multicolumn{8}{|c|}{ Convert from 1296 to (16, 81, SFC1) as the input of layer 9 }\\
\hline
9-SFC1& (16, 81, SFC1)& 32 & 16 & 4 & 15 & (16, 322, SFC1) & tanh\\
\hline
10-SFC1& (16, 322, SFC1)& 32 & 16 & 4 & 15 & (16, 1286, SFC1) & tanh\\
\hline
11-SFC1& (16, 1286, SFC1)& 32 & 16 & 4 & 16 & (16, 5140, SFC1) & tanh\\
\hline
12-SFC1& (16, 5140, SFC1)& 32 & 4 & 4 & 19 & (4, 20550, SFC1) & tanh\\
\hline
\multicolumn{8}{|c|}{Split (4, 20550, SFC1) into two:(2, 20550, SFC1$u$), (2, 20550, SFC1$v$) }\\
\hline
13-SFC1$u$& (2, 20550, SFC1$u$)& 3 Variable ($3\times20550$)& 2 & 1 & 0 & (1, 20550, SFC1$u$) & Identity\\
\hline
13-SFC1$v$& (2, 20550, SFC1$v$)& 3 Variable ($3\times20550$)& 2 & 1 & 0 & (1 20550, SFC1v) & Identity\\
\hline
14-SFC1$u$& (1, 20550, SFC1$u$)& 1 Identity & 1 & 1 & 0 & (1, 20550, FEM1u) & Identity\\
\hline
14-SFC1$v$& (1, 20550, SFC1$v$)& 1 Identity & 1 & 1 & 0 & (1, 20550, FEM1v) & Identity\\
\hline
15-FEM& ($1\times 2$, 20550, FEM)& 1 Identity $\concat$(FEM1$u$,FEM1$v$) & 2 & 1 & 0 & (2, 20550, FEM) & tanh\\
\hline
\end{tabular}
}
\caption{The architecture of the SFC-based CNN with one space-filling curve with nearest-neighbour smoothing for flow past a cylinder. The number of neurons in the bottle neck layer can be adjusted from 128 to any suitable value. Whenever a non-identity activation term is used a bias is added to each neuron.  }
\label{one-SFC-CNN-cylinder} 
\end{table}

\subsubsection[Architecture of a CAE based on two space-filling curves with neighbours]{An autoencoder based on two space-filling curves with nearest-neighbour smoothing}\label{sec:unstructured_CAE_2SFC_NN}

\paragraph{Layers 1 and 2 (sparse layers)}\hfill\\
Given the input ${}^{fem}\bm{x}\in \mathbb{R}^{20550\times 2}$, convert this to the SFC1 and SFC2 ordering ${}^{sfc1}\bm{x_1}, {}^{sfc2}\bm{x_1}$ using the SFC mappings in layer 1. By splitting the channels of ${}^{sfc1}\bm{x_1}$  and $ {}^{sfc2}\bm{x_1}$, we obtain ${}^{sfc1u}\bm{x_1}, \, {}^{sfc2u}\bm{x_1}, \, {}^{sfc1v}\bm{x_1}, \, {}^{sfc2v}\bm{x_1}$ for velocity components~$u$ and~$v$ respectively. The output of the sparse layer~2 is:
\begin{eqnarray}
{}^{sfc1u}\bm{x_{2}}&=&f\left({}^{sfc1u}\bm{w_{2}}\odot \conca2({}^{sfc1u}\bm{x_1})+{}^{sfc1u}\bm{w_{2}^{+}}\odot \conca2({}^{sfc1u}\bm{x^{+}_1})\right.\nonumber\\ &&\hspace{6.5cm}\left.+{}^{sfc1u}\bm{w_{2}^{-}}\odot \conca2({}^{sfc1u}\bm{x^{-}_1})+{}^{sfc1u}\bm{b_{2}}\right)\, ,\\[3mm]
{}^{sfc1v}\bm{x_{2}}&=&f\left({}^{sfc1v}\bm{w_{2}}\odot\conca2({}^{sfc1v}\bm{x_1})+{}^{sfc1v}\bm{w_{2}^{+}}\odot \conca2({}^{sfc1v}\bm{x^{+}_1})\right.\nonumber\\
&&\hspace{6.5cm}\left.+{}^{sfc1v}\bm{w_{2}^{-}}\odot\conca2({}^{sfc1v}\bm{x^{-}_1})+{}^{sfc1v}\bm{b_{2}}\right)\, ,\\[3mm]
{}^{sfc2u}\bm{x_{2}}&=&f\left({}^{sfc2u}\bm{w_{2}}\odot\conca2({}^{sfc2u}\bm{x_1})+{}^{sfc2u}\bm{w_{2}^{+}}\odot\conca2({}^{sfc2u}\bm{x^{+}_1})\right.\nonumber\\ 
&&\hspace{6.5cm}\left.+{}^{sfc2u}\bm{w_{2}^{-}}\odot\conca2({}^{sfc2u}\bm{x^{-}_1})+{}^{sfc2u}\bm{b_{2}}\right)\, ,\\[3mm]
{}^{sfc2v}\bm{x_{2}}&=&f\left({}^{sfc2v}\bm{w_{2}}\odot\conca2({}^{sfc2v}\bm{x_1})+{}^{sfc2v}\bm{w_{2}^{+}}\odot\conca2({}^{sfc2v}\bm{x^{+}_1})\right.\nonumber\\
&& \hspace{6.5cm}\left.+{}^{sfc2v}\bm{w_{2}^{-}}\odot\conca2({}^{sfc2v}\bm{x^{-}_1})+{}^{sfc2v}\bm{b_{2}}\right)\, ,
\end{eqnarray}
where ${}^{sfc1u}\bm{w_{2}}, \, {}^{sfc1v}\bm{w_{2}}, \, , {}^{sfc1u}\bm{w_{2}^{+}}, \, {}^{sfc1v}\bm{w_{2}^{+}}, \, {}^{sfc1u}\bm{w_{2}^{-}}, \, {}^{sfc1v}\bm{w_{2}^{-}}, \, {}^{sfc2u}\bm{w_{2}}, \, {}^{sfc2v}\bm{w_{2}}, \, , {}^{sfc2u}\bm{w_{2}^{+}}, \, {}^{sfc2v}\bm{w_{2}^{+}}, \\ {}^{sfc2u}\bm{w_{2}^{-}}, \, {}^{sfc2v}\bm{w_{2}^{-}}\in \mathbb{R}^{20550\times 2}$ are the weight vectors, ${{}^{sfc1u}\bm{b_{2}}, {}^{sfc1v}\bm{b_{2}}, \, {}^{sfc2u}\bm{b_{2}}, {}^{sfc2v}\bm{b_{2}} \in \mathbb{R}^{20550\times 2}}$ are the bias vectors, ${{}^{sfc1u}\bm{x_{2}}, {}^{sfc1v}\bm{x_{2}}, \, {}^{sfc2u}\bm{x_{2}}, {}^{sfc2v}\bm{x_{2}} \in \mathbb{R}^{20550\times 2}}$, and $f$ is the tanh activation function. Then we use: 
\begin{equation}
{}^{sfc1}\bm{x_{2}} = \concat({}^{sfc1u}\bm{x_{2}}, {}^{sfc1v}\bm{x_{2}}) \in \mathbb{R}^{20550\times 4},
\end{equation}
as the input for layer 3-SFC1, and 
\begin{equation}
{}^{sfc2}\bm{x_{2}}=\concat({}^{sfc2u}\bm{x_{2}}, {}^{sfc2v}\bm{x_{2}})\in \mathbb{R}^{20550\times 4},
\end{equation}
as the input for layer~3-SFC2. 
See Table~\ref{ta:ArchituctureOfTwo} for a description of the following convolutional layers.

\paragraph{Layers 13, 14 and 15 (sparse layers)}\hfill\\ 
The vector ${{}^{sfc1}\bm{x_{12}}\in \mathbb{R}^{20550 \times 4}}$ is split into ${{}^{sfc1u}\bm{x_{12}}}, {{}^{sfc1v}\bm{x_{12}}\in \mathbb{R}^{20550 \times 2}}$, and similarly ${{}^{sfc2}\bm{x_{12}}}$ is split into ${}^{sfc2u}\bm{x_{12}}$, ${}^{sfc2v}\bm{x_{12}}\in \mathbb{R}^{20550 \times 2}$. In SFC ordering, ${{}^{sfc{\cal C}X}\bm{x_{12}^{+}}\in \mathbb{R}^{20550 \times 2}}$ are the neighbours (with an increase in SFC node number) of ${{}^{sfc{\cal C}X}\bm{x_{12}}}$, and ${{}^{sfc{\cal C}X}\bm{x_{12}^{-}}\in \mathbb{R}^{20550 \times 2}}$ are the neighbours (with an decrease in SFC node number) of ${{}^{sfc{\cal C}X}\bm{x_{12}}}$ in which ${\cal C}\in\{1,2\}$ is the SFC number and $X\in\{u,v\}$ indicated the velocity component. Thus: 
\begin{equation}%\label{eq:straight_line}
{}^{sfc1u}\bm{x_{13}}=
\sumchannels(
{}^{sfc1u}\bm{w_{13}}\odot{}^{sfc1u}\bm{x_{12}} + {}^{sfc1u}\bm{w_{13}^{+}}\odot{}^{sfc1u}\bm{x_{12}^{+}} + {}^{sfc1u}\bm{w_{13}^{-}}\odot{}^{sfc1u}\bm{x_{12}^{-}} 
),
\end{equation}
\begin{equation}
{}^{sfc1v}\bm{x_{13}}=
\sumchannels(
{}^{sfc1v}\bm{w_{13}}\odot{}^{sfc1v}\bm{x_{12}} + {}^{sfc1v}\bm{w_{13}^{+}}\odot{}^{sfc1v}\bm{x_{12}^{+}} + {}^{sfc1v}\bm{w_{13}^{-}}\odot{}^{sfc1v}\bm{x_{12}^{-}} 
),
\end{equation}
\begin{equation}
{}^{sfc2u}\bm{x_{13}}=
\sumchannels(
{}^{sfc2u}\bm{w_{13}}\odot{}^{sfc2u}\bm{x_{12}} + {}^{sfc2u}\bm{w_{13}^{+}}\odot{}^{sfc2u}\bm{x_{12}^{+}} + {}^{sfc2u}\bm{w_{13}^{-}}\odot{}^{sfc2u}\bm{x_{12}^{-}} 
),
\end{equation}
\begin{equation}
{}^{sfc2v}\bm{x_{13}}=
\sumchannels(
{}^{sfc2v}\bm{w_{13}}\odot{}^{sfc2v}\bm{x_{12}} + {}^{sfc2v}\bm{w_{13}^{+}}\odot{}^{sfc2v}\bm{x_{12}^{+}} + {}^{sfc2v}\bm{w_{13}^{-}}\odot{}^{sfc2v}\bm{x_{12}^{-}} 
),
\end{equation}
where ${}^{sfc1u}\bm{w_{13}}, \, {}^{sfc1v}\bm{w_{13}}, \, , {}^{sfc1u}\bm{w_{13}^{+}}, \, {}^{sfc1v}\bm{w_{13}^{+}}, \, {}^{sfc1u}\bm{w_{13}^{-}}, \, {}^{sfc1v}\bm{w_{13}^{-}}, \, {}^{sfc2u}\bm{w_{13}}, \, {}^{sfc2v}\bm{w_{13}}, \, {}^{sfc2u}\bm{w_{13}^{+}}, \\ {}^{sfc2v}\bm{w_{13}^{+}}, \, {}^{sfc2u}\bm{w_{13}^{-}}, \, {}^{sfc2v}\bm{w_{13}^{-}}\in \mathbb{R}^{20550\times 2}$ are the weight vectors 
and ${}^{sfc1u}\bm{x_{13}}$, ${}^{sfc1v}\bm{x_{13}}$,${}^{sfc2u}\bm{x_{13}}$, ${}^{sfc2v}\bm{x_{13}} \in \mathbb{R}^{20550\times 1}$ 
. Now convert ${}^{sfc1u}\bm{x_{13}}$ and ${}^{sfc1v}\bm{x_{13}}$ using the inverse of the first space-filling curve ordering to obtain ${}^{fem1u}\bm{x_{14}}$ and ${}^{fem1v}\bm{x_{14}}$, and convert ${}^{sfc2u}\bm{x_{13}}$ and ${}^{sfc2v}\bm{x_{13}}$ using the inverse of the second space-filling curve to obtain ${}^{fem2u}\bm{x_{14}}$ and ${}^{fem2v}\bm{x_{14}}$. Then the output 
${}^{fem}\bm{x_{15}} \in \mathbb{R}^{20550 \times 2}$
of the network is obtained from:
\begin{equation}%\label{eq:straight_line}
{}^{fem}\bm{x_{15}}=f(\concat({}^{fem1u}\bm{x_{14}}, {}^{fem1v}\bm{x_{14}}) + \concat({}^{fem2u}\bm{x_{14}}, {}^{fem2v}\bm{x_{14}}) + {}^{fem}\bm{b_{15}}), 
\end{equation}
in which $f$ is the tanh activation function, 
the biases are ${{}^{fem}\bm{b_{15}}\in \mathbb{R}^{20550 \times 2}}$ and the two channels in  ${{}^{fem}\bm{x_{15}}}$ contain the~$u$ and~$v$ velocity components.

\begin{table}[htbp]
\centering
\resizebox{\textwidth}{!}{%
\begin{tabular}{|c|c|c|c|c|c|c|c|}
\hline
layer ordering&input size \&  ordering&kernel size&channels&stride&padding&output size \& ordering&activation\\
\hline
1-FEM& (2, 20550, FEM)& 1 Identity & 2 & 1 & 0 & (2, 20550, SFC1) & Identity\\
\hline
\multicolumn{8}{|c|}{Split (2, 20550, SFC1) into two:(1, 20550, SFC1$u$), (1, 20550, SFC1$v$) }\\
\hline
2-SFC1$u$& (1, 20550, SFC1$u$)& 3 Variable ($3\times20550$)& 2 & 1 & 0 & (2, 20550, SFC1$u$) & tanh\\
\hline
2-SFC1v& (1, 20550, SFC1$v$)& 3 Variable ($3\times20550$)& 2 & 1 & 0 & (2, 20550, SFC1$v$) & tanh\\
\hline
\multicolumn{8}{|c|}{concatenate:((2, 20550, SFC1$u$), (2, 20550, SFC1$v$)) forming (4, 20550, SFC1) }\\
\hline
3-SFC1& (4, 20550, SFC1)& 32 & 8 & 4 & 16 & (8, 5138, SFC1) & tanh\\
\hline
4-SFC1& (8, 5138, SFC1)& 32 & 8 & 4 & 16 & (8, 1285, SFC1) & tanh\\
\hline
5-SFC1& (8, 1285, SFC1)& 32 & 8 & 4 & 16 & (8, 322, SFC1) & tanh\\
\hline
6-SFC1& (8, 322, SFC1)& 32 & 8 & 4 & 16 & (8, 81, SFC1) & tanh\\
\hline
1-FEM& (2, 20550, FEM)& 1 Identity & 2 & 1 & 0 & (2, 20550, SFC2) & Identity\\
\hline
\multicolumn{8}{|c|}{Split (2, 20550, SFC2) into two:(1, 20550, SFC2$u$), (1, 20550, SFC2$v$) }\\
\hline
2-SFC2$u$& (1, 20550)& 3 Variable ($3\times20550$)& 2 & 1 & 0 & (2, 20550, SFC2$u$) & tanh\\
\hline
2-SFC2$v$& (1, 20550)& 3 Variable ($3\times20550$)& 2 & 1 & 0 & (2, 20550, SFC2$v$) & tanh\\
\hline
\multicolumn{8}{|c|}{concatenate:((2, 20550, SFC2$u$), (2, 20550, SFC2$v$)) forming (4, 20550, SFC2) }\\
\hline
3-SFC2& (4, 20550, SFC2)& 32 & 8 & 4 & 16 & (8, 5138, SFC2) & tanh\\
\hline
4-SFC2& (8, 5138, SFC2)& 32 & 8 & 4 & 16 & (8, 1285, SFC2) & tanh\\
\hline
5-SFC2& (8, 1285, SFC2)& 32 & 8 & 4 & 16 & (8, 322, SFC2) & tanh\\
\hline
6-SFC2& (8, 322, SFC2)& 32 & 8 & 4 & 16 & (8, 81, SFC2) & tanh\\
\hline
\multicolumn{8}{|c|}{Flatten the output data of layer 6-SFC1 and 6-SFC2 to 1D - concatenate 2 sequences}\\
\hline
7& 1296 (= 8$\times$81 + 8$\times$81) & \multicolumn{4}{c|}{ } & 128& tanh\\
\hline
8 & 128 & \multicolumn{4}{c|}{ } & 1296& tanh\\
\hline
\multicolumn{8}{|c|}{Split the data into 2 sequence as the input of layer 9-SFC1 and 9-SFC2, convert from 1296 to (8, 81) }\\
\hline
9-SFC1& (8, 81, SFC1)& 32 & 8 & 4 & 15 & (8, 322, SFC1) & tanh\\
\hline
10-SFC1& (8, 322, SFC1)& 32 & 8 & 4 & 15 & (8, 1286, SFC1) & tanh\\
\hline
11-SFC1& (8, 1286, SFC1)& 32 & 8 & 4 & 16 & (8, 5140, SFC1) & tanh\\
\hline
12-SFC1& (8, 5140, SFC1)& 32 & 4 & 4 & 19 & (4, 20550, SFC1) & tanh\\
\hline
\multicolumn{8}{|c|}{Split (4, 20550, SFC1) into two:(2, 20550, SFC1$u$), (2, 20550, SFC1$v$) }\\
\hline
13-SFC1$u$& (2, 20550, SFC1$u$)& 3 Variable ($3\times20550$)& 2 & 1 & 0 & (1, 20550, SFC1$u$) & Identity\\
\hline
13-SFC1$v$& (2, 20550, SFC1$v$)& 3 Variable ($3\times20550$)& 2 & 1 & 0 & (1, 20550, SFC1$v$) & Identity\\
\hline
14-SFC1$u$& (1, 20550, SFC1$u$)& 1 Identity & 1 & 1 & 0 & (1, 20550, FEM1$u$) & Identity\\
\hline
14-SFC1$v$& (1, 20550, SFC1$v$)& 1 Identity & 1 & 1 & 0 & (1, 20550, FEM1$v$) & Identity\\
\hline
9-SFC2& (8, 81, SFC2)& 32 & 8 & 4 & 15 & (8, 322, SFC2) & tanh\\
\hline
10-SFC2& (8, 322, SFC2)& 32 & 8 & 4 & 15 & (8, 1286, SFC2) & tanh\\
\hline
11-SFC2& (8, 1286, SFC2)& 32 & 8 & 4 & 16 & (8, 5140, SFC2) & tanh\\
\hline
12-SFC2& (8, 5140, SFC2)& 32 & 4 & 4 & 19 & (4, 20550, SFC2) & tanh\\
\hline
\multicolumn{8}{|c|}{Split (4, 20550, SFC2) into two:(2, 20550, SFC2$u$), (2, 20550, SFC2$v$) }\\
\hline
13-SFC2$u$& (2, 20550, SFC2$u$)& 3 Variable ($3\times20550$)& 1 & 1 & 0 & (1, 20550, SFC2$u$) & Identity\\
\hline
13-SFC2$v$& (2, 20550, SFC2$v$)& 3 Variable ($3\times20550$)& 1 & 1 & 0 & (1 20550, SFC2$v$) & Identity\\
\hline
14-SFC2$u$& (1, 20550, SFC2$u$)& 1 Identity & 1 & 1 & 0 & (1, 20550, FEM2$u$) & Identity\\
\hline
14-SFC2$v$& (1, 20550, SFC2$v$)& 1 Identity & 1 & 1 & 0 & (1, 20550, FEM2$v$) & Identity\\
\hline
15-FEM & ($1\times 4$, 20550)& 1 Identity ($\concat$(FEM1$u$,FEM1$v$)+$\concat$(FEM2$u$,FEM2$v$)) & 2 & 1 & 0 & (2, 20550, FEM) & tanh\\
\hline
\end{tabular}
}
\caption{The architecture of the SFC-based CNN with two space-filling curves with nearest-neighbour smoothing for flow past a cylinder. The number of neurons in the bottleneck layer can be adjusted from 128 to any suitable value. Whenever a non-identity activation term is used a bias is added to each neuron.  }
\label{two-SFC-CNN-cylinder}
\end{table}

%\clearpage
\section{Results}\label{sec:results}

Three examples are used to demonstrate the potential of the method proposed in this paper, of applying convolutional networks to data held on any mesh (structured grids or unstructured meshes) by using space-filling curves. Two of the test cases consist of structured grid data; the first data set represents advection of a square wave and the second represents advection of a Gaussian function. The third example uses a data set consisting of solutions of 2D flow past a cylinder solved on an unstructured mesh. 

%For the structured data, we construct (1)~an autoencoder based on one space-filling curve; (2)~an autoencoder based on two space-filling curves; (3)~an autoencoder based on two space-filling curves with smoothing by the nearest neighbours. These results are compared with singular value decomposition and a classical 2D convolutional autoencoder. For the unstructured data, we construct (1)~an autoencoder based on one space-filling curve with nearest-neighbour smoothing; (2)~an autoencoder based on two space-filling curves with nearest-neighbour smoothing. Comparison cannot be made with a classical 2D convolutional autoencoder in this case, as the solutions lie on an unstructured mesh, but we do compare with results from an SVD. 

%The autoencoder lends itself to comparison with singular value decomposition, another method that is often used for compressing data. 

Section~\ref{sec:results_analysis} describes how results from an SVD will be compared with the SFC-based autoencoders, and describes the error measure used to compare the results. Section~\ref{sec:remaining_hyper-parameters} defines some of the remaining hyper-parameters not defined in Section~\ref{sec:architectures}.  Section~\ref{sec:sfc_cnn_structured_applications} explains how the two structured-grid data sets are generated and then presents results for three autoencoders based on space-filling curves, a classical 2D autoencoder and an SVD.  Section~\ref{sec:sfc_cnn_unstructured_applications} describes generation of the results for flow past a cylinder, which make up the unstructured mesh data set. Results are presented for the two autoencoders based on space-filling curves and an SVD.  

\subsection{Analysis of results} \label{sec:results_analysis}
\subsubsection{Singular Value Decomposition}\label{sec:analysis_SVD}
The results from the convolutional autoencoders developed in this paper are compared with results from singular value decomposition (SVD). For a matrix~$\bm{M}$, where each column corresponds to a particular solution or example from the data set, and each row corresponds to a particular node or cell, the SVD is defined as
\begin{equation}
\label{eq:SVD_definition}
\bm{M}=\bm{U} \bm{\Sigma} \bm{V^{*}}\,,
\end{equation}
where ${\bm{\Sigma}}$ is a diagonal matrix containing the singular values of $\bm{M}$ given in descending order, $\bm{U}$ and $\bm{V}$ are the left- and right-singular vectors respectively, and the asterisk denotes the conjugate transpose. The square of each singular value indicates how much information is contained in the corresponding mode or basis function (column of $\bm{U}$). A low-rank approximation of $\bm{M}$ can be formed by retaining only the $N_{\Sigma}$ largest singular values, setting the remaining singular values to zero, and recalculating the product:
\begin{equation}
%\label{eq:straight_line}
\bm{\widetilde{M}}=\bm{U} \bm{\widetilde{\,\Sigma\,}} \bm{V^{*}}\,,
\end{equation}
where the only non-zero terms in ${\bm{\widetilde{\,\Sigma\,}}}$ are the $N_\Sigma$ largest singular values from ${\bm{\Sigma}}$. Comparison is made between convolutional autoencoders which have a latent space of dimension~$N_\Sigma$ and SVDs which have been truncated to $N_\Sigma$ singular values.

\subsubsection{Measuring the error}\label{sec:analysis_error}
To evaluate the error in the various autoencoder and SVD approaches, the mean square error is used: 
\begin{equation}\label{eq:MSE}
MSE = \frac{ \sum_{k=1}^{E}\sum_{i=1}^{N} (\varepsilon_{i}^k)^{2}}{N\, E }\,,
\end{equation}
in which $N$ is the number of input (or output variables), the number of nodes or cells for instance, $k$ is a particular solution or example taken from the data set, $E$ is the number of solutions or examples under consideration, and $\epsilon_i^k$ is the difference between the original solution and its approximation. For the autoencoder, the approximation is formed by encoding and decoding the solutions or examples. For the SVD, the solutions are assembled into a matrix to which SVD is applied. The truncation is performed and the low-rank approximation obtained. For the examples using structured grid data, the examples or solutions are scalar fields so Equation~\eqref{eq:MSE} can be applied directly. For the data set with examples on an unstructured mesh, the solutions are vector fields (velocities), so in this case, $\varepsilon_i^k$ is defined to be the error in speed at the $i$th node for the $k$th example. 

\subsection{Hyper-parameters of the autoencoders} \label{sec:remaining_hyper-parameters}
For each particular autoencoder, some of the hyper-parameters are specified in Section~\ref{sec:architectures}, including kernel size, number of channels, stride, padding, activation function, number of layers and number of neurons per layer. Other hyper-parameters are loss function, batch size, optimiser, learning rate and number of epochs, which are given here, except for the number of epochs which is stated in the results section where appropriate. Through tuning these hyper-parameters, we have optimised the neural networks within this work. Unless explicitly stated, it can be assumed that the values given in Table~\ref{hyper} are used for the SFC-based autoencoders.  
%Figure \ref{fig:Autoencoder} shows the basic architecture of the  two space-filling 
%curve autoencoder. 

\begin{table}[htbp]
\centering
\begin{tabular}{lcccccc}
\toprule
& & batch size & loss function & optimiser &  learning rate  & activation function \\
\midrule
structured data &\quad & 64 & MSE &  Adam & 0.0001 & ReLU \\
unstructured data & &16 & MSE &  Adam & 0.0001 & tanh \\
\bottomrule
\end{tabular}
\caption{The hyper-parameters used when training the SFC-based convolutional autoencoders for the structured and unstructured data sets. Mean square error is written as~MSE; and ReLU represents the rectified linear unit.}
\label{hyper}
\end{table}

%\begin{table}[htbp]
%\centering
%\begin{tabular}{lcccccc}
%\toprule
%& & batch size & kernel size & optimiser &  learning rate  & %activation function  \\
%\midrule
%structured data &\quad & 64 & 32 &  Adam & 0.0001 & ReLU \\
%unstructured data & &16 & 32 &  Adam & 0.0001 & tanh \\
%\bottomrule
%\end{tabular}
%\caption{The hyper-parameters used for training the SFC-based convolutional autoencoders.}
%\label{hyper}
%\end{table}

The kernel size for the 2D classical autoencoder is ${5\times 5=25}$ with a stride of ${2\times 2}$. To make the~1D SFC-based autoencoders approximately analagous to this, a kernel size of 32 is used with a stride of~4. A number of activation functions were experimented with, and ReLU and tanh performed the best. Here we use ReLU for the structured data and tanh for the unstructured data. All the input data for the autoencoders are normalised between ${[0,1]}$ for the structured data and  ${[-1,1]}$ for the unstructured data. The unstructured data has a larger computational cost associated with it, so the batch size was reduced from~64 to~16. The loss function chosen is the mean square error~(MSE) given in Equation~\eqref{eq:MSE}, and errors calculated using this formula are based on the normalised data values.

%The auto-encoder here has the architecture shown in figure \ref{fig:Autoencoder}.  It has a several convolutional layers starting from, in 2D, the vectors of u and v velocity components. The auto-encoder has in its centre fully connected multi-layered percepron layers which connect to the convolutional layers and produce the required bottle neck or compressed number of variables or latent variables.  

%If a single or multiple concentration/temperature/density fields, that are advected with the fluid, are also solved for, then there is a choice on how to compress these. One could compress them uncoupled to the velocity field and even seperate to one another (if there are multiple concentration fields) or one could compress them together with the velocity field in the hope that fewer compressed variables would be needed overall. 

\subsection{Structured Grid Applications}
\label{sec:sfc_cnn_structured_applications}
The generation of the data sets for the square wave and Gaussian function test cases is described. These data sets represent two extremes, from abruptly changing fields to smoothly changing fields. The performance of the new autoencoders based on space-filling curve ordering is analysed and comparisons are made with a classical 2D convolutional autoencoder (CAE) and SVD. 

\subsubsection{Generating square wave data}

In order to generate data that represents advection of a square wave, the time-dependent 2D advection equation is solved:
\begin{equation}\label{eq:advection}
\frac{\partial c}{\partial t} + u \frac{\partial c}{\partial x} + v \frac{\partial c}{\partial y} = 0\,,
\end{equation}
in which ${t\in[0,T]}$~represents time, $x$~and $y$~are the coordinates, $c(x,y,t)$ is the concentration field and the advection velocity vector is given by $(u, v)^T$. The discretisation method used to solve this equation is upwind differencing in space and the backward Euler method in time. The initial condition is given as: 
\begin{equation}
\label{eq:square_wave_initial_condition}
c(x, y, 0) = \left\{
             \begin{array}{lr}
             1, \;\;\; (x,y) \in [x_0, x_0 + 0.5] \times [y_0, y_0 + 0.5]; &  \\

             0, \;\;\; \text{otherwise}, &  
             \end{array}
\right.
\end{equation}
and defined on the domain ${[0,L]\times[0,L]}$. For this problem, $L=3$, $T=0.3$, $u=1$, $v=1$. The domain is discretised using a uniform ${128\times 128}$ grid, the time step is ${\Delta t = 10^{-2}}$ and the variables $x_0$~and $y_0$ are chosen randomly from the interval~$[0,2]$. For this test case, the number of initial conditions created is ${N_s=512}$, each generating ${N_t=30}$~time levels, so the total number of solutions or examples in the data set is~15360. % %Therefore, the set of solutions has the form ${\chi \in [0,1]^{N_x \times N_y \times N_t \times N_s}}$. 
The data is divided randomly into three parts according to the proportion 6:2:2, for training, validation and testing.

\subsubsection{Generating the Gaussian data}
To represent advection of a Gaussian function, the following profile is simply located at different points in the domain. The form of the Gaussian is given by
\begin{equation}\label{eq:Gaussian}
g(x,y) =  \exp \left( -\frac{(x-x_c)^2+(y-y_c)^2}{2 \sigma^2}\right)
\end{equation}
where ${(x_c, y_c)}$ represents the centre of the Gaussian function, whose values are randomly sampled from the domain, which is discretised with a structured ${128\times128}$ grid. The parameter ${\sigma}$, which controls the width of the curve, is uniformly randomly sampled from the interval $[10, 20]$. The data set consists of a total ${N_s=15360}$ examples. % so the data set has the form ${\chi \in [0,1]^{N_x \times N_y \times  N_s}}$. 

\subsubsection{Results for advection of a square wave}

%\subsubsection{Comparison of CNN performance}

\begin{figure}[htbp]
\centering
\includegraphics[width=15cm]{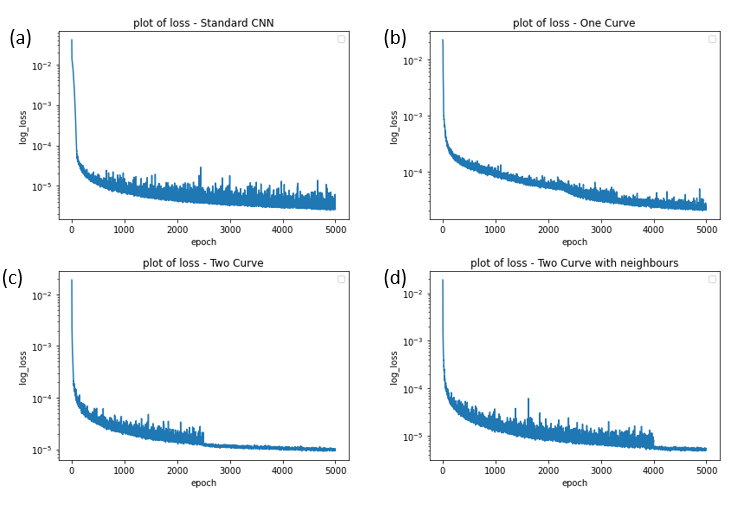}
\caption{Comparison of the training losses of the four convolutional networks: (a)~Classical CAE; (b) CAE-SFC; (c) CAE-2SFC; (d) CAE-2SFC-NN. 
The validation losses mostly overlaid the training losses, and are thus not plotted. The square wave data was compressed to 16~variables and the networks were trained over 5000~epochs.}
\label{fig:GeneralComparison-step}
\end{figure}

In order to assess their relative capabilities, we compare three SFC-based convolutional autoencoders, a classical 2D convolutional autoencoder and the SVD, all applied to the square wave data set. The autoencoders, their abbreviations and the section in which their architesctures are described, are listed in Table~\ref{tab:CAEs_structured_data}. %
\begin{table}[htbp]
\centering
\begin{tabular}{llp{0.5\textwidth}}
\toprule
abbreviation & architecture & description  \\
\midrule %\ref{one-SFC-CNN-cylinder} \ref{two-SFC-CNN-cylinder}  
CAE         & Table~\ref{classicCNN}, Section~\ref{sec:structured_CAE}& classical convolutional autoencoder \\[2mm]
CAE-SFC     & Table~\ref{ta:ArchituctureOfOne}, Section~\ref{sec:structured_CAE_SFC} & a convolutional autoencoder based on ordering from one SFC\\[2mm]
CAE-2SFC    & Table~\ref{ta:ArchituctureOfTwo}, Section~\ref{sec:structured_CAE_2SFC} & a convolutional autoencoder based on ordering from two SFCs \\[2mm]
CAE-2SFC-NN & Section~\ref{sec:structured_CAE_2SFC_NN} & a convolutional autoencoder based on ordering from two SFCs and nearest-neighbour smoothing\\
\bottomrule
\end{tabular}
\caption{Abbreviations of the autoencoders applied to the structured data sets.}
\label{tab:CAEs_structured_data}
\end{table}
The convolutional autoencoders all compress to 16~latent variables having been trained over 5000~epochs and the SVD truncates to 16~variables. The so-called losses, defined here as the mean square error between the inputs and outputs of the autoencoder, see Equation~\eqref{eq:MSE}, are shown in Figure~\ref{fig:GeneralComparison-step} and Table~\ref{tabGeneral}, the latter includes the truncation error of the SVD. Both figure and table show that at 5000 epochs the losses of the classical 2D CAE and CAE-2SFC-NN are smaller than those of the CAE-SFC and CAE-2SFC autoencoders. All the autoencoders have a lower mean square error than the SVD. 

\begin{table}[htbp]
\centering
%\resizebox{\textwidth}{!}{%
\begin{tabular}{ccccc}
\toprule
2D CAE & CAE-SFC & CAE-2SFC & CAE-2SFC-NN & SVD\\
\midrule
${ 2.8487 \times 10^{-6}}$&${ 2.1742 \times 10^{-5} }$&${ 1.0125 \times 10^{-5} }$&${ 5.2978 \times 10^{-6} }$&${ 5.6460 \times 10^{-3}}$\\
\bottomrule
\end{tabular}
%}
\caption{Losses for the validation data of the autoencoders and truncation error of the SVD after application to the square wave data set. The training and test losses are similar to these. For training, 5000 epochs were used and compression is to 16~variables.}
\label{tabGeneral}  
\end{table}

Next, we investigate the effect of the dimension of the latent space by training two autoencoders to compress to both~16 and~8 latent variables. The best-performing SFC-based autoencoder from the previous results (CAE-SFC2-NN) is compared with the classical 2D CAE and the SVD. In training, 2000~epochs were used, and the losses for the training, validation and test data sets can be seen in Table~\ref{ta:StandNeighStep}. 
From this, it can be seen that the classical CAE performs a little better than CAE-2SFC-NN when compressing to 16 variables; both perform better than the SVD. However, the classical autoencoder struggles to converge when compressing to 8~variables. The CAE-2SFC-NN network performs well for both levels of compression presented here. 

\begin{table}[htbp]
\centering
%\resizebox{\textwidth}{!}{%
\begin{tabular}{lcccc}
\toprule
& 2D CAE & CAE-2SFC-NN & SVD & compression \\
\midrule
Training data &${ 4.8715 \times 10^{-6}}$&${ 1.6845 \times 10^{-5} }$&${ 5.6699 \times 10^{-3} }$ &  \\
Validation data &${ 5.0757 \times 10^{-6}}$&${ 1.3293 \times 10^{-5} }$&${ 5.6470 \times 10^{-3} }$ & 16 \\
Test data &${ 5.0195 \times 10^{-6}}$&${ 1.3419 \times 10^{-5} }$&${ 5.6459 \times 10^{-3} }$ &\\[1mm]
\midrule
Training data &${ 2.0209 \times 10^{-2}}$&${ 1.2142 \times 10^{-5} }$&${ 1.0501 \times 10^{-2} }$ & \\
Validation data &${ 2.0206 \times 10^{-2}}$&${ 1.3336 \times 10^{-5} }$&${ 1.0455 \times 10^{-2} }$ & 8\\ 
Test data &${ 2.0211 \times 10^{-2}}$&${ 1.3625 \times 10^{-5} }$&${ 1.0410 \times 10^{-2} }$ &  \\[1mm] 
\bottomrule
\end{tabular}
\caption{Comparing truncation error in the SVD with the losses after 2000~epochs of the classical 2D CAE and the SFC-based autoencoder with two SFCs and nearest-neighbour smoothing (CAE-2SFC-NN). The square wave data set on a structured ${128\times128}$ grid was used. Results are shown for compression to both 16 and~8 variables. }
\label{ta:StandNeighStep}
\end{table}

\begin{figure}[htbp]
\centering
\includegraphics[width=13cm]{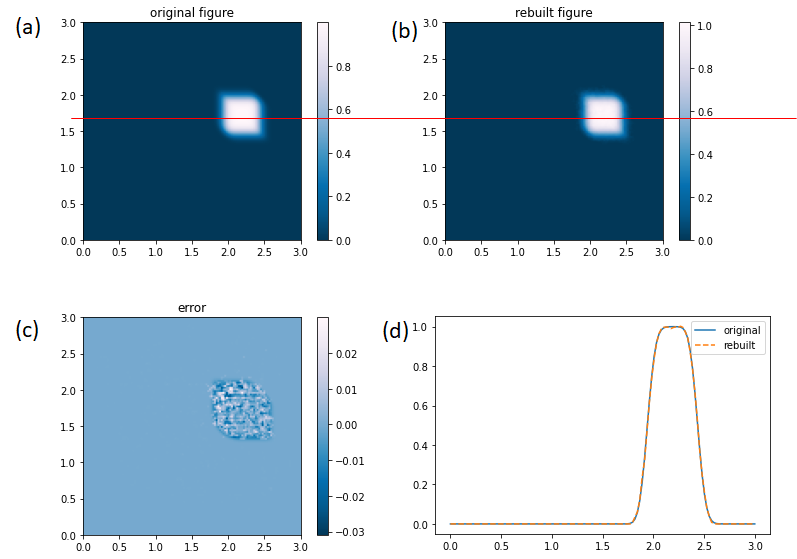}
\caption{For one example %from the  \textcolor{red}{training, validation, test data?} 
of the advection of a square wave (a)~original solution; (b)~solution after applying the autoencoder CAE-2SFC-NN; (c)~the point-wise error; (d)~the original solution and the solution after applying the autoencoder taken at a height of $y=1.7$.}
\label{fig:ErrorNeighStep16}
\end{figure}

%\textcolor{red}{\textit{Yuling, was this example taken from the training, validation, or test data?}} 
For one example, %from the \textcolor{red}{training/validation/test data set},
Figure~\ref{fig:ErrorNeighStep16} shows the performance of the autoencoder based on two SFCs with nearest-neighbour smoothing (CAE-2SFC-NN), which compresses to 16~latent variables and was trained in 5000~epochs. It can be seen that the network performs well and accurately reproduces the profile of the demanding square wave function. Plot~(a) and~(b) show the square wave solution before and after the application of the autoencoder respectively. The red line indicates a height of 1.7, at which plot~(d) compares the original solution with the output of the autoencoder. Plot~(c) shows the pointwise error, which, has an absolute value of, at most,~3\%. 

We also used a finer grid, ${256\times 256}$, to compare the classical 2D CAE and the CAE-2SFC-NN autoencoder for the square wave data set, see Table~\ref{DenseStep} for the losses.  The classical CAE does not converge on this finer grid, with a loss similar in magnitude to the truncation error of the SVD, but the CAE-2SFC-NN autoencoder converges well. 

\begin{table}[htbp]
\centering
%\resizebox{\textwidth}{!}{%
\begin{tabular}{lccc}
\toprule
  & classical CAE & CAE-2SFC-NN & SVD\\
\midrule
Training data &${ 2.5308 \times 10^{-2}}$&${ 1.1520 \times 10^{-4} }$& ${ 1.0180 \times 10^{-2} }$\\
Validation data &${ 2.5308 \times 10^{-2}}$&${ 1.1783 \times 10^{-4} }$& ${ 1.0224 \times 10^{-2} }$\\
Test data &${ 2.5293 \times 10^{-2}}$&${ 1.1813 \times 10^{-4} }$& ${ 1.0152 \times 10^{-2} }$\\
\bottomrule
\end{tabular}
%}
\caption{Losses for the autoencoders and truncation error of the SVD when applied to the square wave data set 
on the finer grid, ${256\times256}$. Results are shown for compression to 16~variables, networks were trained for 1000~epochs. }
\label{DenseStep}
\end{table}

\subsubsection{Results for advection of a Gaussian function}
%\textcolor{red}{\textit{Yuling: do we have losses for this?}}\\
In this section, the architectures are the same as those used in the previous section for the advection of a square wave, as the grid is the same size. First, the performance of the autoencoder based on two SFCs with nearest-neighbour smoothing is compared with the classical CAE and the SVD. The autoencoders are both compressed to~16 and~8 variables, and trained for 2000~epochs. The SVD is truncated to 16 and 8~variables. The losses are given in Table~\ref{ta:StandNeighGauss}, from which it can be seen that the autoencoder based on two SFCs with nearest-neighbour smoothing performs a little better than the classical autoencoder, and both are considerably better than the SVD when compressing to 16~variables. However, the classical convolutional autoencoder cannot converge when compressing to 8~variables. 

\begin{table}[htbp]
\centering
%\resizebox{\textwidth}{!}{%
\begin{tabular}{lcccc}
\toprule
& 2D CAE & CAE-2SFC-NN & SVD & compression\\
\midrule
Training data & ${ 1.8429 \times 10^{-5}}$&${ 5.3081 \times 10^{-6} }$&${ 3.4123 \times 10^{-3} }$ & \\
Validation data & ${ 2.2369 \times 10^{-5}}$&${ 7.3057 \times 10^{-6} }$&${ 3.3651 \times 10^{-3} }$ & 16 \\
Test data & ${ 2.2246 \times 10^{-5}}$&${ 7.2510 \times 10^{-6} }$&${ 3.3677 \times 10^{-3} }$ & \\[1mm]
\midrule
Training data &${ 3.8499 \times 10^{-2}}$&${ 4.7776 \times 10^{-6} }$&${ 9.8718 \times 10^{-3} }$ & \\
Validation data &${ 3.8795 \times 10^{-2}}$&${ 6.0837 \times 10^{-6} }$&${ 9.8222 \times 10^{-3} }$& 8 \\ 
Test data &${ 3.8537 \times 10^{-2}}$&${ 6.0445 \times 10^{-6} }$&${ 9.8445 \times 10^{-3} }$& \\[1mm]
\bottomrule
\end{tabular}
%}
\caption{Comparison of the losses of the classical 2D CAE, the autoencoder based on 2 SFCs with nearest-neighbour smoothing (CAE-2SFC-NN) and truncation error of the SVD when applied to the Gaussian function data set, with compression to~16 and 8~variables. The networks were trained for 2000~epochs.}
\label{ta:StandNeighGauss}
\end{table}

\begin{figure}[htbp]
\centering
\includegraphics[width=13cm]{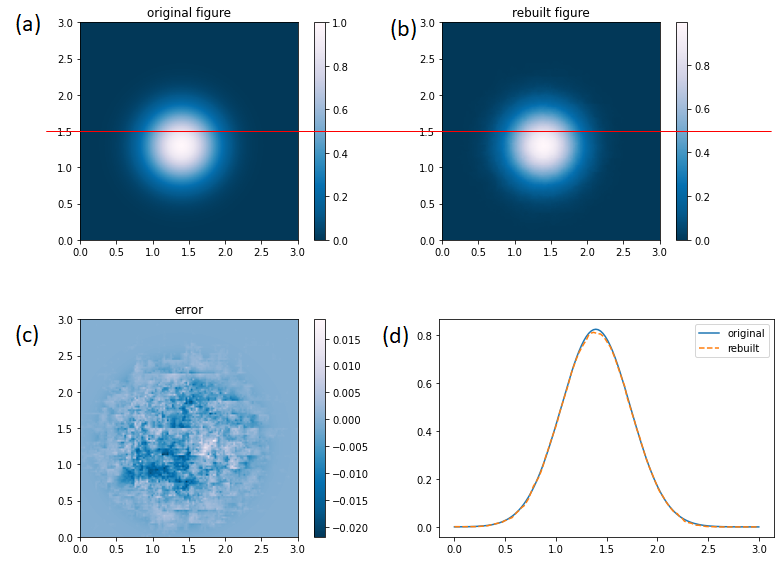}
\caption{For one solution %from the \textcolor{red}{training/ validation/ test data set} of %
for the advection of a Gaussian function (a)~original solution; (b)~solution after applying the autoencoder (CAE-2SFC-NN); (c)~the point-wise error; (d)~the original solution and the solution after applying the autoencoder taken at a height of $y=1.5$.}
\label{fig:ErrorNeighGauss16}
\end{figure}

%\textcolor{red}{\textit{Yuling, was this example taken from the training, validation, or test data?}} 
For one solution, %from the  \textcolor{red}{training/validation/test} data set of the Gaussian function,
Figure~\ref{fig:ErrorNeighGauss16} shows the solution (a)~before and (b)~after being passed through the autoencoder (CAE-2SFC-NN). %compressed to 16 variables and trained over 2000 epochs. 
It can be seen that that the model performs well and accurately reproduces the smooth Gaussian function for this particular example. Plot~(c) shows the pointwise error, which, at most, has an absolute value of 2\%. Plot~(d) compares the original solution with the output of the autoencoder at a height of~1.5. The profiles shown are in very close agreement. 

%\subsubsection{Comparison between classical CNNs and two curves with neighbours SFC-CNN on Gaussian function dataset in 256${\times}$256 structured mesh (see table~\ref{DenseGauss}):}

The SFC-based autoencoder with two SFCs and nearest-neighbour smoothing is also compared with the classical CAE and SVD for a data set generated on a finer grid of ${256\times256}$. The losses are shown in Table~\ref{DenseGauss}. Here we see that, as for the square wave data set (Table~\ref{DenseStep}), the classical CAE does not converge for the ${256\times256}$ Gaussian function data set but the new CAE-2SFC-NN does converge. The data was compressed to 16~variables and the networks trained over 1000~epochs. We thus conclude that the new SFC-based autoencoder has comparable (possibly slightly better) accuracy to the classical 2D convolutional autoencoder; it outperforms the SVD (for these low dimensional spaces); it converges for a wider range of problems; but most importantly, the SFC-based autoencoders can be applied to data from unstructured meshes.

\begin{table}[htbp]
\centering
%\resizebox{\textwidth}{!}{%
\begin{tabular}{lccc}
\toprule
 & 2D CAE &  CAE-2SFC-NN & SVD\\
\midrule
Training data &${ 2.6126 \times 10^{-2}}$&${ 1.0321 \times 10^{-5} }$&${ 3.9307 \times 10^{-3} }$\\
Validation data &${ 2.6239 \times 10^{-2}}$&${ 1.8860 \times 10^{-5} }$&${ 3.8783 \times 10^{-3} }$\\
Test data &${ 2.6734 \times 10^{-2}}$&${ 1.7490 \times 10^{-5} }$&${ 3.7960 \times 10^{-3} }$\\
\bottomrule
\end{tabular}
%}
\caption{Truncation error for the SVD and losses for the classical 2D CAE and the autoencoder based on 2 SFCs with nearest-neighbour smoothing for the Gaussian function data set on the ${256\times 256}$ grid. Results are shown for compression to 16~variables, networks trained for 1000~epochs.}
\label{DenseGauss}
\end{table}

\subsection{Unstructured mesh test case}\label{sec:sfc_cnn_unstructured_applications}
A data set consisting of solutions for 2D flow past a cylinder is created. These results lie on an unstructured mesh which will reveal the effectiveness of the method for this type of data. The SFC-based approach is compared with the SVD.

\subsubsection{Results for 2D flow past a cylinder}
Using the conservation laws, the following system of partial differential equations governing the motion of an incompressible fluid is obtained: 
\begin{eqnarray}
 \bm{\nabla} \cdot \bm{u} & = & 0\,, \label{conmass} \\
\frac{\partial}{\partial t}\left(\rho\bm{u}\right) + \bm{\nabla} \cdot\left(\rho\bm{u}\otimes\bm{u}\right) -\bm{\nabla}\cdot\bm{\tau} & =& -\bm{\nabla}p\,,
\label{conmom}
\end{eqnarray}
where $\rho$ is the density (assumed constant), $\bm{u}$~is the velocity vector, $\bm{\tau}$~contains the viscous terms associated with an isotropic Newtonian fluid, $p$~represents the non-hydrostatic pressure, $t$~is time and the gradient operator $\bm{\nabla}$ is defined as 
\begin{equation}
\bm{\nabla} = \left(\frac{\partial}{\partial x}\,,\ \frac{\partial}{\partial y}\right)^T \,.
\end{equation}
This system of equations is solved as outlined in~\cite{Xie_2016}. For the discretisation, a linear triangular element is adopted with a discontinuous Galerkin discretisation of the velocities and a continuous Galerkin representation of the pressure, often referred to as the P1DG-P1 element. Crank Nicolson is used to discretise in time. Only velocity variables are needed to train the networks, as these fully describe incompressible flow. The Reynolds number for this problem is  
\begin{equation}
    Re=\frac{\rho\, U L }{\mu}=3900
\end{equation} 
in which the inlet velocity is constant,  $U=$\SI{0.039}{\metre\per\second}, the density has value $\rho=$\SI{1000}{\kilogram\per\cubic\metre} and the diameter of the cylinder is $L=\SI{0.1}{m}$. Thus the dynamic viscosity is $\mu=$\SI{e-3}{\kilo\gram\per\metre\per\second}. 
The domain measures \SI{2.2}{m} horizontally, \SI{0.41}{m} vertically, and the centre of the cylinder is located \SI{0.2}{m} from the left boundary on the horizontal centreline of the domain. No slip and no normal flow boundary conditions are applied on the upper and lower walls, and on the cylinder. Zero shear and zero normal stress are applied at the outlet (the right boundary of the domain). In the following results, time is given in seconds, and speeds and velocities are in metres per second.

The data set, formed from the solutions of the above problem, consists of ${N_s = 1000}$ snapshots. Each snapshot has ${N = 20550}$ nodes, and each node has ${N_{uv} = 2}$ features 
(the two velocity components,~$u$ and~$v$, in this two-dimensional problem). %This problem was solved on an unstructured mesh.%Therefore, the dataset has the form ${\chi \in \mathbb{R}^{N_s \times N \times N_{uv}} }$. We shuffle the dataset in the dimension ${N_s}$,  and divide the dataset into three parts with the proportion 8:1:1 for training, validation and testing (has no influence on the training) respectively.
The data set is divided randomly (in the $N_s$ dimension) into three parts according to the proportion 8:1:1 for training, validation and testing.  
%We shuffle the dataset in the dimension ${N_s}$,  and divide the dataset into three parts with the proportion 8:1:1 for training, validation and testing (has no influence on the training) respectively. 

First, in order to determine the effectiveness of multiple SFCs, comparison is made between two SFC-based convolutional autoencoders: one using one space-filling curve and the other using two space-filling curves. Both have nearest-neighbour smoothing.  Following this, the performance of the CAE-2SFC-NN network is investigated for different dimensions of latent space. The architectures of the autoencoders are described in the sections and tables listed in Table~\ref{tab:CAEs_unstructured_data}.

\begin{table}[htbp]
\centering
\begin{tabular}{llp{0.5\textwidth}}
\toprule
abbreviation & architecture & description  \\
\midrule 
CAE-SFC-NN     & Table~\ref{one-SFC-CNN-cylinder}, Section~\ref{sec:unstructured_CAE_SFC_NN} & a convolutional autoencoder based on ordering from one SFC with nearest-neighbour smoothing\\[2mm]
CAE-2SFC-NN & Table~\ref{two-SFC-CNN-cylinder}, Section~\ref{sec:unstructured_CAE_2SFC_NN}  & a convolutional autoencoder based on ordering from two SFCs with nearest-neighbour smoothing\\
\bottomrule
\end{tabular}
\caption{Abbreviations of the autoencoders applied to the unstructured data set.}
\label{tab:CAEs_unstructured_data}
\end{table}

\subsubsection{Comparison between two SFC-based autoencoders} 
\begin{figure}[htbp]
\centering
\includegraphics[width=18cm, trim=20mm 0mm 20mm 0mm,clip]{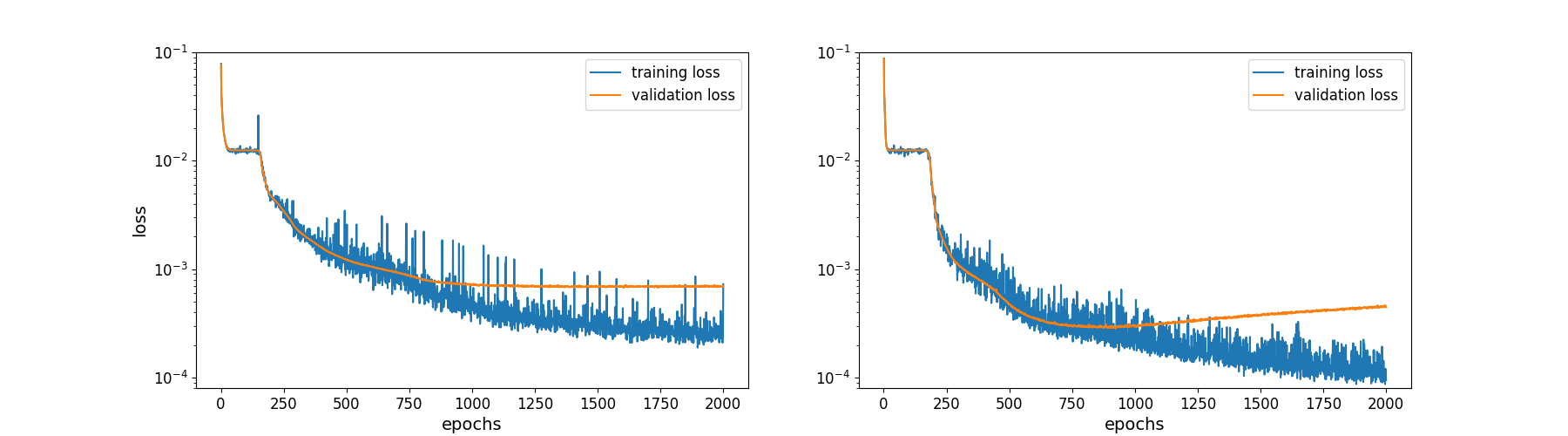}
\caption{Training and validation loss for two autoencoders trained on unstructured mesh data with compression to 128 variables: one space-filling curve with nearest-neighbour smoothing (CAE-SFC-NN) (left); two space-filling curves with nearest-neighbour smoothing (CAE-2SFC-NN) (right).}
\label{fig:OneTwoUnstructured}
\end{figure}

\begin{figure}[htbp]
\centering
\includegraphics[width=8cm]{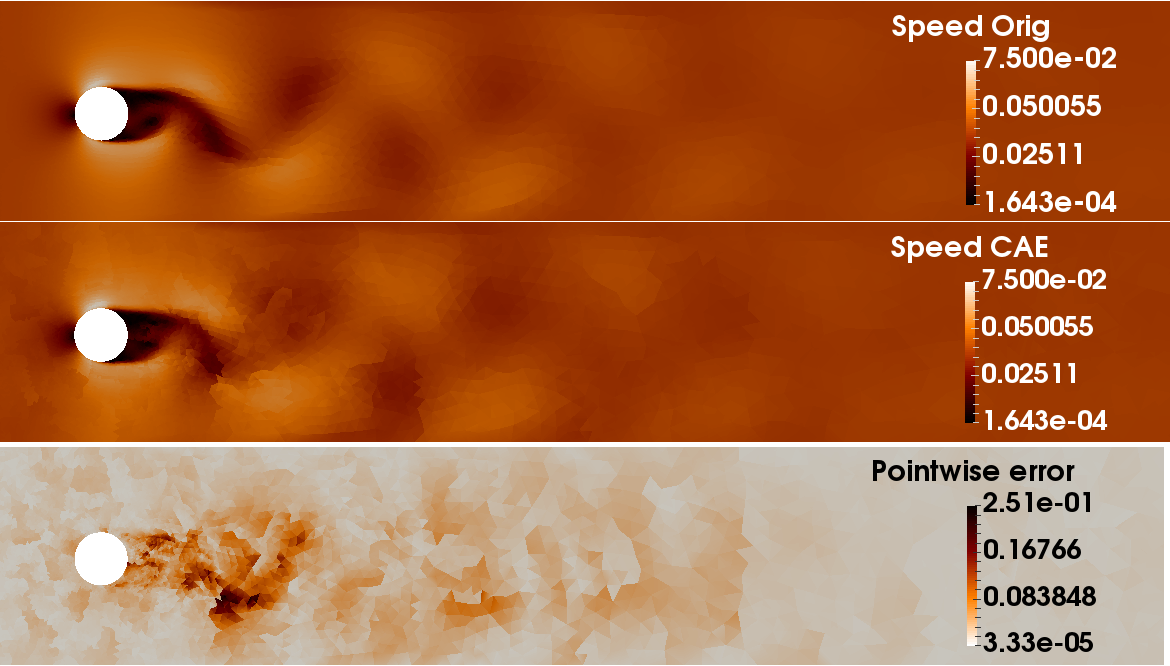}~~%
\includegraphics[width=8cm]{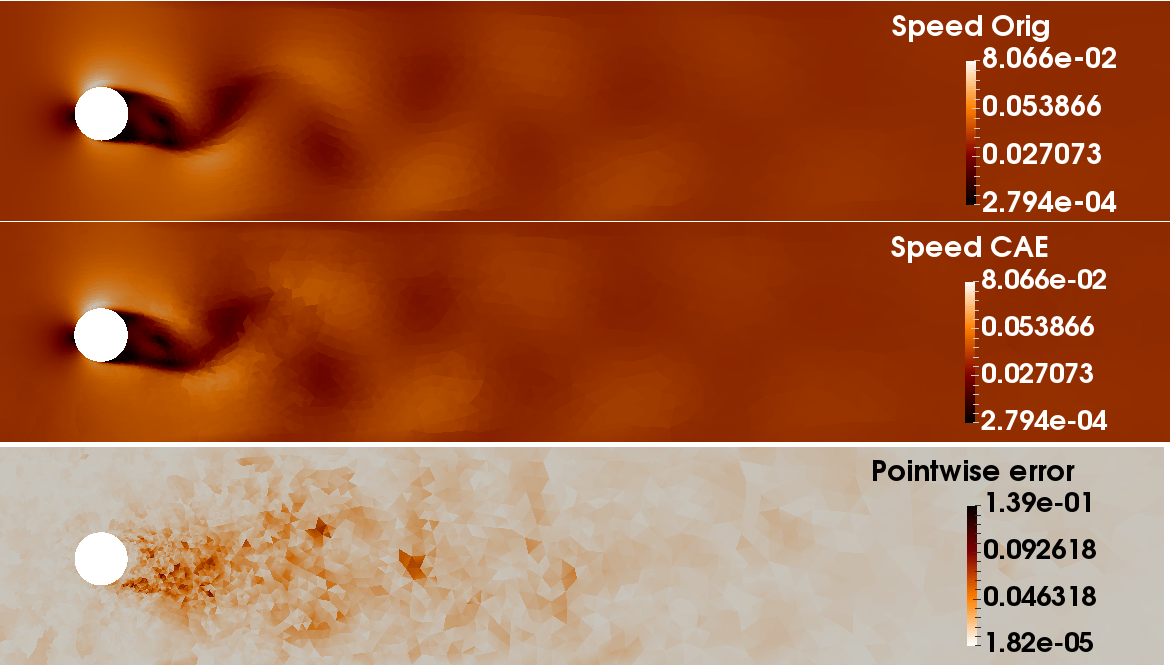}\\[2mm]
\includegraphics[width=8cm]{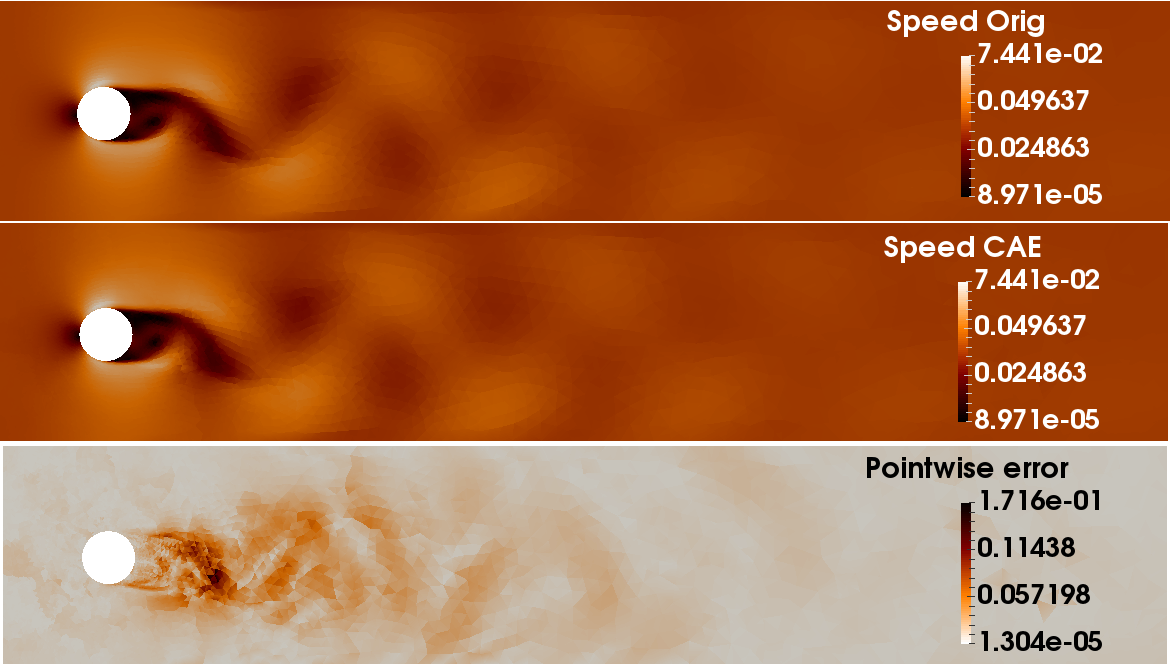}~~%
\includegraphics[width=8cm]{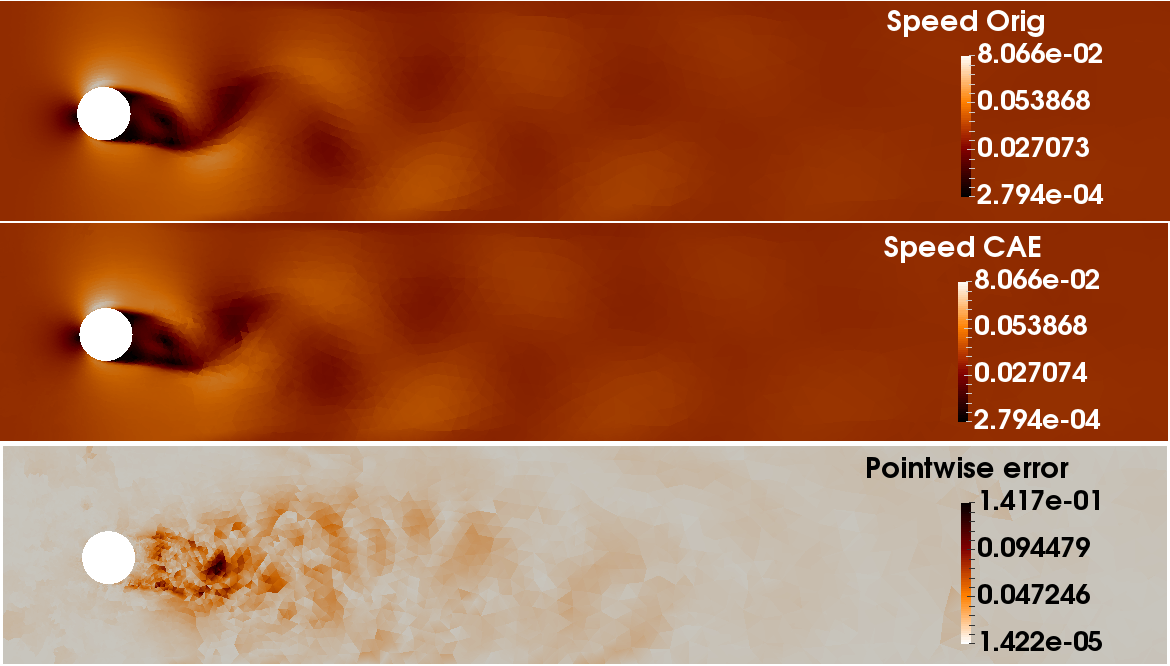}
\caption{Results from the CAE-SFC-NN and CAE-2SFC-NN autoencoders for the data from flow past a cylinder on an unstructured mesh. The data was compressed to 128~variables and the networks were trained for 2000~epochs. For the CAE-SFC-NN network, one example ($t=850$) is shown from the validation data (top left set of three plots), and one ($t=975$) from the test data (top right set of three plots). For the CAE-2SFC-NN network, one example ($t=850$) is shown from the validation data (bottom left), and one ($t=975$) from the test data (bottom right). Each set of plots is divided up into: speed of the original solution~(top); speed after compression and reconstruction by the autoencoder~(middle); pointwise error in speed. }
\label{fig:One-Two128Train-Validate}
\end{figure}

Figure~\ref{fig:OneTwoUnstructured} and Table~\ref{OneTwoUnstructured} show the losses of the two SFC-based autoencoders, one using one space-filling curve and nearest-neighbour smoothing (CAE-SFC-NN), and the other using two space-filling curves and nearest-neighbour smoothing (CAE-2SFC-NN). Judging from the losses, the network that performs slightly better is the autoencoder based on two SFCs, despite the fact that the number of weights in the convolutional layers for this network is approximately half that of the autoencoder with one SFC (CAE-SFC-NN). The reason for its effectiveness is that the two space-filling curves attempt to capture features in all directions across the mesh as the curves sample two different directions (each curve is very roughly orthogonal to the other) at a given node in the mesh. The compression accuracy of the two networks can be seen in Figure~\ref{fig:One-Two128Train-Validate} and are both impressive. On careful inspection of the results based on one SFC, in Figure~\ref{fig:One-Two128Train-Validate} (upper plots), one can see that they are noisier than the results based on two SFCs (lower plots), which highlights a possible advantage of two SFCs leading to more accurate results. However, the network based on one SFC (CAE-SFC-NN) still represents the basic flow features very well.

\begin{table}[htbp]
\centering
%\resizebox{\textwidth}{!}{%
\begin{tabular}{lcc}
\toprule
 & CAE-SFC-NN & CAE-2SFC-NN\\
\midrule
Training data   &${ 2.51 \times 10^{-4} }$    & ${ 2.54 \times 10^{-4}}$\\
Validation data &${ 6.89 \times 10^{-4} }$    & ${ 2.91 \times 10^{-4}}$\\
Test data       &  ${ 6.37 \times 10^{-4} }$  & ${ 2.88 \times 10^{-4}}$\\
\bottomrule
\end{tabular}
%}
\caption{Losses achieved for the convolutional autoencoders based on one and two space-filling curves, both with nearest-neighbour smoothing (CAE-SFC-NN and CAE-2SFC-NN respectively). The data was compressed to 128~variables and the networks were trained for 2000~epochs.  The data set was based on flow past a cylinder solved on an unstructured mesh.}
\label{OneTwoUnstructured}
\end{table}

%\begin{figure}[htbp]
%\centering
%\hbox{ 
%\includegraphics[width=9cm]{figures/One128Train.jpg}
%\hspace{-1.0cm} 
%\includegraphics[width=9cm]{figures/One128Valid.jpg}
%}
%\hbox{ 
%\includegraphics[width=9cm]{figures/Two128Train.jpg}
%\hspace{-1.0cm} 
%\includegraphics[width=9cm]{figures/Two128Valid.jpg}
%}
%\caption{Results from the CAE-SFC-NN and CAE-2SFC-NN autoencoders for the data from flow past a cylinder on an unstructured mesh. The data was compressed to 128~variables and the networks were trained for 2000~epochs. 
%\textcolor{red}{We show the results for CAE-SFC-NN for training data (top left) at t=850 and validation data (top right) at t=950. Similarly we show the results for CAE-2SFC-NN and reproducing the training data (bottom left) at t=750 and validation data (top right) at t=850. The speed of the flow is shown. Each set of figures is divided up into: snapshot from validation data~(top); compressed and reconstructed snapshot from validation data~(middle); point-wise difference between snapshot fluid speed and reconstructed fluid speed. }}
%\label{fig:One-Two128Train-Validate}
%\end{figure}

\subsubsection{Performance of the autoencoder based on two SFCs with nearest-neighbour smoothing for different compression ratios}

\begin{figure}[htbp]
\centering
\includegraphics[width=18cm,trim=20mm 0mm 20mm 0mm, clip]{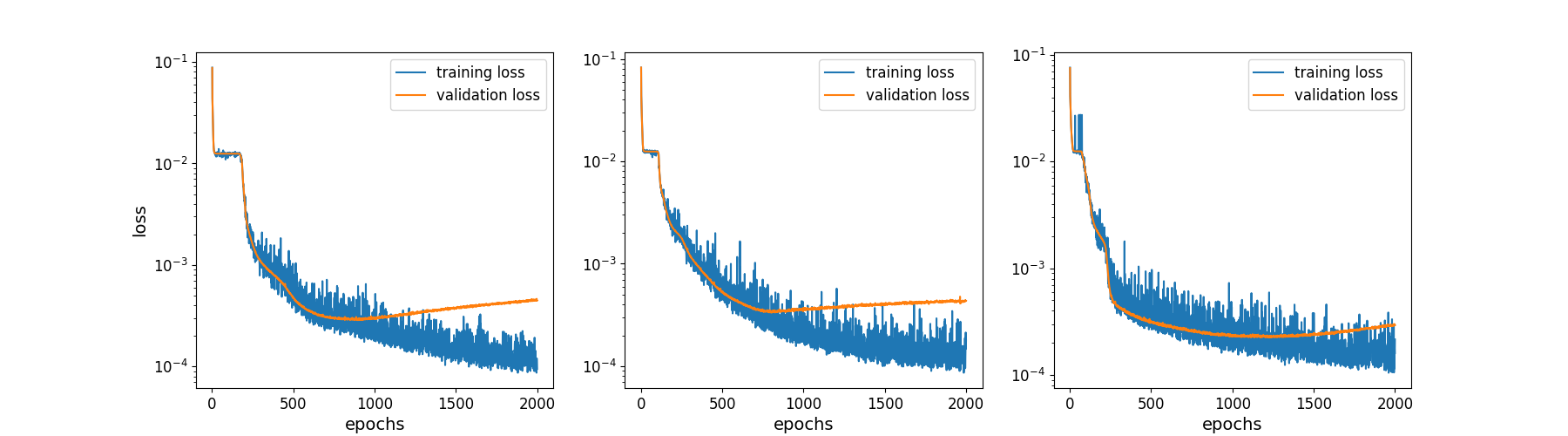}
\caption{Training and validation loss for the CNN based on two space-filling curves with nearest-neighbour smoothing applied to the unstructured mesh data set and compressed to differing numbers of variables: compression size 128 (left); compression size 64 (centre); compression size 32 (right). %\textcolor{red}{Should we include all the losses here? That is, for 16, 8, 4, 2 and 1?}
}
\label{fig:Two1286432}
\end{figure}

\begin{figure}[htbp]
\centering
%\hbox{ 
\includegraphics[width=7.7cm]{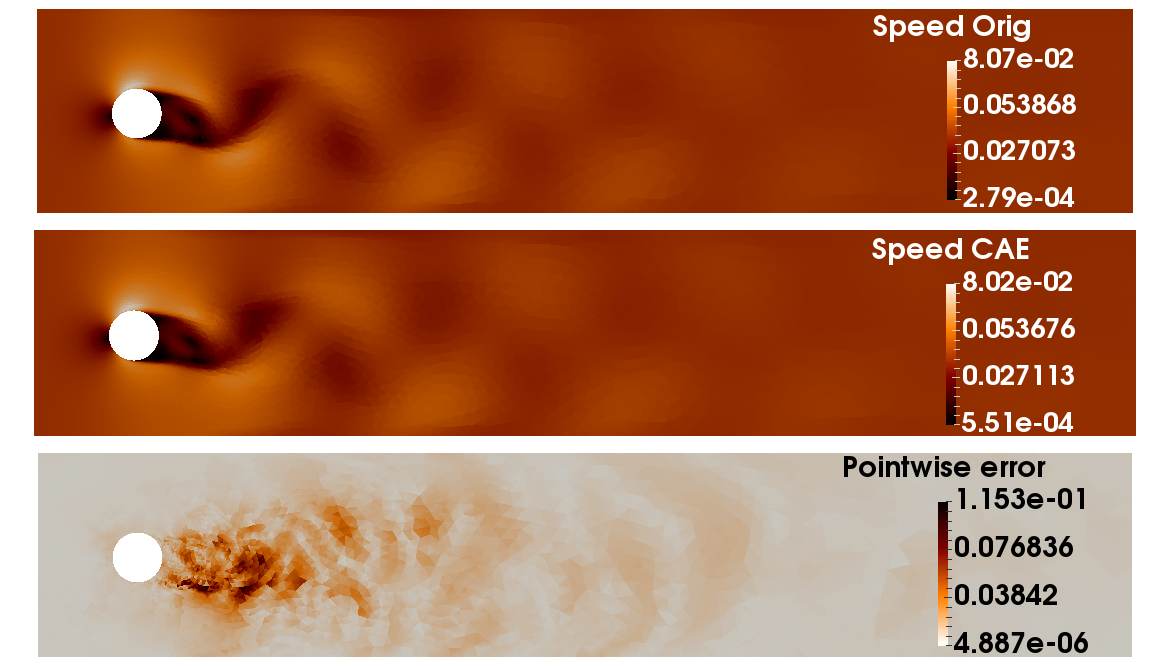}~~~%
\includegraphics[width=8cm]{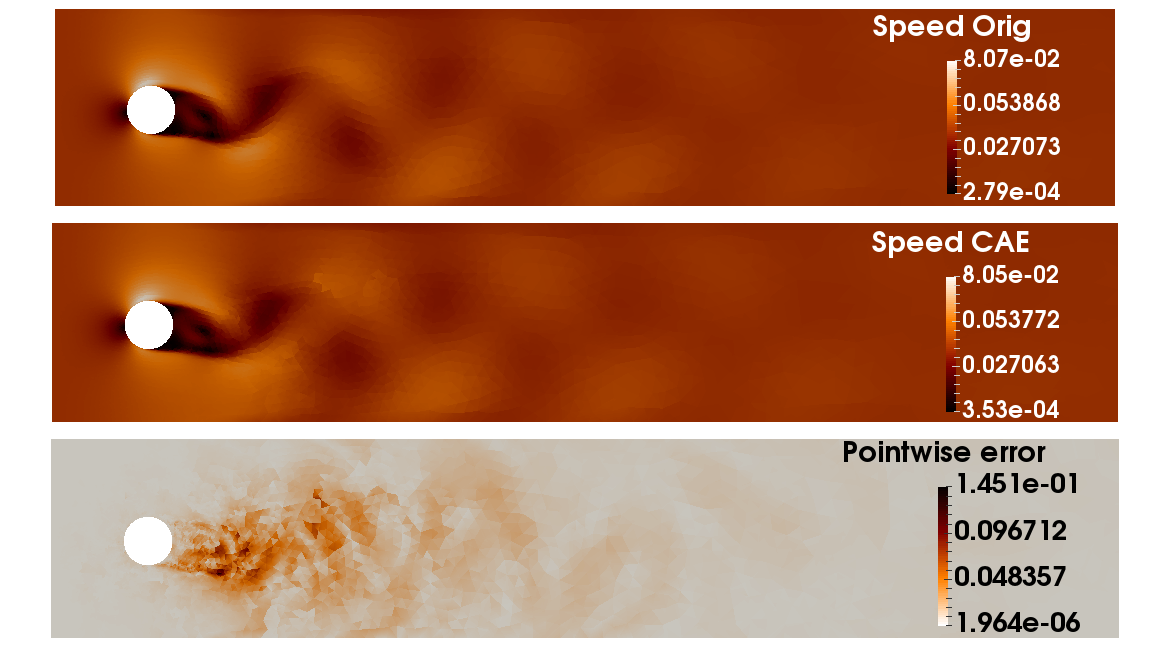}\\[2mm]
%}
%\hbox{ 
\includegraphics[width=8cm]{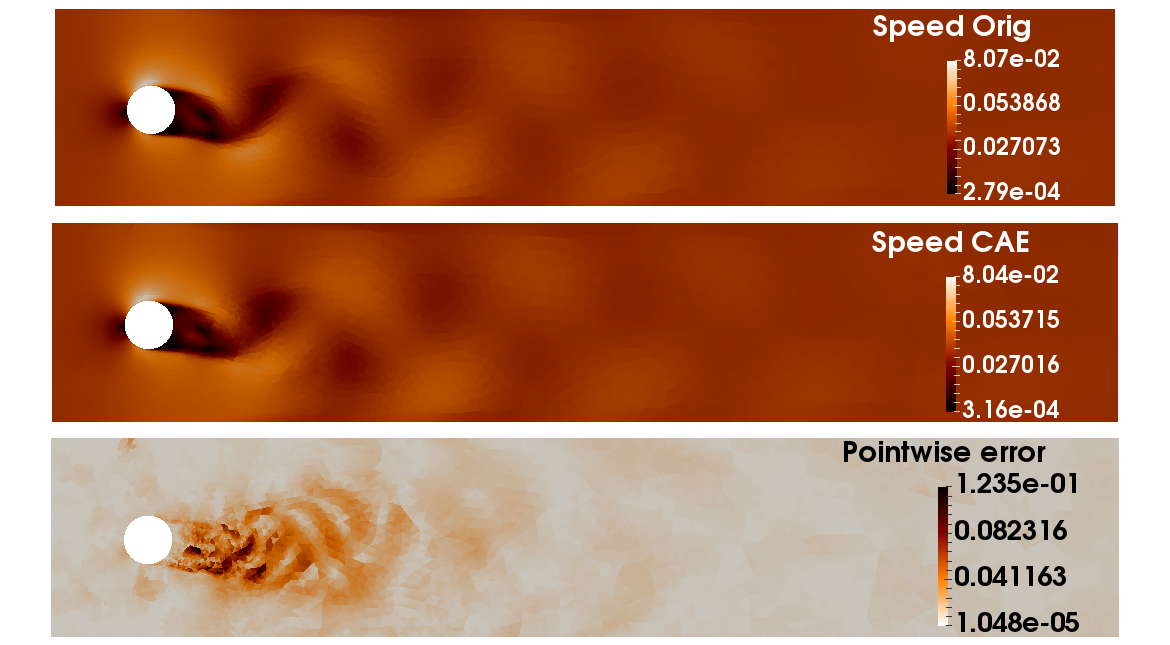}~%
\includegraphics[width=8cm]{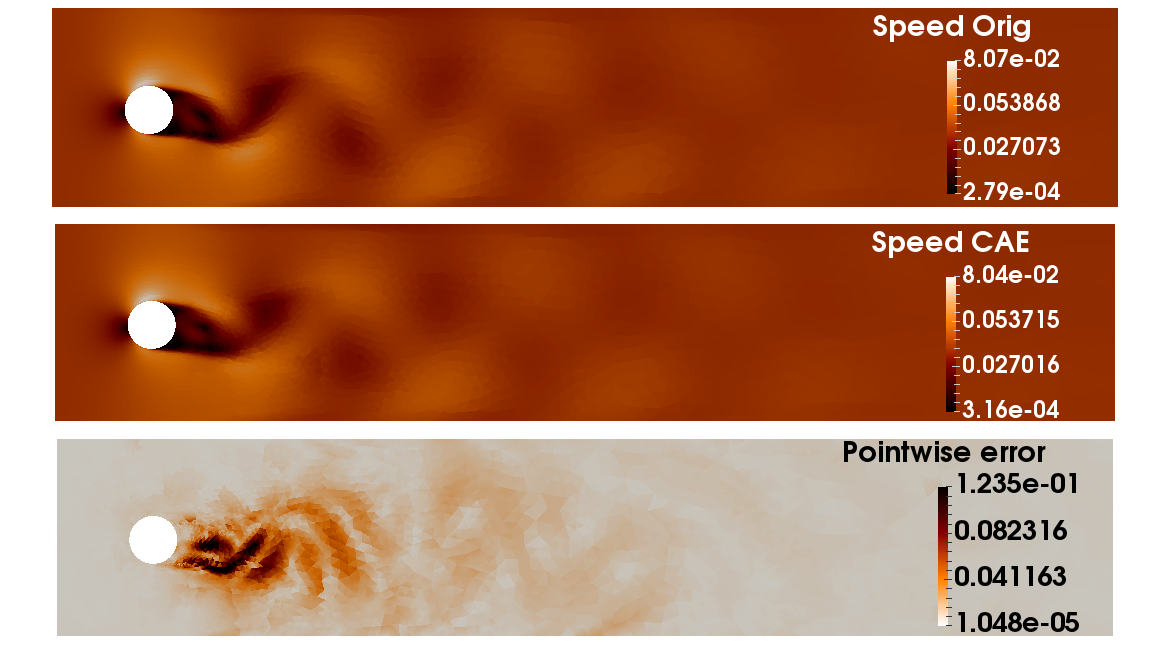}
%}
\caption{Comparison of results of the CAE-SFC2-NN autoencoder when trained four times to compress to four different numbers of latent variables: 64~variables (top left set of three plots), 16~variables (top right set), 4~variables (bottom left set) and 1 variable (bottom right set). The example shown corresponds to time level~975 which lies in the test data. Each set of plots is divided up into: the original solution~(top); solution compressed and reconstructed by the autoencoder~(middle); pointwise error in speed~(bottom). }
\label{fig:comparison-sizes-validate}
\end{figure}

For the convolutional autoencoder based on two space-filling curves with nearest-neighbour smoothing (CAE-2SFC-NN), the effect of dimension of latent space is explored. In the previous section, a dimension of~128 was used. Here, the autoencoder is trained multiple times for~64, 32, 16, 8, 4, 2 and 1~latent variables. After 2000~epochs (4000 in the case of 2 and 1 variables), the losses are shown in Figure~\ref{fig:Two1286432} and Table~\ref{TwoUnstructuredLossTable}, which indicate that, even though we reduce the number of compressed variables, the accuracy is still reasonable. In fact, the spatial distribution of the results for compression to 1~variable looks similar to the spatial distribution for 128~compressed variables in Figure~\ref{fig:One-Two128Train-Validate}. Notice that as the number of compressed variables is reduced the gap between the validation and training loss is reduced. This could be interpreted as a reduction in the tendency for over-fitting as the number of latent or compressed variables is reduced. Also notice from Figures~\ref{fig:Two1286432} and~\ref{fig:comparison-sizes-validate} that the greater the number of compressed variables the shorter the wavelength of the structures in the error fields. Thus, with smaller compression sizes one seems to be able to capture the basic structures and as more variables are introduced finer scale structures are captured.  Note that the imprints of the space-filling curves are not noticeable 
for the smaller number of compressed variables, therefore, the nearest-neighbour smoothing layers of the autoencoder must be performing well. However, one can see this imprint in the error field for larger numbers of compressed variables. 

%\textcolor{red}{We have interpolated the flow past a cylinder results onto a structured mesh and applied a classical CNN but were unable to compress to less than 32 variables with the resulting network which provides evidence that the SFC-CNN can produce more robust results than the classical approach although this needs further investigation. }
Comparing the CAE-2SFC-NN network with the SVD (see Table~\ref{TwoUnstructuredLossTable}), we find that for compressed variable sizes of 128, 64, 32, 16, and~8, the truncation error of the SVD is smaller than the loss of the autoencoders.  However, for compressed variables sizes of 4, 2 and~1, the losses of the autoencoders are smaller than those of the SVD. This is because flow past a cylinder is a smooth problem, periodic in time, for which only a few SVD modes are needed to accurately represent the flow. However, this is still a hard problem for the an autoencoder, and the CAE-2SFC-NN network does manage to more effectively compress to a small number of variables than the SVD approach.  

%\begin{figure}[htbp]
%\centering
%\includegraphics[width=15cm]{figures/Two1286432.jpg}
%\caption{Training and validation loss for the CNN based on two space-filling curves with nearest-neighbour smoothing applied to the unstructured mesh data set and compressed to differing numbers of variables: compression size 128 (left); compression size 64 (centre); compression size 32 (right). %\textcolor{red}{Should we include all the losses here? That is, for 16, 8, 4, 2 and 1?}
%}
%\label{fig:Two1286432}
%\end{figure}

\begin{table}[htbp]
\centering
\resizebox{\textwidth}{!}{%
\begin{tabular}{lcccccccc}
\toprule
 & \multicolumn{8}{c}{number of compressed variables} \\
  \cmidrule{2-9}
 &  128 &  64 &  32&  16&  8&  4&  2&  1\\
\midrule
Training data &${ 3.08 \times 10^{-4}}$&${ 2.54 \times 10^{-4} }$&${ 1.65 \times 10^{-4} }$&${ 4.11 \times 10^{-4} }$&${ 4.37 \times 10^{-4} }$&${ 4.10 \times 10^{-4} }$&${ 1.79 \times 10^{-4} }$&${ 5.99 \times 10^{-4} }$\\

Validation data &${ 2.91 \times 10^{-4}}$&${ 3.33 \times 10^{-4} }$ &${ 2.29 \times 10^{-4} }$&${ 3.32 \times 10^{-4} }$&${ 2.21 \times 10^{-4} }$&${ 2.67 \times 10^{-4} }$&${ 1.96 \times 10^{-4} }$&${ 6.00 \times 10^{-4} }$\\

Test data &${ 2.88 \times 10^{-4}}$&${ 3.32 \times 10^{-4} }$ &${ 2.28 \times 10^{-4} }$&${ 3.63 \times 10^{-4} }$&${ 2.21 \times 10^{-4} }$&${ 2.63 \times 10^{-4} }$&${ 1.98 \times 10^{-4} }$&${ 5.35 \times 10^{-4} }$\\
\midrule
SVD &${ 1.40 \times 10^{-8}}$&${ 3.68 \times 10^{-7} }$ &${ 4.03 \times 10^{-6} }$&${ 2.97 \times 10^{-5} }$&${ 1.26 \times 10^{-4} }$&${ 4.44 \times 10^{-4} }$&${ 5.99 \times 10^{-3} }$&${ 1.23 \times 10^{-2} }$\\
\bottomrule
\end{tabular}
}
\caption{Losses of the CAE-2SFC-NN autoencoder, trained to compress to the following numbers of variables: 128, 64, 32, 16, 8, 4, 2 and~1. Truncation error of the SVD when compressing to the same numbers of variables.}
\label{TwoUnstructuredLossTable}
\end{table}

%\begin{figure}[htbp]
%\centering
%\hbox{ 
%\includegraphics[width=9cm]{figures/Two1Valid.jpg}
%\hspace{-1.0cm} 
%\includegraphics[width=9cm]{figures/Two4Valid.jpg}
%}
%\hbox{ 
%\includegraphics[width=9cm]{figures/Two16Valid.jpg}
%\hspace{-1.0cm} 
%\includegraphics[width=9cm]{figures/Two64Valid.jpg}
%}
%\caption{\textcolor{red}{Comparison of results of the SFC-CNN with two-space filling curves and nearest neurbour smoothing applied to randomly selected validation data for different compression sizes: 1 variable (top left) at t=860; 4 variables (top right) at t=846; 16 variables (bottom left) at t=800; 64 variables (bottom right) at t=805. }Each set of figures is divided up into: snapshot from validation data~(top); compressed and reconstructed snapshot from validation data~(middle); difference  between snapshot fluid speed and reconstructed fluid speed~(bottom). }
%\label{fig:comparison-sizes-validate}
%\end{figure}

\begin{figure}[htbp]
\centering
\includegraphics[width=\textwidth]{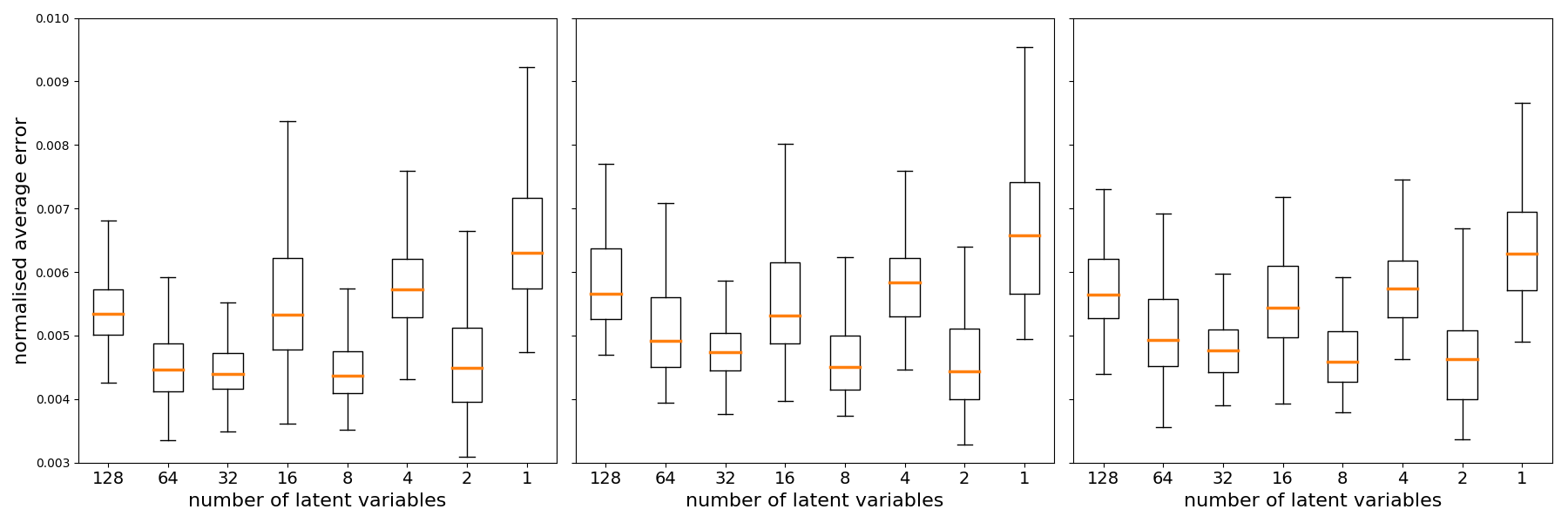}
\caption{Box and whisker plots of the normalised average error of the training data (left), validation data (centre) and the test data (right), for the CAE-2SFC-NN autoencoder.}
\label{fig:box_whisker}
\end{figure}

Finally, a box and whisker plot can be seen in Figure~\ref{fig:box_whisker}. Here, the average absolute error in space is calculated for each example in the data sets for each compression ratio. Then, for each compression ratio and for each data set (training, validation, test), the median error is shown in orange, the upper~(Q3) and lower~(Q1) quartiles are represented by the box and the whiskers indicate the range of the data, taken as the data points closest to Q1-1.5IQR and Q3+1.5IQR, where the inter-quartile range (IQR) is Q3-Q1. Approximately 1\% of the data points are considered to be outliers and are not shown. The variation over training (left), validation (centre) and test (right) data is similar due to the time-periodic problem studied here. Remarkably, as the number of latent variables is reduced, the range of errors becomes slightly larger, but there is no obvious upward trend. This agrees with the plots shown in  Figure~\ref{fig:comparison-sizes-validate}, where the accuracy seems not to be affected by the number of latent variables. 

%\clearpage
\section{Future work} 
%The success of the transformation of multi-dimensional data by space-filling curves in advance of the application of CNNs has been established in this paper.
The SFC-based approach to applying CNNs to data on unstructured meshes has been established in this paper. Furthermore, this approach has been applied to data on an unstructured mesh with differing areas of resolution. A next big challenge will be to apply the SFC-based CNN approach to solutions on adaptive meshes, that is meshes that change their topology and resolution as time evolves in order to optimally represent the physics. The SFC-based approach has the potential to be able to deal with even this, as the filters do not necessarily depended on the mesh. One may, for example, be able to take into account the connectivity by also including the coordinates as inputs to the 1D CNNs or to the SFC-based CNN. Once these issues are dealt with then it may be possible to apply (at least as a starting model for training) the weights from one example problem to the weights of a different fluids problem. In addition, the approach may be used to interpolate from one mesh to a different mesh while increasing the numerical resolution and sharpness of the features especially when interpolating onto meshes or mesh parts that have finer resolution than the original mesh. This interpolation approach is very much akin to the deblurring GAN-CNN or DeblurGAN~\cite{Kupyn2018}.  

Further work will also be needed to explore the application of the method to 3D, with, potentially, the use of three space-filling curves. 

Another area of further study is the inclusion of the time dimension within the CNN.  Rather than treating the time dimension independently~\cite{Wu_P_2020}, a space-filling curve could be applied simultaneously to four-dimensional solutions (three spatial dimensions and one temporal dimension).% This could be achieved by putting 1D space-filling curves through the four-dimensional space. However, it may also be beneficial to treat the time dimension independently, for example, if one has several or more time levels to compress simultaneously~\cite{Wu_P_2020}. 

%\textcolor{blue}{\textit{Move to Future work? One challenge remains for CNNs when compressing data from incompressible flow problems however. For SVD-based methods, the resulting basis functions span the snapshot velocities and, therefore, the incompressibility constraint will be satisfied providing the snapshot velocities satisfy incompressibility. }}

\section{Conclusions} The space-filling curve (SFC) approach shows great promise in  enabling the application of convolutional networks to data from unstructured meshes. The approach has several features that make it ideal for application to convolutional neural networks (CNNs). This includes the aptitude of SFCs for automatic coarsening of meshes. %These features have enabled the SFC-based CNN to show great promise in being able to compress information from unstructured meshes containing the numerical solution of differential equations.  % Space-Filling Curve CNN (SFC-CNN) has several features that make it ideal for  neural networks. This includes automatic coarsening of the meshes during the one dimensional CNN filter operations. 
%The Space Filling Curve Convolutional Neural Network 

%In fact, it is hoped that, it might be used whenever machine learning methods are used in combination with data from unstructured mesh finite element or control volume solution methods. 
We demonstrate the approach by compressing the results of solutions of the Navier-Stokes equations for incompressible flow past a cylinder and show that it is able to compress of the order of 50,000 solution variables with a complex discontinuous Galerkin stencil down to between 128~variables and 1~variable, while maintaining accuracy. We also show that it can be beneficial to use two space-filling curves to help increase the accuracy of the CNN. We also found it important to reduce the noise in the outputs of the SFC-based autoencoders by introducing sparse smoothing layers near the output (and input) of the autoencoders. On structured mesh data sets the accuracy of the new SFC-based autoencoders are similar to classical autoencoders and it seems worth exploring their comparative capabilities even for data on structured grids. 
%Most importantly this is the first time that space filling curves have been used 
%to construct a 1D vector that is able to be used with a classical (standard approach in the literature) 
%CNN methods and thus it 
%maintains the multi-scale feature capturing abilities of classical CNN approaches. 
%We hope this new SFC-CNN will become important 
%for unstructured mesh modelling combined with machine learning  
%in the near future.  

%\textcolor{blue}{\textit{Use this text in the conclusions?} The Space-Filling Curve CNN (SFC-CNN) has several features that make it ideal for  neural networks. This includes automatic coarsening of the meshes during the one dimensional CNN filter operations. Although space filling curves have been 
%used previously with CNNs this is the first time they have been applied to unstructured mesh finite element solutions and good accuracy of the method is achieved partly through the use of multiple space filling curves and partly through the use of sparse smoothing layers near the output of the SFC-CNNs.  }

%%%%%%%%%%%%%%%%%%%%%%%%%%%%%%%%%%%%%%%%%%%%%%%%%
\section*{Acknowledgements}
The authors would like to acknowledge the following EPSRC grants: MUFFINS (EP/P033180/1); PREMIERE (EP/T000414/1); MAGIC (EP/N010221/1) and INHALE (EP/T003189/1). 
%%%%%%%%%%%%%%%%%%%%%%%%%%%%%%%%%%%%%%%%%%%%%%%%%  

%\section{References}
\bibliographystyle{elsarticle-num}
\bibliography{references}

\begin{thebibliography}{10}
\expandafter\ifx\csname url\endcsname\relax
  \def\url#1{\texttt{#1}}\fi
\expandafter\ifx\csname urlprefix\endcsname\relax\def\urlprefix{URL }\fi
\expandafter\ifx\csname href\endcsname\relax
  \def\href#1#2{#2} \def\path#1{#1}\fi

\bibitem{Minnen_2018}
D.~Minnen, J.~Ball{\'e}, G.~Toderici, {High-Fidelity Generative Image
  Compression}, arXiV: 1809.02736 (2018).

\bibitem{Mentzer_2020}
F.~Mentzer, G.~Toderici, M.~Tschannen, E.~Agustsson, {High-Fidelity Generative
  Image Compression}, arXiV: 2006.09965v2 (2020).

\bibitem{Krizhevsky_2012}
A.~Krizhevsky, I.~Sutskever, G.~E. Hinton, {ImageNet Classification with Deep
  Convolutional Neural Networks}, in: F.~Pereira, C.~J.~C. Burges, L.~Bottou,
  K.~Q. Weinberger (Eds.), Advances in Neural Information Processing Systems
  25, Curran Associates, Inc., 2012, pp. 1097--1105.

\bibitem{Yim2015}
J.~Yim, J.~Ju, H.~Jung, J.~Kim, Image classification using convolutional neural
  networks with multi-stage feature, in: J.-H. Kim, W.~Yang, J.~Jo, P.~Sincak,
  H.~Myung (Eds.), Robot Intelligence Technology and Applications 3, Springer
  International Publishing, Cham, 2015, pp. 587--594.

\bibitem{Gonzalez_2018}
F.~J. Gonzalez, M.~Balajewicz, {Deep convolutional recurrent autoencoders for
  learning low-dimensional feature dynamics of fluid systems}, arXiV:
  1808.01346 (2018).

\bibitem{lee2020model}
K.~Lee, K.~T. Carlberg, Model reduction of dynamical systems on nonlinear
  manifolds using deep convolutional autoencoders, Journal of Computational
  Physics 404 (2020) 108973.

\bibitem{Walton_2017}
S.~Walton, O.~Hassan, K.~Morgan, Advances in co-volume mesh generation and mesh
  optimisation techniques, Computers \& Structures 181 (2017) 70--88.

\bibitem{Xie_2016}
Z.~Xie, D.~Pavlidis, P.~Salinas, J.~R. Percival, C.~C. Pain, O.~K. Matar, A
  balanced-force control volume finite element method for interfacial flows
  with surface tension using adaptive anisotropic unstructured meshes,
  Computers \& Fluids 138 (2016) 38--50.

\bibitem{Kampitsis_2020}
A.~E. Kampitsis, A.~Adam, P.~Salinas, C.~C. Pain, A.~H. Muggeridge, M.~D.
  Jackson, Dynamic adaptive mesh optimisation for immiscible viscous fingering,
  Computational Geosciences 24 (2020) 1221--1237.

\bibitem{Alauzet_2016}
F.~Alauzet, A.~Loseille, {A Decade of Progress on Anisotropic Mesh Adaptation
  for Computational Fluid Dynamics}, Computer-Aided Design 72 (2016) 13--39.

\bibitem{Mack2020}
J.~Mack, R.~Arcucci, M.~Molina-Solana, Y.-K. Guo, {Attention-based
  Convolutional Autoencoders for 3D-Variational Data Assimilation}, Computer
  Methods in Applied Mechanics and Engineering 372 (2020) 113291.

\bibitem{Carlberg2020}
K.~Lee, K.~T. Carlberg, Model reduction of dynamical systems on nonlinear
  manifolds using deep convolutional autoencoders, Journal of Computational
  Physics 404 (2020) 108973.

\bibitem{Xu2020}
J.~Xu, K.~Duraisamy, Multi-level convolutional autoencoder networks for
  parametric prediction of spatio-temporal dynamics, Computer Methods in
  Applied Mechanics and Engineering 372 (2020) 113379.

\bibitem{Wu2020}
Z.~Wu, S.~Pan, F.~Chen, G.~Long, C.~Zhang, P.~Yu, {A Comprehensive Survey on
  Graph Neural Networks}, IEEE Transactions on Neural Networks and Learning
  Systems (2020) 1--21.

\bibitem{kipf2017}
T.~N. Kipf, M.~Welling, {Semi-Supervised Classification with Graph
  Convolutional Networks}, arXiV: 1609.02907 (2017).

\bibitem{PyTorch_2019}
A.~Paszke, S.~Gross, F.~Massa, A.~Lerer, J.~Bradbury, G.~Chanan, T.~Killeen,
  Z.~Lin, N.~Gimelshein, L.~Antiga, A.~Desmaison, A.~Kopf, E.~Yang, Z.~DeVito,
  M.~Raison, A.~Tejani, S.~Chilamkurthy, B.~Steiner, L.~Fang, J.~Bai,
  S.~Chintala, {PyTorch: An Imperative Style, High-Performance Deep Learning
  Library}, in: H.~Wallach, H.~Larochelle, A.~Beygelzimer, F.~d\textquotesingle
  Alch\'{e}-Buc, E.~Fox, R.~Garnett (Eds.), Advances in Neural Information
  Processing Systems 32, Curran Associates, Inc., 2019, pp. 8024--8035.

\bibitem{Hanocka2019}
R.~Hanocka, A.~Hertz, N.~Fish, R.~Giryes, S.~Fleishman, D.~Cohen-Or, {MeshCNN:
  A Network with an Edge}, ACM Trans. Graph. 38~(4) (2019).

\bibitem{Dargaville_2020}
S.~Dargaville, A.~G. Buchan, R.~P. Smedley-Stevenson, P.~N. Smith, C.~C. Pain,
  A comparison of element agglomeration algorithms for unstructured geometric
  multigrid, arXiV: 2005.09104 (2020).

\bibitem{LU202042}
Q.~Lu, C.~Chen, W.~Xie, Y.~Luo, {PointNGCNN:} {D}eep convolutional networks on
  {3D} point clouds with neighborhood graph filters, Computers \& Graphics 86
  (2020) 42--51.

\bibitem{Pain_2001}
C.~C. Pain, A.~P. Umpleby, C.~R.~E. de~Oliveira, A.~J.~H. Goddard, Tetrahedral
  mesh optimisation and adaptivity for steady-state and transient finite
  element calculations, Computer Methods in Applied Mechanics and Engineering
  190~(29) (2001) 3771--3796.

\bibitem{Sanchez-Gonzalez2020}
A.~Sanchez-Gonzalez, J.~Godwin, T.~Pfaff, R.~Ying, J.~Leskovec, P.~W.
  Battaglia, {Learning to Simulate Complex Physics with Graph Networks}, arXiV:
  2002.09405v2 (2020).

\bibitem{Pfaff2020}
T.~Pfaff, M.~Fortunato, A.~Sanchez-Gonzalez, P.~W. Battaglia, {Learning
  Mesh-Based Simulation with Graph Networks}, arXiV: 2010.03409 (2020).

\bibitem{Peano_1890}
G.~Peano, Sur une courbe, qui remplit toute une aire plane, Mathematische
  Annalen (in French) 36 (1890) 157--160.

\bibitem{Hilbert_1891}
D.~Hilbert, Ueber die stetige abbildung einer line auf ein
  fl{\"a}chenst{\"u}ck, Mathematische Annalen 38 (1891) 459--460.

\bibitem{Sagan}
H.~Sagan, {Space-Filling Curves}, Springer, 1994.

\bibitem{Bader}
M.~Bader, {Space-Filling Curves: An Introduction with Applications in
  Scientific Computing}, Springer, 2013.

\bibitem{Bohm2018}
C.~{B{\"o}hm}, M.~{Perdacher}, C.~{Plant}, {A Novel Hilbert Curve for
  Cache-locality Preserving Loops}, IEEE Transactions on Big Data (2018) 1--15.

\bibitem{Bungartz_2006}
H.-J. Bungartz, M.~Mehl, T.~Weinzierl, A parallel adaptive cartesian pde solver
  using space–filling curves, in: W.~Nagel, W.~Walter, W.~Lehner (Eds.),
  Euro-Par 2006 Parallel Processing, Springer Berlin Heidelberg, Berlin,
  Heidelberg, 2006, pp. 1064--1074.

\bibitem{Behrens_2000}
J.~Behrens, J.~Zimmermann, Parallelizing an unstructured grid generator with a
  space-filling curve approach, in: A.~Bode, T.~Ludwig, W.~Karl,
  R.~Wism{\"u}ller (Eds.), Euro-Par 2000 Parallel Processing, Springer Berlin
  Heidelberg, Berlin, Heidelberg, 2000, pp. 815--823.

\bibitem{Sprecher_2002}
D.~A. Sprecher, S.~Draghici, {Space-filling curves and Kolmogorov
  superposition-based neural networks}, Neural Networks 15 (2002) 57--67.

\bibitem{Papamarkos_1999}
N.~Papamarkos, {Color Reduction Using Local Features and a Kohonen
  Self-Organized Feature Map Neural Network}, International Journal of Imaging
  Systems and Technology 10 (1999) 404--409.

\bibitem{Anjum_2019}
M.~M. Anjum, I.~A. Tahmid, M.~S. Rahman, {CNN Model With Hilbert Curve
  Representation of DNA Sequence For Enhancer Prediction}, bioRxiv:
  10.1101/552141 (2019).

\bibitem{Yin_2018}
B.~Yin, M.~Balvert, D.~Zambrano, A.~Sch{\"o}nhuth, S.~M. Bohte, {An Image
  Representation Based Convolutional Network for DNA Classification}, 6th
  International Conference on Learning Representations (2018).

\bibitem{Skubalska_1997}
E.~Skubalska-Rafaj{\l}owicz, Applications of the space-filling curves with data
  driven measure-preserving property, Nonlinear Analysis: Theory, Methods \&
  Applications 30~(3) (1997) 1305--1310.

\bibitem{Corcoran_2018}
T.~Corcoran, R.~Zamora-Resendiz, X.~Liu, S.~Crivelli, {A Spatial Mapping
  Algorithm with Applications in Deep Learning-Based Structure Classification},
  arXiv: 1802.02532v2 (2018).

\bibitem{Phillips_2020}
T.~Phillips, C.~E. Heaney, P.~N. Smith, C.~C. Pain, {An autoencoder-based
  reduced-order model for eigenvalue problems with application to neutron
  diffusion problems}, arXiv:2008.10532 [math.NA],$\ $ (2020).

\bibitem{Hinton_2011}
G.~E. Hinton, A.~Krizhevsky, S.~Wang, {Transforming Auto-encoders}, in:
  T.~Honkela, W.~Duch, M.~Girolami, S.~Kaski (Eds.), Artificial Neural Networks
  and Machine Learning -- ICANN 2011, Springer, 2011, pp. 44--51.

\bibitem{Hesthaven_2018}
J.~Hesthaven, S.~Ubbiali, Non-intrusive reduced order modeling of nonlinear
  problems using neural networks, Journal of Computational Physics 363 (2018)
  55--78.

\bibitem{Xiao_BE_2019}
D.~Xiao, C.~E. Heaney, L.~Mottet, F.~Fang, W.~Lin, I.~M. Navon, Y.~Guo, O.~K.
  Matar, A.~G. Robins, C.~C. Pain, A reduced order model for turbulent flows in
  the urban environment using machine learning, Building and Environment 148
  (2019) 323--337.

\bibitem{Xiao_CAF_2019}
D.~Xiao, C.~Heaney, F.~Fang, L.~Mottet, R.~Hu, D.~Bistrian, E.~Aristodemou,
  I.~Navon, C.~Pain, A domain decomposition non-intrusive reduced order model
  for turbulent flows, Computers \& Fluids 182 (2019) 15--27.

\bibitem{Ahmed2019}
S.~E. Ahmed, S.~M. Rahman, O.~San, A.~Rasheed, I.~M. Navon, Memory embedded
  non-intrusive reduced order modeling of non-ergodic flows, Physics of Fluids
  31~(12) (2019) 126602.

\bibitem{Buluc2016}
A.~Bulu{\c{c}}, H.~Meyerhenke, I.~Safro, P.~Sanders, C.~Schulz, {Recent
  Advances in Graph Partitioning}, Springer International Publishing, 2016, pp.
  117--158.

\bibitem{Fitzgerald2019}
A.~P. Fitzgerald, B.~Kochunas, S.~Stimpson, T.~Downar, Spatial decomposition of
  structured grids for nuclear reactor simulations, Annals of Nuclear Energy
  132 (2019) 686--701.

\bibitem{Kawamoto2019}
T.~Kawamoto, M.~Tsubaki, T.~Obuchi, Mean-field theory of graph neural networks
  in graph partitioning, Journal of Statistical Mechanics: Theory and
  Experiment 12 (2019) 124007.

\bibitem{Karypis1998_JPDC}
G.~Karypis, V.~Kumar, Multilevel $k$-way {P}artitioning {S}cheme for
  {I}rregular {G}raphs, J. Parallel Distrib. Comput. 48~(1) (1998) 96--129.

\bibitem{Karypis1998_SIAM}
G.~Karypis, V.~Kumar, {A Fast and High Quality Multilevel Scheme for
  Partitioning Irregular Graphs}, SIAM J. Sci. Comput. 20~(1) (1998) 359--392.

\bibitem{Moulitsas2001}
I.~{Moulitsas}, G.~{Karypis}, Multilevel algorithms for generating coarse grids
  for multigrid methods, in: SC '01: Proceedings of the 2001 ACM/IEEE
  Conference on Supercomputing, 2001, pp. 15--15.

\bibitem{Pain_1999}
C.~C. Pain, C.~R.~E. De~Oliveira, A.~J.~H. Goddard, A neural network graph
  partitioning procedure for grid-based domain decomposition, International
  Journal for Numerical Methods in Engineering 44~(5) (1999) 593--613.

\bibitem{Kupyn2018}
O.~{Kupyn}, V.~{Budzan}, M.~{Mykhailych}, D.~{Mishkin}, J.~{Matas}, {DeblurGAN:
  Blind Motion Deblurring Using Conditional Adversarial Networks}, in: 2018
  IEEE/CVF Conference on Computer Vision and Pattern Recognition, 2018, pp.
  8183--8192.

\bibitem{Wu_P_2020}
P.~Wu, J.~Sun, X.~Chang, W.~Zhang, R.~Arcucci, Y.~Guo, C.~C. Pain, Data-driven
  reduced order model with temporal convolutional neural network, Computer
  Methods in Applied Mechanics and Engineering 360 (2020) 360--373.

\end{thebibliography}

\end{document}